\documentclass[reqno,a4paper,11pt,oneside,final,draft]{amsart}
\usepackage{fixltx2e}
\usepackage{float}
\usepackage{wrapfig}
\usepackage{amssymb,amsmath,amsthm}
\usepackage{mathtools}
\usepackage{fixmath}
\usepackage{braket}
\usepackage{caption}
\usepackage{caption}
\usepackage{subcaption}
\usepackage{enumitem}
\usepackage{units}
\usepackage[final]{graphicx}
\usepackage{xcolor}
\usepackage{booktabs}
\usepackage[english]{babel}
\usepackage[T1]{fontenc}
\usepackage{lmodern}
\usepackage{algorithm}
\usepackage{algorithmic}
\usepackage{amssymb,amsmath,exscale}
\usepackage[final,stretch=10]{microtype}

\usepackage[obeyDraft,textsize=scriptsize,textwidth=2.7cm]{todonotes}
\usepackage{showlabels}

\DeclareMathOperator{\supp}{supp}

\DeclareMathOperator*{\argmin}{arg\,min}

\newcommand{\Lzwoboch}{L^2_\mu (Y_0;W_\infty)}
\newcommand{\Lzwobochm}{L^2_\mu (Y_0;L^2(I;\R^m))}
\newcommand{\Lzwobochn}{L^2_\mu (Y_0;L^2(I;\R^n))}
\newcommand{\Linfboch}{L^\infty_\mu (Y_0;W_\infty)}
\newcommand{\Lzwobochnorm}[1]{\|#1\|_{L^2_\mu (Y_0;W_\infty)}}
\newcommand{\Lzwobochnormm}[1]{\|#1\|_{L^2_\mu (Y_0;L^2(I;\R^m))}}
\newcommand{\Lzwobochnormn}[1]{\|#1\|_{L^2_\mu (Y_0;L^2(I;\R^n))}}
\newcommand{\Linfbochnorm}[1]{\|#1\|_{L^\infty_\mu (Y_0;W_\infty)}}
\newcommand{\abd}{\mathbf{a}}
\newcommand{\Cc}{\mathcal{C}}
\newcommand{\F}{\mathcal{F}}
\newcommand{\N}{\mathbb{N}}
\newcommand{\R}{\mathbb{R}}
\newcommand{\de}{\mathrm{d}}
\newcommand{\scalarrn}[2]{\left(#1,#2\right)_{\R^n}}

\newcommand{\wnorm}[1]{\|#1\|_{W_\infty}}

\newcommand{\tr}[1]{tr(C(1)^{-1})}

\newcommand{\scalarl}[2]{(#1,#2)_{L^2(I;\R^n)}}

\newcommand{\eps}{\varepsilon}

\newcommand{\vertiii}[1]{{\left\vert\kern-0.25ex\left\vert\kern-0.25ex\left\vert #1
    \right\vert\kern-0.25ex\right\vert\kern-0.25ex\right\vert}}

\newcommand{\Ya}{\mathcal{Y}_{ad}}

\makeatletter
\newcommand*\rel@kern[1]{\kern#1\dimexpr\macc@kerna}
\newcommand*\widebar[1]{%
  \begingroup
  \def\mathaccent##1##2{%
    \rel@kern{0.8}%
    \overline{\rel@kern{-0.8}\macc@nucleus\rel@kern{0.2}}%
    \rel@kern{-0.2}%
  }%
  \macc@depth\@ne
  \let\math@bgroup\@empty \let\math@egroup\macc@set@skewchar
  \mathsurround\z@ \frozen@everymath{\mathgroup\macc@group\relax}%
  \macc@set@skewchar\relax
  \let\mathaccentV\macc@nested@a
  \macc@nested@a\relax111{#1}%
  \endgroup
}
\makeatother

\numberwithin{equation}{section}

\definecolor{darkred}{rgb}{.7,0,0}

\definecolor{green}{rgb}{0,0.7,0}

\usepackage{pgfplots}
\usepackage{tikz}

\usepackage[colorlinks=true,linkcolor=blue!70!black,urlcolor=black,citecolor=green!60!black,pdftex]{hyperref}
\hypersetup{final} 

\theoremstyle{plain}
\newtheorem{theorem}{Theorem}
\newtheorem{prop}[theorem]{Proposition}
\newtheorem{lemma}[theorem]{Lemma}
\newtheorem{coroll}[theorem]{Corollary}
\theoremstyle{definition}
\newtheorem{mydef}{Definition}

\newtheorem{assumption}{Assumption}
\theoremstyle{remark}
\newtheorem{remark}{Remark}
\newtheorem{example}{Example}

\newcommand{\tcr}[1]{\textcolor{black}{#1}}
\begin{document}
\title[]{Semiglobal optimal Feedback stabilization of autonomous systems via deep neural network approximation}

\author {Karl Kunisch\textsuperscript{$*$}}
\thanks{\textsuperscript{$*$}University of Graz, Institute of Mathematics and Scientific
Computing, Heinrichstr. 36, A-8010 Graz, Austria and Johann Radon Institute for Computational and Applied Mathematics
(RICAM), Austrian Academy of Sciences, Altenberger Stra\ss{}e 69, 4040 Linz, Austria, research supported by the ERC advanced grant 668998 (OCLOC) under the EU’s H2020 research program, ({\tt
karl.kunisch@uni-graz.at}).}
\author{Daniel Walter\textsuperscript{$\dagger$}}
\thanks{\textsuperscript{$\dagger$}Johann Radon Institute for Computational and Applied Mathematics
(RICAM), Austrian Academy of Sciences, Altenberger Stra\ss{}e 69, 4040 Linz, Austria,({\tt
daniel.walter@oeaw.ac.at}).}
\maketitle

\begin{abstract}
A learning approach for optimal feedback gains for nonlinear continuous time control systems is proposed and analysed. The goal is to establish a rigorous framework  for computing approximating  optimal feedback gains using neural networks.  The approach rests on two main ingredients. First, an optimal control formulation involving an ensemble of trajectories with 'control' variables given by the feedback gain functions. Second, an approximation to the feedback functions via realizations of neural networks. Based on universal approximation properties we prove the existence and convergence of optimal stabilizing neural network feedback controllers.
\end{abstract}

{\em{ Keywords:}}
optimal feedback stabilization, neural networks, Hamilton-Jacobi-Bellman equation, reinforcement learning.

{\em{AMS classification:}}
49J15,   
49N35, 
68Q32,   	
93B52,   	
93D15.   	

\section{Introduction} \label{sec:Introduction}
Closed loop optimal feedback control of nonlinear dynamical systems
remains to be a major challenge. Ultimately this involves the solution
of a Hamilton Jacobi Bellman (HJB) equation, a hyperbolic system whose
dimension is that of the state space. Thus, inevitably one is confronted
with a curse of dimensionality, possibly first stated in \cite{bell61}.
In the case of linear dynamics, the absence of control and state
constraints and a quadratic cost functional the HJB equation simplifies
to a Riccati equation. This Riccati equation has received a tremendous
amount of attention in the control literature, both from the point of
view of analysis as well as efficient numerical realization. One of the
frequently used approaches to nonlinear problems is therefore based on
applying the Riccati framework locally to  linear-quadratic
approximations of genuine nonlinear problems. This can be very
effective, but optimal controls based on the HJB equation differ from
their Riccati-based approximation especially in the transient phase.
Moreover linearization based feedback laws are often more expensive or
may even fail to reach the control objective, see e.g. \cite{bkp18,
kk18}. This suggests the search for deriving methods to obtain good
approximations to optimal feedback solutions for nonlinear systems,
which go beyond local linear-quadratic approximations.
In this paper we propose and analyze such an approach.

To present the main ideas and features of the methodology we focus on
stabilization problems of the form
\begin{equation} \label{def:refproblem}\tag{$P$}
\left\{
\begin{aligned}
\quad &\inf_{u \in L^2(0,\infty; \R^m)} \frac{1}{2}\int_{0}^\infty \left
( |Qy(t)|^2+ \beta |u(t)|^2 \right)~\de t \\
&s.t. \quad \dot{y}= {f}(y)+{B}u, \quad y(0)=y_0,
\end{aligned}
\right.
\end{equation}
where  $f$ describes the nonlinear dynamics, \tcr{$B\in \R^{n\times m}$}, is
the control operator, $Q \in \R^{n\times n} $ \tcr{is positive semi-definite}, $\beta>0$, and the optimal control is
sought in feedback from. The value function $V$ associated  to
\eqref{def:refproblem}, which assigns to each initial condition $y_0$
the value of the optimal cost, satisfies a HJB equation. Once
available,  the optimal control to \eqref{def:refproblem} can be
expressed in feedback form as $u^*(t)=- B^\top \nabla V(y(t))$. If $M$
degrees of freedom are used to discretize the HJB equation in each of
the directions $y_i$, then this results in a system of dimension $M^n$.
Except for small dimensions $n$ of the state equation this is unfeasible
and alternatives must be sought.
In this paper we propose to replace the control $u$ in
$\eqref{def:refproblem}$ by $F(y)$ and minimize with respect to the
feedback $F$ over an admissible set $\mathcal{H}$.
It can be expected that the effectiveness of such a procedure depends on
the location of the orbit $\mathcal {O}= \{y(t;y_0): t\in (0,\infty)\}$
within the state space $R^n$. To accommodate the case that $\mathcal{O}$
does not 'cover' the state-space sufficiently well, we propose to look
at the ensemble of orbits departing from a compact set $Y_0$ of initial
conditions and reformulate the problem accordingly. For this purpose we
introduce a probability measure $\mu$ on $Y_0$ and replace
\eqref{def:refproblem} by
\begin{equation} \label{def:refproblemensemle}\tag{$P_{Y_0}$}
\left\{
\begin{aligned}
\quad &\inf_{F\in \mathcal{H}} \frac{1}{2} \int_{Y_0} \int_{0}^\infty
\left ( |Qy(t;y_0)|^2+ \beta |F(y(t;y_0))|^2 \right)~\de t \,d\mu \\
&s.t. \quad \tcr{\dot{y}}(y_0)= {f}(y(y_0))+{B}F(y), \quad y(0)=y_0,
\end{aligned}
\right.
\end{equation}
Mathematical precision to this formulation will be given below. We note
that \eqref{def:refproblemensemle} represents a learning problem for the
feedback function $F$ on the basis of the data $\mathcal {O}_{ext}=
\{y(t;y_0): t\in (0,\infty), y_0 \in Y_0\}$.

For the  practical realization of \eqref{def:refproblemensemle} the
discretization or parametrization of the feedback function $F$ must be
addressed. To this end, due to their excellent approximation properties
in practice, we propose the use of  neural networks. Let $\theta $
denote the string of network parameters, and denote by $F_{\theta}$ the
approximation of $F$ by the realization of a neural network. Using this
parameterization in \eqref{def:refproblemensemle} we arrive at the
problem which represents the focus of the investigations in this paper:
\begin{equation} \label{def:reftheta}\tag{$P_{Y_0,\theta}$}
\left\{
\begin{aligned}
\quad &\inf_{\theta} \frac{1}{2} \int_{Y_0} \int_{0}^\infty \left (
|Qy(t;y_0)|^2+ \beta |F_\theta(y(t;y_0))|^2   \right)~\de t \, d\mu \\
&s.t. \quad \tcr{\dot{y}}(y_0)= {f}(y(y_0))+{B}F_\theta(y) , \quad y(0)=y_0.
\end{aligned}
\right.
\end{equation}
While \eqref{def:reftheta} conceptually represents well the goal of our
investigations, some amendments, as for example a properly chosen
regularization term, will be necessary. This will be described in the
later sections.
One of the main goals of this paper consists in the analysis of the
approximation properties of \eqref{def:refproblemensemle} by
\eqref{def:reftheta}, and in demonstrating the numerical feasibility of
the approach.

Let us very briefly and with no pretense for completeness mention some
of the methodologies which have been proposed and analysed to
approximate optimal feedback gains. First, there is the vast literature
on solving the HJB equations directly, see e.g. \cite{ff14, kkr18} and
obtaining the feedback by means of the verification theorem.
With the aim of increasing the dimension of computable HJB equations a
sparse grid approach was presented in \cite{gk17}.
More recently  significant progress was made in solving  high
dimensional HJB equations \tcr{by the use of} policy iterations (Newton method
applied to HJB) combined with tensor calculus techniques, \cite{ dkk19,
kk18}. The use of Hopf formulas was proposed in e.g. \cite{lr86,
cloy18}. Interpolation techniques, utilizing ensembles of open loop
solutions have been analyzed in the work of \cite{ngk19}, for example.

 From the machine learning point of view, the HJB approach is intimately
related to reinforcement learning (RL). For discrete time systems,
overviews on this subject are presented in the inspiring monographs
\cite{bert19,sutton}. See also the survey articles \cite{recht18,
vls14}. Reinforcement learning comprises e.g. value \tcr{function-based}
feedback approaches such as value and policy iterations, including the
Q-factor versions. The concept of RL rollout
algorithms,~\cite{bertmult}, is closely related to receding horizon
control. Another type of RL methods which do not rely on value function
approximation can be summarized under~\textit{parametrized policy
learning}. These include e.g. interpolation of feedback laws from
optimal open loop state-control pairs also known as expert supervised or
imitation learning, \cite{ROB-053}. The~parametrized policy method which
is closest to our approach is referred to as training by cost
optimization, \cite[Chapter 4]{bert19}, or policy
gradient,~\cite[Chapter 13]{sutton}, see also \cite {PeSa08}. We also
mention structural similarities between~\eqref{def:reftheta} and optimal
control problems under uncertainty,~\cite{garreis,kouri}.

While there are conceptual parallels between our research and selected
machine learning methodologies, here the focus lies on a rigorous
mathematical analysis of the approach sketched above. We shall also
provide examples illustrating the numerical feasibility of the proposed
method but we do not aim for sophistication in this respect within the
present paper.

The manuscript is structured as follows. In Section 2 we briefly
summarize our notation, Section 3 contains some aspects of optimal
feedback stabilization which are pertinent for our work. In Section 4 we
present the proposed learning formulation of the optimal feedback
stabilization problem. We then aim for an approximation result of the
optimal feedback law by neural networks. For this purpose we collect and
derive necessary results for neural networks in Section 5. Finally, in
Section 6, the set of admissible feedback laws in the learning problem
is replaced by realizations of neural networks. Existence of optimal
neural networks as well as convergence of the associated feedback laws
is argued. Amongst other issues, main difficulties that need to be
overcome arise due to the infinite time horizon and
the fact that the dynamics of the system as well as the feedback
functions are only assumed to be locally Lipschitz. Some aspects for the
practical realization of our approach are given in Section 7. Section 8
contains a brief description of numerical examples. The Appendix details
the proofs of several necessary technical results. \tcr{ In particular in Appendix D we give a version of the universal approximation theorem in $C^1$, with bounds on the coefficients, which we were not able to find in the  literature.}

\section{Notation}
Let a complete measure space~$(\Omega, \mathcal{A}, \mu)$ as well as a Banach space~$Y$ with norm~$\|\cdot\|_Y$ be given. The topological dual space of~$Y$ is denoted by~$Y^*$ and the corresponding duality paring is abbreviated as~$\langle \cdot, \cdot \rangle_{Y,Y^*}$. Following~\cite[p. 41]{DU77} we call~$\mathbf{y}\colon \Omega \to Y$~$\mu$-\textit{measurable} or~\textit{strongly measurable} if there exists a sequence~$\{\mathbf{y}_i\}_{i\in \N}$ of simple functions with
\begin{align*}
\lim_{i \rightarrow \infty}\|\mathbf{y}_i(\omega)-\mathbf{y}(\omega)\|=0 \quad \text{for}~\mu-a.e.~\omega \in \Omega.
\end{align*}
Here a mapping~$\mathbf{y}_i \colon \Omega \to Y$ is~\textit{simple} if~$\mathbf{y}_i=\sum^{N_i}_{j=1} y_j \chi_{\mathcal{A}_j}$ for some~$y_1,\dots,y_{N_i} \in Y$ and characteristic functions~$\chi_{A_j}$,~$A_j \in \mathcal{A}$.
If the mapping
\begin{align*}
\widehat{\mathbf{y}}\colon \Omega \to Y, \quad \widehat{\mathbf{y}}(\omega)=\langle \mathbf{y}(\omega),y^* \rangle_{Y,Y^*}
\end{align*}
is~$\mu$-measurable for all~$y^*\in Y^*$ we term~$\mathbf{y}$ weakly measurable. Strong and weak measurability are stable under pointwise almost everywhere convergence i.e. if~$\{\mathbf{y}_i\}_{i\in\N}$ is a sequence of~strong (weak) measurable functions and~$\mathbf{y}$ is defined by~$\mathbf{y}(\omega)=\lim_{i\rightarrow \infty} y_i(\omega)$ for~$\mu$-a.e.~$\omega \in \Omega$ then~$\mathbf{y}$ is strong (weak) measurable, respectively,~\cite[p. 41]{DU77}. If~$Y$ is separable, Pettis' Theorem, cf.~\cite[p. 42]{DU77}, states the equivalence of strong and weak measurability.

We further introduce~$L^p_\mu(\Omega, \mathcal{A}; Y)$,~$p \in [1, \infty]$, as the Bochner space of equivalence classes of strongly measurable functions~$\mathbf{y} \colon \Omega \to Y$ such that the associated norms
\begin{align*}
\|\mathbf{y}\|_{L^p_\mu(\Omega, \mathcal{A}; Y)}= \left(\int_{\Omega} \|\mathbf{y}(\omega)\|^p_Y~\de \mu(\omega)\right)^{\frac{1}{p}}, \quad p \in [1, \infty)
\end{align*}
as well as
\begin{align*}
\|\mathbf{y}\|_{L^\infty_\mu(\Omega, \mathcal{A}; Y)}= \inf_{\substack{O \in \mathcal{A}, \\ \mu(O)=0}} \sup_{\omega \in \Omega \setminus O } \|\mathbf{y}(\omega)\|_Y
\end{align*}
are finite. If the measure~$\mu$ and/or the~$\sigma$-algebra~$\mathcal{A}$ are clear from the context we omit them in the description of the space and simply write~$L^p_\mu(\Omega; Y)$ and~$L^p(\Omega; Y)$, respectively, for shortage. If~$Y$ is reflexive and~$p \in (1,\infty)$ we identify
\begin{align*}
L^p_\mu(\Omega, \mathcal{A}; Y)^* \simeq L^q_\mu(\Omega, \mathcal{A}; Y^*),
\end{align*}
where~$1/q+1/p=1$. For all~$\mathbf{y} \in L^p_\mu(\Omega, \mathcal{A}; Y)$ and~$\mathbf{y}' \in L^q_\mu(\Omega, \mathcal{A}; Y^*)$ we have
\begin{align*}
\langle \mathbf{y}, \mathbf{y}' \rangle_{L^p_\mu(\Omega, \mathcal{A}; Y), L^q_\mu(\Omega, \mathcal{A}; Y^*)}= \int_{\Omega} \langle \mathbf{y}(\omega), \mathbf{y}'(\omega) \rangle_{Y,Y^*}~ \de \mu(\omega).
\end{align*}
If~$Y$ is a Hilbert space with respect to the norm induced by the inner product~$(\cdot, \cdot)_Y$ then~$L^2_\mu(\Omega, \mathcal{A}; Y)$ is also a Hilbert space with inner product
\begin{align*}
(\mathbf{y}_1, \mathbf{y}_2)_{L^2_\mu(\Omega, \mathcal{A}; Y)}= \int_{\Omega} (\mathbf{y}_1(\omega), \mathbf{y}_2(\omega))_Y~\de \mu(\omega) \quad \forall \mathbf{y}_1, \mathbf{y}_2 \in L^2_\mu(\Omega, \mathcal{A}; Y).
\end{align*}
For further references see \cite{DU77,Ed65}.

For a metric space~$X$ we denote the space of bounded continuous functions between~$X$ and~$Y$ by~$\mathcal{C}_b(X;Y)$. It forms a Banach space together with the supremum norm
\begin{align*}
\|\varphi\|_{\mathcal{C}_b(X;Y)}= \sup_{x \in X} \|\varphi(x)\|_Y \quad \text{for } \varphi \in \mathcal{C}_b(X; Y).
\end{align*}
Throughout this paper, let~$Y_0 \subset \R^n$  be compact, and set~$I:=[0,\infty)$. Define
\begin{align*}
W_\infty= \{\,y \in L^2(I; \R^n)\;|\;\dot{y}\in L^2(I; \R^n)\,\}.
\end{align*}
Here the temporal derivative is understood in the distributional sense. We equip~$W_\infty$ with the norm induced by the inner product
\begin{align*}
(y_1,y_2)_{W_\infty}=(\dot{y}_1,\dot{y}_2)_{L^2(I;\R^n)}+(y_1,y_2)_{L^2(I;\R^n)} \quad \text{for } y_1, y_2 \in W_\infty,
\end{align*}
making it a Hilbert space.
We frequently use the following properties of functions in $W_\infty$. First,~$W_\infty$ continuously embeds into~$\Cc_b({I;\R^n})$, i.e. there exists a constant $C>0$ such that
\begin{equation}\label{eq:kk18}
\|y\|_{\mathcal{C}_b(I; \R^n)} \leq C \wnorm{y}, \text{ for all } y\in W_\infty,
\end{equation}
see e.g. \cite[Theorem 3.1]{LioM72}. Second, we have
\begin{equation}\label{eq:kk19}
\lim_{t\to \infty} |y(t)| = 0  , \text{ for all } y\in W_\infty,
\end{equation}
see \cite[pg. 521] {CK17}.

Given two Banach spaces~$X$ and~$Y$, the space of linear and continuous functionals between~$X$ and~$Y$ is denoted by~$\mathcal{L}(X,Y)$. It is a Banach space with respect to the canonical operator norm given by
\begin{align*}
\|A\|_{\mathcal{L}(X,Y)}= \sup_{\|x\|_X=1} \|Ax\|_Y \quad \forall A \in \mathcal{L}(X,Y).
\end{align*}
We further abbreviate~$\mathcal{L}(X):=\mathcal{L}(X,X)$.
\section{Preliminaries}
\label{sec:formulation}
In the following we consider a controlled autonomous  dynamical system of the form
\begin{align}
\label{eq:openloop}
\dot{y}= \mathbf{f}(y)+\mathcal{B}u \quad \text{in}~L^2(I;\R^n), \quad y(0)=y_0,
\end{align}
where~$y_0 \in \R^n$ is a given initial condition. The nonlinearity is described by a Nemitsky operator
\begin{align*}
\mathbf{f} \colon W_\infty \to L^2(I;\R^n), \quad \mathbf{f}(y)(t)= f(y(t)) \quad \text{for a.e.}~t \in I,
\end{align*}
where~$f \colon \R^n \to \R^n$. The requirements on $f$ will be specified in Assumption~\ref{ass:feedbacklaw} below. The dynamics of the system can be influenced by choosing a control function~$u \in L^2(I; \R^m)$ which enters via a bounded linear operator
\begin{align*}
 \mathcal{B} \colon L^2(I;\R^m) \to L^2(I;\R^n), \quad \mathcal{B}u(t)=Bu(t) \quad \text{for a.e.}~t \in I
\end{align*}
induced by a matrix~$B\in\R^{n \times m}$.
\subsection{Open loop optimal control}
Our interest lies in determining a control~$u^*\in L^2(I; \R^m)$ which keeps the associated solution~$y^* \in W_\infty$ to~\eqref{eq:openloop} small in a suitable sense and steers it to~$0$ as~$t \rightarrow \infty$. This can be achieved by computing an optimal pair~$(\bar{y}, \bar{u})$ to the \textit{open loop} infinite horizon optimal control problem
%
\begin{equation} \label{def:openloopproblem}\tag{$P_\beta^{y_0}$}
\left\{
\begin{aligned}
\quad &\inf_{y \in W_\infty,\, u \in L^2(I; \R^m)} \frac{1}{2}\int_{I} \left ( |Qy(t)|^2+ \beta |u(t)|^2   \right )~\de t \\
&s.t. \quad \dot{y}= \mathbf{f}(y)+\mathcal{B}u, \quad y(0)=y_0,
\end{aligned}
\right.
\end{equation}
for some positive semi-definite matrix~$Q \in \R^{n \times n}$ and a fixed cost parameter~$\beta>0$. In the following we silently assume that the infimum in~\eqref{def:openloopproblem} is attained for every~$y_0 \in \R^n$. For abbreviation we further set
\begin{align*}
J \colon W_\infty \times L^2(I; \R^m) \to \R, \quad J(y,u)=\frac{1}{2}\int_{I} \left ( |Qy(t)|^2+ \beta |u(t)|^2   \right )~\de t.
\end{align*}
The open loop approach comes with several drawbacks. First, open loop controls do not take into account possible  perturbations of the dynamical system.  Second, determing the control action for a new initial condition requires to solve~\eqref{def:openloopproblem}  again.

\subsection{Semiglobal optimal feedback stabilization}\label{sec3.2}
The aforementioned disadvantages of open loop optimal controls serve as a motivation to study the \textit{semiglobal optimal feedback problem} associated to~\eqref{def:openloopproblem} and to construct the optimal control (at any time $t\ge 0$) as a function of the associated state at time $t$. More precisely, given a compact set~$Y_0 \subset \R^n $ containing the origin, we are looking for a function~$F^* \colon \R^n \to\R^m$ which induces a Nemitsky operator
\begin{align*}
\F^* \colon W_\infty \to L^2(I; \R^m) , \quad \F^*(y)(t)=F^*(y(t)) \quad \text{for a.e.}~t \in I,
\end{align*}
such that for every~$y_0 \in Y_0$ the~\textit{closed loop system}
\begin{align} \label{eq:cloloop}
\dot{y}= \mathbf{f}(y)+ \mathcal{B} \mathcal{F}^*(y), \quad y(0)=y_0
\end{align}
admits a unique solution~$y^*( y_0 ) \in W_\infty$ and~$(y^*( y_0), \mathcal{F}^*(y^*(y_0 )))$ is a minimizing pair of~\eqref{def:openloopproblem}.

As a starting point to its solution we define the value function associated to the open loop problem~\eqref{def:openloopproblem} by
\begin{align} \label{def:valuefunc}
V(y_0):= \min_{y \in W_\infty, u \in L^2(I; \R^m)} J(y,u) \quad s.t. \quad \dot{y}= \mathbf{f}(y)+\mathcal{B}u, \quad y(0)=y_0
\end{align}
for~$y_0\in \R^n$. If~$V$ is continuously differentiable in a neighborhood of~$y_0 \in \R^n$ it solves
the stationary~\textit{Hamilton-Jacobi-Bellman equation} (HJB)
\begin{align} \label{eq:HJB}
(f(y_0), \nabla V(y_0))_{\R^n}- \frac{1}{2\beta}|B^*\nabla V(y_0)|^2+ \frac{1}{2} |Q y_0|^2=0
\end{align}
in the classical sense, see e.g. \cite[pg 6]{FS06}.
Moreover, once computed, an optimal control for~\eqref{def:openloopproblem} in feedback form is given by the verification theorem as
$u=-\frac{1}{\beta}B^T\nabla V(y)$, see \cite[Theorem I.7.1]{FS06}. This leads to the closed loop system in optimal feedback form
\begin{align*}
\dot{y}=\mathbf{f}(y)- \frac{1}{\beta}\mathcal{B}\mathcal{B}^* \nabla \mathcal{V}(y), \; y(0)=y_0,\; \text{where } \mathcal{B}^* \nabla \mathcal{V}(y)(t)=B^\top \nabla V(y(t)), \, \text{for a.e.}~t \in I,
\end{align*}
yielding a function~$y^*(y_0)\in W_\infty$. Thus
\begin{align*}
\left(y^*(y_0), -\frac{1}{\beta}\mathcal{B}^* \nabla \mathcal{V}(y^*(y_0))\right) \in \argmin \eqref{def:openloopproblem}
\end{align*}
and the function
\begin{align*}
F^*=-\frac{1}{\beta} B^\top \nabla V
\end{align*}
solves the semiglobal optimal feedback problem as formulated above.

Realizing this closed loop system  requires a solution of the HJB equation~\eqref{eq:HJB}, which is a partial  differential equation in $\R^n $. This makes the HJB based approach numerically challenging  or infeasible if ~$n\in \N$ is large.

In this paper we adopt a different viewpoint for the semiglobal feedback problem which does not involve the solution of the HJB equation. An optimal feedback law is determined  as a minimizer to a suitable optimal control problem which incorporates the closed loop system as a constraint. The following assumptions on the dynamical system~\eqref{eq:openloop}, the open loop problems~\eqref{def:openloopproblem}, and the value function~$V$, are assumed throughout.

\begin{assumption} \label{ass:feedbacklaw}
\begin{itemize}
\item[\textbf{A.1}] The function~$f \colon \R^n \to \R^n$ is continuously differentiable with $f(0)=0$. Its Jacobian~$Df \colon \R^n \to \R^{n \times n} $ is Lipschitz continuous on compact sets.

\item [\textbf{A.2}] \tcr{ There exists a neighborhood   ~$\mathcal{N}(Y_0)\subset \R^n$  of ~$Y_0$, a function ~$F^* \colon \mathcal{N}(Y_0) \to \R^m$ with~$F^*(0)=0$   as well as the induced Nemitsky operator~$\F^* \colon W_\infty \to L^2(I;\R^m)$ such that the closed loop system}
\begin{align*}
\dot{y}= \mathbf{f}(y)+ \mathcal{B} \F^*(y), \quad y(0)=y_0
\end{align*}
admits a unique solution~$y^*(y_0)\in W_\infty$ for every~$y_0 \in \mathcal{N}(Y_0)$. Moreover we have
\begin{align*}
(y^*(y_0), \F^*(y^*(y_0))) \in \argmin \eqref{def:openloopproblem} \quad \forall y_0 \in \mathcal{N}(Y_0).
\end{align*}
\item [\textbf{A.3}] The mapping
\begin{align*}
\mathbf{y}^* \colon \mathcal{N}(Y_0) \to W_\infty, \quad y_0 \mapsto y^*(y_0)
\end{align*}
is continuously differentiable. Moreover there is a constant~$M>0$ \tcr{such that}
\begin{align} \label{eq:apriorioptfeedbackass}
\wnorm{\mathbf{y}^*(y_0)}\leq M |y_0| \quad \forall y_0 \in \mathcal{N}(Y_0).
\end{align}

\item [\textbf{A.4}]
There holds
\begin{align*}
F^* \in \mathcal{C}^1\left(\bar{B}_{ 2\widehat{M}}(0);\R^m\right), \quad DF^* \in \operatorname{Lip}\left(\bar{B}_{ 2\widehat{M}}(0);\R^{m \times n } \right).
\end{align*}
\end{itemize}
\end{assumption}

We close this section with several remarks concerning Assumption~\ref{ass:feedbacklaw}.
\begin{remark} \label{rem:assumptionsforvalue}
Following the previous discussions, the canonical candidate for the optimal feedback law in \tcr{~$\mathbf{A.2}$} is given by~$F^*=-\frac{1}{\beta} B^\top \nabla V $ provided the value function~$V$ is sufficiently smooth. A discussion of the regularity assumptions according to~$\mathbf{A.4}$, in this context can be found in Appendix~\ref{app:smoothnessofvalue}.
\end{remark}

\begin{remark}\label{rem:semiglobal} Sufficient conditions on the problem data in \eqref{def:refproblem} which imply $\mathbf{A.2}$ and $\mathbf{A.3}$  are given, for example,  by the following scenario: $f(y)=Ay + g(y)$, where $A\in \R^{n\times n}$ and $g:\R^n\to \R^n$, with $g(0)=0$, where we suppose that
\begin{equation} \label{semiglass1}
\left\{
\begin{array}l
(A,B) \text{ is stabilizable in the sense that there exists } K\in \R^{m\times n}\\ \text{ and  } \mu <0, \text{ such that }
((A+BK)y,y)_{\R^n} \le -\mu |y|^2, \text{ for all }  y \in\R^n,
\end{array}
\right.
\end{equation}
and
\begin{equation} \label{semiglass2}
\left\{
\begin{array}l
g \in C^2(\R^n,\R^n) \text{ with } \|Dg(x)\| \le L \text{ for some } L \in (0,\mu) \text{ and all } x\in\mathcal{N}(Y_0). \\
\text{Moreover } \|Dg(x)\| \le \tilde L \text{ for some }  \tilde L \text{ independent of } x \in \R^n.
\end{array}
\right.
\end{equation}
The proof is given in Appendix~\ref{app:semiglobal}.
In Appendix \ref{app:smoothnessofvalue} a sufficient condition is provided which guarantees that under less stringent conditions $\mathbf{A.2}$ and  $\mathbf{A.3}$ hold in a neighborhood of the origin.

\end{remark}

\begin{remark}
In~$\mathbf{A.3}$ we assume the continuous Fr\'echet differentiability of the solution operator~$\mathbf{y}^*$ to the closed loop system. Utilizing the regularity assumptions on~$f$ and~$F^*$ specified in~$\mathbf{A.1}$ and~$\mathbf{A.4}$,  we conclude the Fr\'echet differentiability of the induced Nemitsky operators~$\mathbf{f}$ and~$\F^*$ at~$\mathbf{y}^*(y_0)$ with~$y_0 \in \mathcal{N}(Y_0)$. Consequently,~denoting the Fr\'echet derivative of~$\mathbf{y}^*$ at~$y_0 \in \mathcal{N}(Y_0)$ by~$\delta \mathbf{y}^*(y_0)$, the function~$\delta y:= \delta\mathbf{y}^*(y_0)\delta y_0\in W_\infty$ fulfills
\begin{equation}\label{eq:linearized}
\dot{\delta y}= D \mathbf{f}(\mathbf{y}^*(y_0)) \delta y+ \mathcal{B} D \F^*(\mathbf{y}^*(y_0))\delta y , \quad \delta y(0)=\delta y_0
\end{equation}
for every~$y_0 \in \mathcal{N}(Y_0),~\delta y_0 \in \R^n$.
Note that there exists a constant~$c>0$ with
\begin{align*}
\|\delta \mathbf{y}^*(y_0)\|_{\mathcal{L}(\R^n,W_\infty)} \leq c \quad \forall y_0 \in Y_0
\end{align*}
due to the continuity of~$\delta \mathbf{y}$ and the compactness of~$Y_0$.
\end{remark}

Next we discuss a sufficient condition for~\eqref{eq:apriorioptfeedbackass} to hold.
\begin{prop}
Assume that $\mathbf{A.1}$ and $\mathbf{A.3}$ hold,
 and that the value function $V$ associated to $(P_\beta^{y_0})$ is $C^1$ on $\mathbb{R}^n$, and satisfies $V(y)\le \bar{\kappa}|y|^2 $  \tcr{ for some~$\bar{\kappa}>0$} for all  $y \in \mathbb{R}^n$. Further assume \tcr{that}
\begin{itemize}
\item[a)] $f$ is globally Lipschitz continuous on $\mathbb{R}^n$, or
\item[b)] there exists $\underline{\kappa}_1 > 0, \underline{\kappa}_2 \in \mathbb{R}^n$ such that $V(y) \ge \underline{\kappa_1} |y| - \underline{\kappa_2}$ for all $y \in \mathcal{N}(Y_0)$.
\end{itemize}
Then~\eqref{eq:apriorioptfeedbackass} holds.
\end{prop}
\begin{proof}
The HJB equation associated to $(P_\beta^{y_0})$ is given by
\begin{equation} \label{aux10}
(f(y), \nabla V(y))_{\mathbb{R}^n} - \frac{1}{2\beta} |B^T\nabla V(y)|^2 + \frac{1}{2} |y|^2 = 0 \text{ for } y \in \mathbb{R}^n.
\end{equation}
For $y_0 \in \mathcal{N}(Y_0)$ the optimal closed loop system is given by
\begin{equation} \label{aux11}
\dot{y} = f(y) + BF(y), y(0) = y_0,
\end{equation}
with $F(y) = -\frac{1}{\beta}B^T \nabla V(y)$, where $y = y(t), t > 0$ (and we drop the superscript *). Testing \eqref{aux11} with $\nabla V (y(t))$, we obtain
\begin{equation} \label{aux12}
\frac{d}{dt} V(y(t)) = (f(y(t)), \nabla V(y(t))) - \frac{1}{2\beta} |B^T \nabla V(y(t))|^2 - \frac{1}{2\beta} |B^T \nabla V(y(t))|^2.
\end{equation}
Utilizing \eqref{aux10} this implies that
\begin{equation} \label{aux13}
\frac{d}{dt} V(y(t)) = - \frac{1}{2} |y(t)|^2 - \frac{1}{2\beta} |B^T \nabla V(y(t))|^2\tcr{,}
\end{equation}
from which we deduce that
\begin{equation} \label{aux14}
2 V(y(t)) + \int_0^t |y(s)|^2 ds + \frac{1}{\beta} \int_0^t |B^T \nabla V(y(s))\tcr{|}^2 ds \le 2 V (y_0),
\end{equation}
for all $t \ge 0$. If $f$ is globally Lipschitz with Lipschitz-constant $L$ then by \eqref{aux11} and \eqref{aux14}
\begin{equation*}
\begin{aligned}
\int_0^\infty |\dot{y}(t)|^2 dt & \le L^2 \int_0^\infty |y(t)|^2 dt + \frac{1}{\beta^2} \int_0^\infty |BB^T \nabla V(y(t))|^2 dt \\
& \le L^2 \int_0^\infty |y(t)|^2 dt + \frac{\|B\|^2}{\beta^2} \int_0^\infty |B^T \nabla V(y(t))|^2 dt \\
& \le \sup(L^2, \frac{1}{\beta}\|B\|^2) \int_0^\infty (|y(t)|^2 + \frac{1}{\beta} |B^T \nabla V(y(t))|^2 dt \\
& \le 2 ~\sup(L^2, \frac{1}{\beta}\|B\|^2) V(y_0).
\end{aligned}
\end{equation*}
This implies that
\begin{equation*}
\begin{aligned}
\tcr{\|\mathbf{y}(y_0)\|}^2_{W_\infty} & =  \int_0^\infty (|y|^2 + |\dot{y}|^2) dt \le\tcr{ 2 (1 + \sup(L^2, \frac{1}{\beta}\|B\|^2))\,V(y_0)}\\
& \le 2 \bar{\kappa} \big(1 + \sup(L^2, \frac{1}{\beta}\|B\|^2)\big) |y_0|^2,
\end{aligned}
\end{equation*}
which gives the desired estimate for $\mathbf{A.2}$ to hold with $M = \sqrt{2 \bar{\kappa} (1 + \sup(L^2, \frac{1}{\beta}\|B\|^2)}$. Under condition (b) from \eqref{aux14}
$$\sup_{t \in [0,\infty)} \underline{\kappa}_1 |y(t)| - \underline{\kappa}_2 \le V(y_0) \le \bar{V} = \sup \{ V(y_0): y_0 \in \mathcal{N}(Y_0)\}.$$
We choose $L$ as the Lipschitz constant of $f$ on the ball with radius $\frac{1}{\underline{\kappa}_1} (\bar{V} + \underline{\kappa}_2)$, and center $0$. Now we can proceed as for condition (a) to arrive at the desired estimate.
\end{proof}

In the case of a linear-quadratic optimal control problem a lower bound as in assumption (b) above is satisfied with $V(y) \ge \underline{\kappa_1} |y|^2$ if the control-system is observable. 

\section{Optimal Feedback stabilization as learning problem}
This section leads up to the formulation of the optimal feedback problem as ensemble-optimal control problem. For this purpose we choose  a complete probability space ~$(Y_0,\mathcal{A},\mu)$ for the initial conditions.
\tcr{The measure~$\mu$ describes the "training set" of initial conditions from which we will learn a suitable feedback law. For example,
the choice of~$\mu$ as the Lebesgue measure treats the case when all the elements
$y_0\in Y_0$ contribute to the learning of ${\mathcal {F}}$. On the contrary, choosing~$\mu$ as a conic combination of Dirac delta functionals describes the case of finitely many
pre-selected initial conditions in the training set. However the following results are not restricted to these canonical examples}.
We recall the definition of $\mathcal{Y}_{ad}$ in Assumption 1 and  endow the set  with the norm of~$W_\infty$.
By $\mathbf{A.4}$ every~$y \in \mathcal{Y}_{ad}$ fulfills~$\|{y}\|_{\mathcal{C}_b(I;\R^n)}\leq 2 \widehat{M}$. According to Proposition~\ref{prop:nemitskylip} in Appendix \ref{app:nemitsky}, every function~$F \in \operatorname{Lip}\left(\bar{B}_{ 2\widehat{M}}(0);\R^m\right)$with $F(0)=0$ induces a Lipschitz continuous Nemitsky operator on~$\mathcal{Y}_{ad}$.
Following this observation the set of admissible feedback laws is defined as
\begin{align*}
\mathcal{H}= \left\{\,\F \colon \mathcal{Y}_{ad} \to L^2(I;\R^m)\;|\; \F(y)(t)=F(y(t)) ~ \forall t \in I,~ F \in \operatorname{Lip}\left(\bar{B}_{ 2\widehat{M}}(0);\R^m\right),~F(0)=0\,\right\}.
\end{align*}
We further introduce the closed loop system associated to~$\F \in \mathcal{H}$ as
\begin{align} \label{eq:closedloopstrong}
\dot{y}=\mathbf{f}(y)+ \mathcal{B}\F(y)\quad \text{in}~L^2(I;\R^n), \quad y(0)=y_0,
\end{align}
where~$y_0 \in Y_0$ and the state~$y$ will be  searched for in~$W_\infty$. By the Gronwall lemma it is imminent that~\eqref{eq:closedloopstrong} admits at most one solution in~$\mathcal{Y}_{ad}$.

In order to facilitate the analysis, we will work with an equivalent variational reformulation of \eqref{eq:closedloopstrong} which couples the closed loop dynamics and the initial condition. That is, we require the state~$y \in \tcr{\Ya}$ to fulfill
\begin{align} \label{eq:closedloopweak}
a(\F)(y,y_0, \phi)=0 \quad \forall \phi \in W_\infty,
\end{align}
where the semilinear form~$a(\cdot)(\cdot,\cdot ,\cdot) \colon \mathcal{H}\times \tcr{\Ya}\times \R^n \times W_\infty \to \R$ is defined by
\begin{align*}
\, \\
a(\mathcal{F})(y,&y_0, \phi)\\&= \scalarl{\dot{y}}{\phi}-\scalarl{\mathbf{f}(y)}{\phi}-\scalarl{\mathcal{B}\F(y)}{\phi}+\scalarrn{y(0)-y_0}{\phi(0)}.
\\
\end{align*}
The next proposition establishes the equivalence of~\eqref{eq:closedloopstrong} and~\eqref{eq:closedloopweak}.
\begin{prop} \label{prop:equivalofweakandstrong}
Let $y_0 \in Y_0$ and $\F \in \mathcal{H}$ be given.
A function~$y \in \mathcal{Y}_{ad}$ fulfills
\begin{align*}
\dot{y}=\mathbf{f}(y)+ \mathcal{B}\F(y)\quad \text{in}~L^2(I;\R^n),~y(0)=y_0,
\end{align*}
if and only if
\begin{align*}
a(\F)(y,y_0, \varphi)=0 \quad \forall \varphi \in W_\infty.
\end{align*}
In particular,~\eqref{eq:closedloopweak} admits at most one solution in~$\mathcal{Y}_{ad}$.
\end{prop}
\begin{proof}
First note that for fixed~$y \in \Ya$,~$y_0 \in Y_0$ and~$\F \in \mathcal{H}$ the mapping
\begin{align*}
h \colon L^2(I; \R^n) \to \R, \quad \phi \mapsto (\dot{y}-\mathbf{f}(y)-\mathcal{B}\F(y),\phi)_{L^2(I;\R^n)}
\end{align*}
defines a linear and continuous functional on~$L^2(I;\R^n)$.
Obviously, if~$y \in \Ya$ fulfills~\eqref{eq:closedloopstrong} we have
\begin{align*}
h(\phi)=0, \quad \scalarrn{y(0)-y_0}{\phi(0)}=0
\end{align*}
for every~$\phi \in W_\infty$. Therefore~\eqref{eq:closedloopweak} holds.

Conversely assume that~$y \in \Ya$ satisfies~\eqref{eq:closedloopweak}. Choosing a test function~$\phi \in \mathcal{C}^\infty(I; \R^n)$ with compact support we have~$\phi(0)=0$ and consequently
\begin{align*}
0=a(\mathcal{F})(y,y_0, \phi)= \scalarl{\dot{y}}{\phi}-\scalarl{\mathbf{f}(y)}{\phi}-\scalarl{\mathcal{B}\F(y)}{\phi}.
\end{align*}
By a density argument we conclude
\begin{align*}
0= \scalarl{\dot{y}}{\phi}-\scalarl{\mathbf{f}(y)}{\phi}-\scalarl{\mathcal{B}\F(y)}{\phi}
\end{align*}
for every~$\phi \in L^2(I;\R^n)$. This proves
\begin{align*}
\dot{y}= \mathbf{f}(y)+ \mathcal{B}\F(y) \quad \text{in}~L^2(I; \R^n).
\end{align*}
Due to~$W_\infty \hookrightarrow L^2(I;\R^n)$ we now obtain
\begin{align*}
0=\scalarrn{y(0)-y_0}{\phi(0)} \quad \forall \phi \in W_\infty.
\end{align*}
Since the mapping~$\phi \to \phi(0)$ from~$W_\infty$ to $\R^n$ is surjective,~$y(0)=y_0$ has to hold. This finishes the proof of the first statement.

If~\eqref{eq:closedloopweak} admits a solution in~$\mathcal{Y}_{ad}$ it is unique due to the equivalence of~\eqref{eq:closedloopstrong} and~\eqref{eq:closedloopweak}.

\end{proof}
Observe that the variational formulation~\eqref{eq:closedloopweak} and thus also its solution, are parameterized as functions of the initial condition~$y_0 \in Y_0$. We introduce a suitable solution concept for the closed loop system which takes this dependence into account. For this purpose consider the mapping
\begin{align*}
\abd \colon \mathcal{H} \times L^2_\mu(Y_0;\Ya) \times \Lzwoboch\to \R
\end{align*}
with
\begin{align*}
 \abd(\mathcal{F})(\mathbf{y}, \Phi)= \int_{Y_0} a(\F)(\mathbf{y}(y_0),y_0,\Phi(y_0))~\de \mu (y_0),
\end{align*}
for every triple~$(\F, \mathbf{y}, \Phi)\in \mathcal{H} \times L^2_\mu(Y_0;\tcr{\Ya}) \times \Lzwoboch$.

\begin{mydef}\tcr{Given a feedback} law~$\F\in \mathcal{H}$, an element~$\mathbf{y}\in \Linfboch$ is called an ensemble solution of~\eqref{eq:closedloopstrong} if  ~$\Linfbochnorm{\mathbf{y}}\leq 2 M_{Y_0}$ and
\begin{align} \label{eq:parametrizedstate}
\abd(\mathcal{F})(\mathbf{y}, \Phi)=0 \quad \forall \Phi \in \Lzwoboch.
\end{align}
\end{mydef}
We next argue that ensemble solutions are unique and fulfill~\eqref{eq:closedloopweak} in an almost everywhere sense.
\begin{prop} \label{prop:parametrizedstate}
Let~$\F \in \mathcal{H}$ be given and let Assumption~\ref{ass:feedbacklaw} hold. Then there exists at most one ensemble solution~$\mathbf{y}\in  \Linfboch$ with~$\Linfbochnorm{\mathbf{y}}\leq 2 M_{Y_0}$, to~\eqref{eq:parametrizedstate}. If it exists then
\begin{align} \label{eq:parametrizedpointwise}
a(\F)(\mathbf{y}(y_0),\phi)=0 \quad \forall \phi \in W_\infty
\end{align}
for~$\mu$-a.e.~$y_0 \in Y_0$. Conversely, if there exists~$\mathbf{y}\in \Linfboch$ such that~$\Linfbochnorm{\mathbf{y}}\leq 2M_{Y_0}$ and~$\mathbf{y}(y_0)$ satisfies~\eqref{eq:parametrizedpointwise} for~$\mu$-a.e.~$y_0\in Y_0$, then it is the unique ensemble solution to~\eqref{eq:parametrizedstate}.
\end{prop}
\begin{proof}
Let~$\F \in \mathcal{H}$ be given.
We first prove the second claim and assume that~\eqref{eq:parametrizedstate} admits an ensemble solution~$\mathbf{y}\in  \Linfboch$. Consequently there exists~$O \in \mathcal{A}$,~$\mu(O)=0$, with~$\mathbf{y}(y_0) \in \mathcal{Y}_{ad}$ for all~$y_0 \in Y_0 \setminus O$. Utilizing the Lipschitz continuity of~$\mathbf{f}$ and~$\F$ on~$\mathcal{Y}_{ad}$, see Proposition~\ref{prop:nemitskylip}, as well as~$\mathbf{f}(0)=\F(0)=0$ we estimate
\begin{align*}
|\scalarl{\mathbf{f}(\mathbf{y}(y_0))}{\phi}&+\scalarl{\mathcal{B}\F(\mathbf{y}(y_0))}{\phi}| \\&\leq (L_{f,2\widehat{M}}+\|B\|L_{F,2 \widehat{M}})\|\mathbf{y}(y_0)\|_{L^2(I;\R^n)}\wnorm{\phi}
\end{align*}
for~$y_0\in Y_0 \setminus O$ and~$\phi \in W_{\infty}$. This implies that for~$\mu$-a.e.~$y_0 \in Y_0$ the mapping
\begin{align*}
\phi \mapsto a(\F)(\mathbf{y}(y_0),y_0, \phi) \text{ for } \phi \in W_\infty
\end{align*}
defines an element of the dual space~$W^*_\infty$.
Further note that
\begin{align*}
\Phi^* \colon Y_0 \to W^*_\infty, \quad y_0 \mapsto a(\F)(\mathbf{y}(y_0),y_0,\cdot),
\end{align*}
defined in a~$\mu$-almost everywhere sense, is~$\mu$-measurable and
\begin{align*}
\langle \Phi^*(y_0), \phi \rangle_{W^*_\infty, W_\infty}= a(\F)(\mathbf{y}(y_0),\phi) \leq c \wnorm{\mathbf{y}(y_0)} \wnorm{\phi},
\end{align*}
with a constant~$c>0$ only depending on~$B$, the Lipschitz constants~$L_{f,2\widehat{M}}$ and~$L_{F,2\widehat{M}}$ of $f$ and $F$ on $\bar B_{2\widehat M}(0)$, but not on~$y_0$. Thus~$\Phi^* \in L^2_\mu(Y_0; W^*_\infty)$. Since~$\mathbf{y}$ is an ensemble solution to~\eqref{eq:parametrizedstate} we get
\begin{align*}
\langle \Phi^*, \Phi \rangle= \int_{Y_0} a(\F)(\mathbf{y}(y_0),y_0,\Phi(y_0))~\de \mu(y_0)= \abd(\F)(\mathbf{y},\Phi)=0 \quad \forall \Phi \in \Lzwoboch
\end{align*}
where~$\langle \cdot, \cdot\rangle$ denotes the duality pairing between~$\Lzwoboch$ and~$L^2_\mu(Y_0, W^*_\infty)$. In particular, we conclude
\begin{align*}
0=\langle \Phi^*(y_0), \phi  \rangle_{W^*_\infty, W_\infty}=a(\F)(\mathbf{y}(y_0),\phi) \quad \forall \phi \in W_\infty
\end{align*}
and~$\mu$-a.e.~$y_0 \in Y_0$.
This gives~~\eqref{eq:parametrizedpointwise}.

Next we prove the uniqueness of ensemble solutions. If~$\mathbf{y}_i \in \Linfboch$,~$i=1,2$, are two ensemble solutions of~\eqref{eq:parametrizedstate} there exist sets~$O_i\in \mathcal{A}$,~$\mu(O_i)=0$,~$i=1,2$, with
\begin{align*}
a(\F)(\mathbf{y}_i(y_0),\phi)=0 \quad \forall \phi \in W_\infty,~y_0 \in Y_0 \setminus O_i.
\end{align*}
From Proposition~\ref{prop:equivalofweakandstrong} we thus conclude~$\mathbf{y}_1(y_0)=\mathbf{y}_2(y_0)$ for all~$y_0 \in Y_0 \setminus O_1 \cup O_2$. Since we have~$\mu (O_1 \cup O_2)=0$ this yields~$\mathbf{y}_1= \mathbf{y}_2$ in~$\Linfboch$.

The last statement follows immediately from the definition of the form~$\abd(\cdot)(\cdot,\cdot)$.
\end{proof}
\begin{example}
Let~$(\F^*, \mathbf{y}^*)$ denote the pair of optimal feedback law and  solution mapping to the associated closed loop system according to\tcr{~$\mathbf{A.2}$ and~$\mathbf{A.3}$. Due to~$\mathbf{A.3}$ and~$\mathbf{A.4}$ }we have~$\mathbf{y}^* \in \mathcal{C}_b(Y_0;W_\infty)$ and~$\F^* \in \mathcal{H}$. Thus~$\mathbf{y}^*$ is~$\mu$-measurable and there holds~$\mathbf{y}^* \in \Linfboch$ with~$\Linfbochnorm{\mathbf{y}^*}\leq M_{Y_0}$. Moreover for every~$y_0 \in Y_0$ we have
\begin{align*}
a(\F^*)(\mathbf{y}^*(y_0),y_0,\phi)=0 \quad \forall \phi \in W_\infty,
\end{align*}
see Proposition~\ref{prop:equivalofweakandstrong}.
In particular, for an arbitrary~$\Phi \in \Lzwoboch$, this implies that
\begin{align*}
a(\F^*)(\mathbf{y}^*(y_0),y_0,\Phi(y_0))=0 \quad \text{for~$\mu$-a.e.}~y_0 \in Y_0,
\end{align*}
i.e.~$\mathbf{y}^*$ is the unique ensemble solution to~\eqref{eq:parametrizedstate} given~$\F^*$. $\square$
\end{example}

We are now prepared to introduce the optimization problem which
constitutes the learning of the feedback function  ${\mathcal {F}}$.
This will involve the ensemble of initial conditions $Y_0$.

For this purpose we define  the ensemble objective functional
\begin{align*}
j \colon \Linfboch \times \mathcal{H} \to \R, \quad j( \mathbf{y},\F)= \int_{Y_0} J(\mathbf{y}(y_0), \F(\mathbf{y}(y_0)))~\de \mu(y_0),
\end{align*}
which is understood as an extended real valued functional,
and the associated constrained minimization problem

\begin{equation} \label{def:contproblem}\tag{$\mathcal{P}$}
\left\{
\begin{aligned}
\quad
&\min_{\substack{\F \in \mathcal{H},
\mathbf{y} \in \Linfboch}} j(\mathbf{y},\F) \\
&s.t.\quad \abd(\mathcal{F})(\mathbf{y}, \Phi) =0 \quad \forall \Phi \in \Lzwoboch,~\Linfbochnorm{\mathbf{y}}\leq 2 M_{Y_0}.
\end{aligned}
\right.
\end{equation}

Problem \eqref{def:contproblem} together with a particular choice for the approximation of the elements $\F \in \mathcal{H}$  constitutes the approach that we propose to determine  semiglobal optimal feedback laws.

By Lemma~\ref{lem:Nemitsky} the cost $j$ is finite on the admissible set.
We next address the existence of a solution to~\eqref{def:contproblem}, and  verify that  every minimizing pair to~\eqref{def:contproblem} solves the semiglobal optimal feedback stabilization problem for~\tcr{$\mu$-}a.e. $y_0$.
\begin{prop} \label{prop:optimalFstar}
Let Assumption~\ref{ass:feedbacklaw} hold. Then~the pair~$(\F^*,\mathbf{y}^*)\in \mathcal{H}\times \Linfboch$ is a global minimizer of problem~\eqref{def:contproblem}. Vice versa, if~$(\bar{\F}, \bar{\mathbf{y}})\in \mathcal{H} \times \Linfboch$ is an optimal solution of~\eqref{def:contproblem}, then there holds
\begin{align*}
(\bar{\mathbf{y}}(y_0), \bar{\F}(\bar{\mathbf{y}}(y_0)) ) \in \argmin \eqref{def:openloopproblem}
\end{align*}
for~$\mu$-a.e.~$y_0 \in Y_0$.
\end{prop}
\begin{proof}
By~$\mathbf{A.2}$ we have
\begin{align*}
J(\mathbf{y}^*(y_0),\F^*(\mathbf{y}^*(y_0)))= V(y_0)\quad \forall y_0 \in Y_0.
\end{align*}
The optimality of~$(\mathbf{y}^*,\F^*)$ thus follows directly from
\begin{align} \label{eq:optestimate}
j(\mathbf{y}, \F)= \int_{Y_0} J(\mathbf{y}(y_0),\F(\mathbf{y}(y_0)))~\de \mu(y_0) \geq \int_{Y_0} V(y_0)~\de \mu(y_0)=j(\mathbf{y}^*, \F^*)
\end{align}
for every admissible pair~$(\mathbf{y},\F)\in \Linfboch \times \mathcal{H}$.

The second statement is obtained in a similar way. If~$(\bar{\F}, \bar{\mathbf{y}})$ is an arbitrary but fixed minimizing pair of~\eqref{def:contproblem}, then we immediately get
\begin{align*}
\int_{Y_0} \left \lbrack J(\bar{\mathbf{y}}(y_0), \bar{\F}(\bar{\mathbf{y}}(y_0)) )-V(y_0) \right \rbrack~\de \mu(y_0)=0
\end{align*}
as well as
\begin{align*}
J(\bar{\mathbf{y}}(y_0), \bar{\F}(\bar{\mathbf{y}}(y_0)) ) \geq V(y_0) \quad \text{for}~y_0 \in Y_0~\mu-a.e.
\end{align*}
from~\eqref{eq:optestimate}, Proposition~\ref{prop:parametrizedstate} and the definition of the value function. Thus, there holds
\begin{align*}
J(\bar{\mathbf{y}}(y_0),\bar{\F}(\bar{\mathbf{y}}(y_0)))=V(y_0) \quad \text{for}~y_0 \in Y_0~\mu-a.e.
\end{align*}
This yields the optimality of~$(\bar{\mathbf{y}}(y_0),\bar{\F}(\bar{\mathbf{y}}(y_0)))$ for~\eqref{def:openloopproblem} and~$\mu$-almost every~$y_0 \in Y_0$.
\end{proof}

For the following considerations we need to recall the notion of support for the measure $\mu$:
\begin{align*}
\supp \mu = Y_0 \setminus  \bigcup  \left \{\,O \in \mathcal{A} \;|\;\mu(O)=0,~O~\text{is open}\,\right\},
\end{align*}
where we assume that~$ \mathcal{A}$ contains the Borel~$\sigma$-algebra on~$Y_0$. Then the support of~$\mu$
is compact. In fact,  by assumption,~$\mu$ is a regular Borel measure. Hence its support is closed, see e.g. \cite[pg. 252-254]{Ro63}.  Since~$\supp \mu \subset Y_0$ and~$Y_0$ is bounded, the compactness of~$\supp \mu$ follows.

Let us note that the explicit form of $j$ arising in \eqref{def:contproblem} is given by
$$
j(y,\mathcal{F})= \frac{1}{2}\int_{Y_0} \, \int_I (|Q y(y_0)(t)|^2 +
\beta|F(y(y_0)(t))|^2 dt \, d\mu(y_0).
$$
Thus for each $y_0\in \supp \mu$ the feedback function is learned along all the
trajectory $\{y(y_0)(t): t\in I\}$. In view of the  Bellman principle
this can also be interpreted in such a manner  that for each $y_0\in \supp \mu$
separately, the cost described in the training of  $\mathcal{F}$ by means of the functional $j$  is influenced by the whole trajectory  $\{y(y_0)(t): t\in I\}$. The following
corollary expresses these observations in a formal manner, and shows how the optimal solutions to \eqref{def:contproblem} naturally can be extended to the set of all points which are met by optimal trajectories starting in $Y_0$.

\begin{coroll} \label{cor:largerset}
Let~$(\bar{\F}, \bar{\mathbf{y}})$ be an optimal solution to~\eqref{def:contproblem} and assume that~$ \mathcal{A}$ contains the Borel~$\sigma$-algebra on~$Y_0$. If~$\bar{\mathbf{y}} \in \mathcal{C}_b(\supp \mu; W_\infty)$ the set of trajectories
\begin{align*}
\widehat{Y}_0:= {\left \{\, \bar{\mathbf{y}}(y_0)(t) \ \;|\;y_0 \in \supp \mu,~t \in [0, +\infty)\,\right\}} \subset \bar{B}_{2 \widehat{M}}(0)
\end{align*}
\tcr{is compact}. Furthermore there exists a unique function
\begin{align*}
\widehat{\mathbf{y}}\in \mathcal{C}_b(\widehat{Y}_0;W_\infty), \quad \|\widehat{\mathbf{y}}(y_0)\|_{\mathcal{C}_b(\widehat{Y}_0;W_\infty)}\leq 2M_{Y_0}
\end{align*}
which fulfills
\begin{align} \label{eq:optimalityextended}
\tcr{a(\bar{\F})}(\widehat{\mathbf{y}}(y_0),y_0, \phi)=0 \quad \forall \phi \in W_\infty, \quad (\widehat{\mathbf{y}}(y_0),\bar{\F}(\widehat{\mathbf{y}}(y_0))) \in \argmin \eqref{def:openloopproblem}
\end{align}
for every~$y_0 \in \widehat{Y}_0$.

\end{coroll}
\begin{proof}
We first show the compactness of~$\widehat{Y}_0$. Let an arbitrary sequence~$\{y_k\}_{k\in\N} \subset \widehat{Y}_0$ be given. Then there exists~$\{y^k_0,t_k\}_{k\in\N} \subset \supp \mu \times I$ such that~$y_k=\bar{\mathbf{y}}(y^k_0)(t_k)$,~$k\in \N$. From the compactness of~$\supp \mu$ we further get the existence of a subsequence, denoted by the same symbol, as well as an element~$\bar{y}_0 \in \supp \mu$ with~$y^k_0 \rightarrow \bar{y}_0$.

We distinguish two cases. First assume that~$\{t_k\}_{k\in \N}$ is bounded. Then, possibly again by selecting a subsequence, there exists~$\bar{t} \in I$ with~$(y^k_0, t_k) \rightarrow (\bar{y}_0, \bar{t})$
as~$k \rightarrow \infty$. Thus we also conclude~$\bar{\mathbf{y}}(y^k_0)(t_k) \rightarrow \bar{\mathbf{y}}(\bar{y}_0)(\bar{t})\in \widehat{Y}_0$ in~$\R^n$.
Second, if~$\{t_k\}_{k\in\N}$ is unbounded, we claim~$\bar{\mathbf{y}}(y^k_0)(t_k)\rightarrow 0 \in \widehat{Y}_0$. By
\begin{align*}
|\bar{\mathbf{y}}(y^k_0)(t_k)| &\leq |\bar{\mathbf{y}}(y^k_0)(t_k)-\bar{\mathbf{y}}(\bar{y}_0)(t_k)|+|\bar{\mathbf{y}}(\bar{y}_0)(t_k)| \\&\leq \|\bar{\mathbf{y}}(y^k_0)-\bar{\mathbf{y}}(\bar{y}_0)\|_{\mathcal{C}_b(I;\R^n)}+|\bar{\mathbf{y}}(\bar{y}_0)(t_k)|.
\end{align*}
The first term on the right hand side vanishes due to~$\bar{\mathbf{y}}\in \mathcal{C}_b(\supp \mu;\mathcal{C}_b(I;\R^n))$ while the second term goes to~$0$ since~$\bar{\mathbf{y}}(\bar{y}_0)\in W_{\infty}$.
Since~$\{y_k\}_{k\in\N} \subset \widehat{Y}_0$ was chosen arbitrarily we conclude that every sequence in~$\widehat{Y}_0$ admits a convergent subsequence whose limit is again in~$\widehat{Y}_0$. Thus,~$\widehat{Y}_0$ is compact in~$\R^n$.

Let~$\widehat{y}_0 \in \widehat{Y}_0 $ be given. Then we have~$\widehat{y}_0=\bar{\mathbf{y}}(\bar{y}_0)(\bar{t})$ for some~$\bar{y}_0 \in \supp \mu$ and~$\bar{t}\in I$. Define the function~$\widehat{\mathbf{y}}(\widehat{y}_0)= \bar{\mathbf{y}}(\bar{y}_0)(\cdot+\bar{t}) $. Clearly, there holds~$\widehat{\mathbf{y}}(\widehat{y}_0) \in \mathcal{Y}_{ad}$. Moreover~$\widehat{\mathbf{y}}(\widehat{y}_0)$ is the unique solution in~$\mathcal{Y}_{ad}$ to
\begin{align*}
\dot{y}=\mathbf{f}(y) +\mathcal{B} \bar{\F}(y), \quad y(0)=\widehat{y}_0,
\end{align*}
and the mapping
\begin{align*}
\widehat{\mathbf{y}} \colon \widehat{Y}_0 \to  \tcr{\Ya} \subset W_\infty, \quad y_0 \mapsto \widehat{\mathbf{y}}(y_0)
\end{align*}
is continuous. By Bellman's dynamical programming principle, we conclude that
\begin{align*}
(\widehat{\mathbf{y}}(y_0),\bar{\F}(\widehat{\mathbf{y}}(y_0))) \in \argmin (P^{y_0}_\beta)
\end{align*}
for every~$y_0 \in \widehat{Y}_0$ and thus~\eqref{eq:optimalityextended} holds.

\end{proof}
\begin{example} \label{exp:possiblemeasures}
To close this section we briefly discuss two canonical examples for the measure space~$(Y_0, \mathcal{A},\mu)$. In both cases we assume that~$Y_0 \subset \R^n$ is the closure of a non-empty domain containing~$0$.

\textbf{a)} As a first example, choose~$\mathcal{A}$ as the Lebesgue $\sigma$-algebra and~$\mu= \lambda(\,\cdot \,\cap Y_0)/\lambda(Y_0)$, where~$\lambda$ denotes the Lebesgue measure on~$\R^n$. Then~$(Y_0, \mathcal{A},\mu)$ is complete and~$\supp \mu=Y_0$. Consequently, any minimizing pair~$(\bar{\F}, \bar{\mathbf{y}})\in \mathcal{H}\times L^\infty_\mu(Y_0; W_\infty)$ to~\eqref{def:contproblem} fulfills
\begin{align*}
(\bar{\mathbf{y}}(y_0), \bar{\F}(\bar{\mathbf{y}}(y_0))) \in \argmin \eqref{def:openloopproblem}
\end{align*}
for all~$y_0 \in Y_0$ with the exception of a Lebesgue zero set.
If the prerequisites of Corollary~\ref{cor:largerset} are met, that is if~$\bar{\mathbf{y}}\in \mathcal{C}_b(Y_0, W_\infty)$, we can extend~$\bar{\mathbf{y}}$ to a continuous function~$\widehat{\mathbf{y}}\in \mathcal{C}(\widehat{Y}_0;W_\infty)$. The pair~$(\widehat{\mathbf{y}}(y_0), \bar{\F}(\widehat{\mathbf{y}}))$ provides an optimal solution to~\eqref{def:openloopproblem} for every~$y_0 \in \widehat{Y}_0$. Note that~$Y_0 \subset \widehat{Y}_0$ where the inclusion can be strict since the optimal trajectories~$\bar{\mathbf{y}}(y_0)$,~$y_0 \in Y_0$, possibly leave the set~$Y_0$.

\textbf{b)} Next we consider a measure~$\mu$ given by a convex combination of finitely many Dirac delta functions on~$Y_0$ i.e.~$\mu= \sum^N_{i=1} \alpha_i \delta_{y^i_0} $,~$\alpha_i >0$,~$y^i_0 \in \R^n$ as well as~$\sum^N_{i=1}\alpha_i=1$. This example is of particular importance for the numerical realization of~\eqref{def:contproblem}, see Section \ref{sec:practicalrealization}. We choose~$\mathcal{A}$ as the completion of the Borel $\sigma$-algebra on~$Y_0$ with respect to~$\mu$. Then we have~$\supp \mu= \{y^i_0\}^N_{i=1}$. We point out that ~$L^2_\mu(Y_0; W_\infty)\simeq W^N_\infty$ as well as~$L^\infty_\mu(Y_0; W_\infty)\simeq \mathcal{C}_b(\supp \mu;W_\infty)$. In particular, given a feedback law~$\F \in \mathcal{H}$, a function~$\mathbf{y}\in L^\infty_\mu (Y_0;W_\infty)$ is an ensemble solution to~\eqref{eq:parametrizedstate} if and only if the vector~$(\mathbf{y}(y^1_0),\dots,\mathbf{y}(y^N_0)) \in W^N_\infty$ fulfills~$\wnorm{\mathbf{y}(y^i_0)}\leq 2M_{Y_0}$ and
\begin{align*}
a(\F)(\mathbf{y}(y^i_0),y^i_0,\phi)=0 \quad \forall \phi \in W_\infty
\end{align*}
for~$i=1,\dots,N$. Thus, problem~\eqref{def:contproblem} can be equivalently reformulated as
\begin{align*}
\min_{\substack{\F \in \mathcal{H}, \\ y_i \in \mathcal{Y}_{ad},~i=1,\dots, N}} \sum^N_{i=1} \frac{\alpha_i}{2} \left \lbrack \int^\infty_0 \left( |Qy_i(t)|^2+ \beta |\F(y_i(t))|^2\right )~\de t \right\rbrack
\end{align*}
subject to
\begin{align*}
a(\F)(y_i,y^i_0, \phi)=0 \quad \forall \phi \in W_\infty,~i=1, \dots,N.
\end{align*}
Let an optimal feedback pair~$(\bar{\F}, \bar{\mathbf{y}})$ obtained by solving~\eqref{def:contproblem} be given. Due to the dynamic programming principle, see Corollary~\ref{cor:largerset}, the associated closed loop system~\eqref{eq:closedloopweak} admits a unique solution~$\widehat{\mathbf{y}}(y_0)$ for every~$y_0\in \widehat{Y}_0 = \bigcup^N_{i=1}\bigcup_{t \in I} \bar{\mathbf{y}}(y^i_0)(t)$ which further fulfills
\begin{align*}
(\widehat{\mathbf{y}}(y_0),\bar{\F}(\widehat{\mathbf{y}}(y_0))) \in \argmin \eqref{def:openloopproblem}.
\end{align*}
\end{example}

\section{Pertinent results on neural networks} \label{sec:neuralnet}
In order to solve the semiglobal feedback stabilization problem~\eqref{def:contproblem}, we replace the set of admissible feedback laws~$\mathcal{H}$ by a family of finitely parameterized ones. For this purpose we use neural networks. In this section we summarize the concepts which are relevant for the present  paper.

\subsection{Notation and network structure}

Let~$L\in \N$,~$L\geq 2$, as well as~$N_{i}\in \N$,~$i=1, \dots,L-1$ be given. We  set~$N_0=n$ and~$N_L=m$. Furthermore define
\begin{align*}
\mathcal{R}= \bigtimes^L_{i=1} \left ( \R^{N_{i} \times N_{i-1}} \times \R^{N_i} \right ).
\end{align*}
Note that the space~$\mathcal{R}$ is uniquely determined by its~\textit{architecture}
\begin{align*}
\text{arch}(\mathcal{R})=\left( N_0, N_1, \dots , N_L\right)\in \N^{L+1}.
\end{align*}
An L-tupel of parameters~$\theta \in \mathcal{R}$ given by
\begin{align*}
\theta=\left( W_{1}, b_1,\dots, W_{L}, b_L \right)
\end{align*}
is called a~\textit{neural network} with~$L$~\textit{layers}. We equip the space~$\mathcal{R}$ with the canonical product norm
\begin{align*}
\|\theta\|_{\mathcal{R}}= \sqrt{\sum^L_{i=1} \lbrack \|W_i\|^2 + |b_i|^2 \rbrack} \quad \forall \theta \in \mathcal{R}
\end{align*}
where~$\|\cdot\|$ denotes the Frobenius norm.
Moreover let~$\sigma \in \mathcal{C}(\R)$ with~$\sigma(0)=0$ be given. A function~$F^{\sigma}_\theta \colon \R^n \to \R^m$ is called the~\textit{realization} of~$\theta$ with~\textit{activation function}~$\sigma$ if
\begin{align}\label{eq:shift}
F^\sigma_\theta(x)= f^{\sigma}_{L, \theta}  \circ f^{\sigma}_{L-1, \theta} \circ \cdots \circ f^{\sigma}_{1, \theta}(x)-f^{\sigma}_{L, \theta}  \circ f^{\sigma}_{L-1, \theta} \circ \cdots \circ f^{\sigma}_{1, \theta}(0) \quad \forall x \in \R^n
\end{align}
where
\begin{align*}
f^{\sigma}_{L, \theta}(x)= W_{L} x+ b_L \quad \forall x \in \R^{N_{L-1}}
\end{align*}
as well as
\begin{align*}
f^{\sigma}_{i, \theta}(x)= \sigma(W_{i}x+b_i) \quad  \forall x \in \R^{N_{i-1}},~i=1, \dots, L-1 .
\end{align*}
Here the application of~$\sigma$ should be understood componentwise i.e. given an index~$i \in \{1,\dots,L-1\}$ and~$x \in \R^{N_i}$ we set
\begin{align*}
\sigma(x)=(\sigma(x_1), \dots, \sigma(x_{N_i}))^t.
\end{align*}
By construction we have that $F^\sigma_\theta(0)=0$. Subsequently, we define the set of associated neural network feedback laws as
\begin{align*}
\mathcal{H}^\sigma_\mathcal{R}= \left\{\,\F \colon \mathcal{Y}_{ad} \to L^2(I;\R^m)\;|\; \F(y)(t)=F^\sigma_\theta(y(t)),~\theta \in \mathcal{R}\,\right\}.
\end{align*}
The following standing assumption is made throughout this paper.
\begin{assumption} \label{ass:smoothofact}
There holds~$\sigma \in \mathcal{C}^1(\R)$.
\end{assumption}
The next proposition is imminent.
\begin{prop} \label{prop:conformity}
Let Assumption~\ref{ass:smoothofact} be satisfied. Then there holds~$\mathcal{H}^\sigma_\mathcal{R} \subset \mathcal{H}$.
\end{prop}
\begin{proof}
Let an arbitrary~$\theta \in \mathcal{R}$ be given. Denote by~$F^\sigma_\theta$ the realization of~$\theta$ with activation function~$\sigma \in \mathcal{C}^1(\R)$. We already observed that ~$F^\sigma_\theta(0)=0$. Next note that every layer~$s$ is either affine linear or the composition of the vectorized function~$\sigma$ and a continuous affine linear mapping. \tcr{In particular, due to the smoothness of~$\sigma$, we conclude the Lipschitz continuity} of~$f^\sigma_{i,\theta}$,~$i=1, \dots,L$, on compact sets. Since compactness is preserved under continuous mappings it is straightforward to verify that~$F^\sigma_\theta=f^\sigma_{L, \theta} \circ \cdots \circ f^\sigma_{1, \theta}$ is also Lipschitz continuous on compact sets. Thus, in particular, we have~$F^\sigma_\theta \in \operatorname{Lip}(\bar{B}_{2 \widehat{M}}(0), \R^m)$. This finishes the proof.
\end{proof}

As a consequence elements in $\mathcal{H}^\sigma_\mathcal{R}$ have the same Lipschitz properties as those in  $\mathcal{H}$.

\subsection{A neural network density result in $\mathcal{C}^1$} \label{subsec:density}

The aim of this subsection is to ascertain the following approximation result for the optimal feedback law~$\mathcal{F}^*$ appearing in Assumption~\ref{ass:feedbacklaw} by realizations of bounded neural networks. For $v\in \mathbb{R}^n$ we denote by $|v|_{\infty} = sup\{|v_i|: i = 1,\ldots,n\}$, and $\|W\|_{\infty}$ is the matrix norm subordinate to $\mathbb{R}^n$ and $\mathbb{R}^m$ endowed with the $|\cdot|_{\infty}$-norms.
\begin{theorem}\label{thm:network}
Let Assumptions~\ref{ass:feedbacklaw} and~\ref{ass:smoothofact} hold.
Let $ \eta_1 > 0,\eta_2 > 0$, and assume that the activation function~$\sigma$ is not a polynomial. Then for each $\varepsilon > 0$ and every~$L_\eps \in\N$,~there exist~$\operatorname{arch}(\mathcal{R}_\eps)\in \N^{L_\eps +1}$ and a neural network
\begin{align*}
\theta_\eps = (W^\eps_1, b^\eps_1, \dots, W^\eps_{L_\eps}, b^\eps_{L_\eps}) \in \mathcal{R}_\eps
\end{align*}
such that~$\|W^\eps_1\|_\infty \leq \eta_1$,~$|b^\eps_i|_\infty \leq \eta_2$,~$i=1,\dots,L$, as well as
\begin{align*}
|F^*(x)-F^\sigma_{\theta_\eps}(x)|+\|DF^*(x)-DF^\sigma_{\theta_\eps}(x)\|
\leq \eps
\end{align*}
for all~$|x|\leq 2M_{Y_0}$.
\end{theorem}

\tcr{If the activation function~$\sigma$ were a polynomial of degree~$N$, then each component of~$F^\sigma_{\theta_\eps}$ would be a polynomial of at most degree~$L_\eps N$, and it would not be possible to approximate an arbitrary function with accuracy $\epsilon$ in $C^1$.}
The proof of Theorem \ref{thm:network} depends on several preparatory steps and will be given at the end of this section. For $\psi:\mathbb{R}\to \mathbb{R}$ we define~\tcr{the set~$\Sigma_n$ of mappings from~$\R^n$ to~$\R$} as
$$
\Sigma_{n} = span \{\, \Psi \;|\; \Psi(x)= \psi( w \cdot x + b) : w \in \mathcal{W}\setminus \{0\},b \in B \, \}
$$
where $\mathcal{W} \subset \mathbb{R}^n$ contains~\tcr{a neighborhood of 0}, and $B\subset \mathbb{R}$ contains an open interval $B_0=( \underline{b}, \bar b)$.
Below $\mathcal{C}^1(\mathbb{R}^n,\mathbb{R})$ is endowed with the topology of uniform convergence of its elements and its derivatives on compact subsets of $\mathbb{R}^n$. Further we call a set of functions $F\subset W_{loc}^{1,\infty}(\mathbb{R}^n, \R)$ {\em{ dense in $\mathcal{C}^1(\mathbb{R}^n, \R)$}} if for every $f \in \mathcal{C}^1(\R^n,\R)$ and every compact set $K\subset \R^n$ there exists a sequence of functions $\{f_n\}\subset W^{1,\infty}(\mathbb{R}^n, \R)$ such that
 $$
 \lim_{n\to \infty} \|f_n-f\|_{W^{1,\infty}(K, \R)}=0.
 $$
 We note that the approximation of  $\mathcal{C}^1$-regular functions $f$
  needs to be  achieved by functions which are possibly only $W_{loc}^{1,\infty}(\mathbb{R}^n, \R)$-
  regular.

\begin{prop}\label{C1app}
If $\psi \in W^{1,\infty}(\mathbb{R})$ is not a polynomial then $\Sigma_{n}$ is dense in $\mathcal{C}^1(\mathbb{R}^n,\mathbb{R})$.
\end{prop}
Since the proof of this universal approximation property is rather technical we postpone it to Appendix~\ref{app:universal}. Proposition \ref{C1app} allows to derive the following approximation result for neural networks.
\begin{prop}\label{prop:network}
Let $L \ge 2, \eta_1 > 0,\eta_2 > 0$, and let the assumptions of Proposition \ref{C1app} hold. Then for each $\epsilon > 0, K \subset \mathbb{R}^n$ compact, and $f \in C^1(K,\mathbb{R}^m)$ there exist $\{N_i\}_{i=1}^{L-1}$ and a neural network $\theta = (W_1, b_1, \dots, b_{L-1}, W_L, 0) \in \mathcal{R}$ with $arch({\mathcal R}) = (n, N_1,\dots, N_{L-1}, m)$, and \tcr{$\|W_1\|_{\infty} \le \eta_1, |b_i|_{\infty} \le \eta_2$, $i=1,\dots,L-1$},  such that
\begin{equation} \label{estim1}
\|f - f_{L,\theta}^\sigma \circ \dots \circ f_{2,\theta}^\sigma \circ f_{1,\theta}^\sigma \|_{\mathcal{C}^1(K,\mathbb{R}^m)} < \epsilon.
\end{equation}
\end{prop}

\begin{proof}

{\em Step 1.} Here we treat the case $L=2$.
By Proposition \ref{C1app} there exist dimensions $(k_1, \ldots, k_m)$ and $W^j_1 \in \mathbb{R}^{k_j\times n}, W^j_2 \in \mathbb{R}^{1\times k_j}, b^j_1 \in \mathbb{R}^{k_j}$ with $\|W_1^j\|_{\infty} \le \eta_1, |b_1^j|_{\infty} \le \eta_2, j = 1, \ldots, m,$ such that
\begin{equation} \label{Kaux1}
\|W_2^j ~\sigma (W_1^j \cdot +\, b_1^j) - f_j\|_{\mathcal{C}^1(K,\mathbb{R}^1)} < \frac{\epsilon}{\sqrt{m}}.
\end{equation}
In fact, for $j=1$, for example, we obtain from Proposition \ref{C1app} the existence of $m_1$ and $(W_1^1)_i^t \in \mathbb{R}^n, (b^1)_i$, where $i=1,\dots, m_1,$  and $(W_2^1)^t \in \mathbb{R}^{m_1}$ with $\|(W_1^1)_i\|_{\infty} \le \eta_1, |(b^1)_i| \le \eta_2$ such that
\begin{equation}
\|W_2^1 ~\sigma (W_1^1 \cdot + \, b_1) - f_1\|_{C^1(K,\mathbb{R})} < \frac{\epsilon}{\sqrt{m}},
\end{equation}
and $W^1_1=col((W_1^1)_1, \dots,(W_1^1)_{m_1})$.
Since this holds for all coordinates of the vector-valued function $f$ we have \eqref{Kaux1}.\\

To treat the vector-valued case we define
\begin{equation*}
\begin{array}l
N_1 = \Pi_{j =1}^m k_j, \qquad  \qquad  \qquad \quad \; \; W_1 = col (W_1^1, \ldots, W_1^m) \in \mathbb{R}^{N_1 \times n}, \\[1.7ex]
 b_1 = col (b_1^1, \ldots, b_1^m) \in \mathbb{R}^{N_1}, \qquad  W_2 = col (\tilde{W}_2^1, \ldots, \tilde{W}_2^m) \in \mathbb{R}^{m \times N_1}, \\[1.7ex]
  \tilde{W}_2^j = ((O^j_L)^t, (W_2^j)_1, \ldots,  (W_2^j)_{k_j}, (O^j_R)^t),
\end{array}
\end{equation*}
where $O^j_L$ is the vector of zeros in $\mathbb{R}^{\Pi_{i = 1}^{j-1} k_i}, O^j_R$ is the vector of zeros in $\mathbb{R}^{\Pi_{i = j+1}^m{k_i}}$, and $\Pi_{i = 1}^0 k_i = \Pi_{i = m+1}^m k_i$ are empty sets.
By construction \eqref{estim1} follows from \eqref{Kaux1}.

{\em Step 2.} We proceed by induction and suppose that the result has been verified for all  levels up to $L-1$ with $L\ge 3$. The proof is inspired by \cite[Corollary 2.6]{Ho89}. By induction hypothesis there exists for arbitrary $\epsilon$ and each compact subset $K$ in $\mathbb{R}^n$, a space ${\mathcal R}_{L-1}$   with $arch({\mathcal R}_{L-1}) = (n, N_1,\dots, N_{L-2},m)$ and a network  $\theta_{L-1}=( W_1, b_1, \dots, b_{L-2}, \hat W_{L-1}, 0) \in \mathcal{R}_{L-1}$, with
$\|W_1\|_{\infty} \le \eta_1, |b_i|_{\infty} \le \eta_2$, $i=1,\dots,L-2$,  such that
\begin{equation} \label{estim2}
\|f - \hat f_{L-1,\theta}^\sigma \circ f_{L-2,\theta}^\sigma \circ\dots \circ f_{2,\theta}^\sigma \circ f_{1,\theta}^\sigma \|_{\mathcal{C}^1(K,\mathbb{R}^m)} < \frac{\epsilon}{2},
\end{equation}
where $\hat f_{L-1,\theta}^\sigma (x)= \hat W_{L-1}x$,  $f_{i,\theta}^\sigma (x)=  \sigma (W_{i}x +b_i)$ for $i=1,\dots,L-2$. Without loss of generality we may assume that the norm on $\mathcal{C}^1(K, \mathbb{R}^m)$ is chosen as $\| f\|_{\mathcal{C}^1(K,\mathbb{R}^m)}=  \| f\|_{\mathcal{C}(K,\mathbb{R}^m)} + \| f'\|_{\mathcal{C}(K,\mathbb{R}^m)}$.

Let us set $\varphi_1= \hat f_{L-1,\theta}^\sigma \circ f_{L-2,\theta}^\sigma \circ\dots \circ f_{2,\theta}^\sigma \circ f_{1,\theta}^\sigma$ and $\hat K = \overline {\varphi_1 (K)} \cup \overline {\varphi'_1 (K)}$, which is a compact set in $\mathbb{R}^m$. By induction hypothesis there exist $N_{L-1}$, $W_L \in \mathbb{R}^{m\times N_{L-1}},
\tilde W_{L-1} \in \mathbb{R}^{ N_{L-1}\times m}$, $b_{L-1}\in  \mathbb{R}^{ N_{L-1}}$ with $|b_{L-1}|_\infty \le \eta_2$ such that
\begin{equation} \label{estim3}
\|Id -  f_{L,\theta}^\sigma \circ \tilde f_{L-2,\theta}^\sigma \|_{\mathcal{C}^1(\hat K,\mathbb{R}^m)} < \frac{\epsilon}{6},
\end{equation}
where $f_{L,\theta}^\sigma (x)=  W_{L}x$,  $\tilde f_{L-1,\theta}^\sigma (x)=  \sigma (\tilde W_{L-1}x +b_{L-1})$ and we set $\varphi_2= f_{L,\theta}^\sigma \circ \tilde f_{L-1,\theta}^\sigma$.

Next we concatenate the networks $\varphi_1$ and $\varphi_2$ to a network $\varphi= \varphi_2 \circ \varphi_1$ with $arch(\mathcal{R})= (n,N_1,N_{L-1}, m)$. By construction its elements satisfy the dimensions and bounds specified in the theorem and $W_{L-2} =  \tilde W_{L-2} \hat W_{L-2} \in \mathbb{R}^{N_{L-2} \times N_{L-1}}$. By \eqref{estim2} and \eqref{estim3} we have the estimate
\begin{align*}
\|f - &\varphi\|_{\mathcal{C}^1(K,\mathbb{R}^m)} < \|f - \varphi_1\|_{\mathcal{C}^1(K,\mathbb{R}^m)} + \|(Id -\varphi_2)\circ  \varphi_1 \|_{\mathcal{C}^1(K,\mathbb{R}^m)} \\
&< \frac{\epsilon}{2} + \|Id -\varphi_2\|_{\mathcal{C}(\hat K,\mathbb{R}^m)} +  \|(Id -\varphi_2) \varphi_1'\|_{\mathcal{C}(K,\mathbb{R}^m)} + \|\varphi_2' \varphi_1\|_{\mathcal{C}( K,\mathbb{R}^m)} <\epsilon,
\end{align*}
as desired.
\end{proof}
\begin{proof}
[Proof of Theorem~\ref{thm:network}]
Fix~$\eps>0$ and choose an arbitrary but fixed~$L_\eps \in \N$,~$L_\eps \geq 2$. According to Proposition~\ref{prop:network} there exist an architecture~$\operatorname{arch}(\mathcal{R}_\eps)\in \N^{L_\eps+1} $ as well as a neural network
\begin{align*}
\theta_\eps = (W^\eps_1, b^\eps_1, \dots, W^\eps_L, b^\eps_L) \in \mathcal{R}_\eps
\end{align*}
satisfying~$\|W^\eps_1\|_\infty \leq \eta_1$,~$|b^\eps_i|_\infty\leq \eta_2$,~$i=1,\dots,L_\eps$, as well as
\begin{align*}
\|F^* - f_{L_\eps,\theta_\eps}^\sigma \circ \dots \circ f_{2,\theta_\eps}^\sigma \circ f_{1,\theta_\eps}^\sigma \|_{\mathcal{C}^1(\bar{B}_{2\widehat{M}}(0),\mathbb{R}^m)} < \frac{\eps}{3}.
\end{align*}
In particular this implies
\begin{align*}
\|DF^*(x)-DF^\sigma_{\theta_\eps}(x)\|< \frac{\eps}{3}
\end{align*}
and, utilizing~$F^*(0)=0$ and \eqref{eq:shift},
\begin{align*}
|F^*(x)-F^\sigma_{\theta_\eps}(x)| &\leq |F^*(0)- f_{L_\eps,\theta_\eps}^\sigma \circ \dots \circ f_{2,\theta_\eps}^\sigma \circ f_{1,\theta_\eps}^\sigma(0)| \\
&+|F^*(x)- f_{L_\eps,\theta_\eps}^\sigma \circ \dots \circ f_{2,\theta_\eps}^\sigma \circ f_{1,\theta_\eps}^\sigma(x)|
< \frac{\eps}{3}+\frac{\eps}{3}=  \frac{2\eps}{3}
\end{align*}
for all~$x\in \bar{B}_{2\widehat{M}}(0)$. This finishes the proof.
\end{proof}
\section{Approximation of the learning problem by neural networks: Well-posedness \& Convergence} \label{sec:wellposed}

We commence this section by defining neural network based approximations to  ~\eqref{def:contproblem}. Let us observe that
by construction the feedback laws in~$\mathcal{H}^\sigma_\mathcal{R}$ can be parameterized as a function of the associated finite-dimensional neural network~$\theta \in \mathcal{R}$.
For this purpose we introduce the surjective mapping
\begin{align*}
\F^\sigma_\cdot \colon \mathcal{R} \to \mathcal{H}^\sigma_\mathcal{R} \quad \text{where} \quad \mathcal{F}^\sigma_\theta(y)(t)=F^\sigma_\theta(y(t))
\end{align*}
with~$y \in W_\infty$ and~$t\in I$. Fixing constants~$\eta_1,~\eta_2>0$, and choosing a number of layers~$L\in \N$ as well as an architecture~$\operatorname{arch}(\mathcal{R})\in \N^{L+1}$,
 we propose to approximate~\eqref{def:contproblem} by
\begin{equation}\label{def:neuralnetproblem}\tag{$\mathcal{P}^\sigma_\mathcal{R}$}
\left\{
\begin{aligned}
\quad & \min_{\substack{\theta \in \mathcal{R}, \\ \mathbf{y} \in \Linfboch}} j(\mathbf{y},\F^\sigma_\theta)+ \mathcal{G}_\mathcal{R}(\theta)
\\
&s.t. \quad  \abd(\F^\sigma_\theta)(\mathbf{y}, \Phi) =0 \quad \forall \Phi \in \Lzwoboch,~\Linfbochnorm{\mathbf{y}}\leq 2 M_{Y_0}.
\end{aligned}
\right.
\end{equation}
Here we introduced an additional regularization term $\mathcal{G}_\mathcal{R}\colon \mathcal{R} \to \R \cup \{\infty\}$  in the objective functional. It is used to guarantee existence  for \eqref{def:neuralnetproblem}. It involves the set of admissible controls, reflecting the bounds suggested by Theorem \ref{thm:network}, given by
\begin{align*}
\mathcal{R}_{ad}:= \{\,\theta=(W_1,b_1,\dots, W_L,b_L) \in \mathcal{R}\;|\; \|W_1\|_\infty \leq \eta_1,~|b_i|_\infty\leq \eta_2,~i=1, \dots,L\,\} \subset \mathcal{R},
\end{align*}
and a penalty on the remaining network parameters:
\begin{align}\label{eq:regularizer}
\tcr{\mathcal{G}_{\mathcal{R}}(\theta)= I_{\mathcal{R}_{ad}}(\theta)+ {G}_{\mathcal{R}}(\theta)}
\end{align}
where~$I_{\mathcal{R}_{ad}}\colon \mathcal{R} \to \R \cup \{\infty\}$ denotes the convex indicator function of~$\mathcal{R}_{ad}$ i.e.
\begin{align*}
I_{\mathcal{R}_{ad}}(\theta)= \left\{
\begin{array}{ll}
0 & \theta \in \mathcal{R}_{ad}\\
+\infty & \, else \\
\end{array}
\right.
\end{align*}
and
\begin{align*}
{G}_{\mathcal{R}}(\theta)=\alpha_{\mathcal{R}} \sum^L_{i=2} \|W_i\|^2
\end{align*}
for some~$\alpha_{\mathcal{R}}>0$.
\begin{remark} \label{rem:regparameter}
In contrast to~$\eta_1, \eta_2$, the regularization parameter~$\alpha_{\mathcal{R}}$ depends on~$\mathcal{R}$. This will be crucial to show convergence of the proposed approximation scheme. For further details we refer to Section~\ref{sec:convergence}.
\end{remark}

 The remainder of this section is  dedicated to the investigation of the connection of ~\eqref{def:contproblem} to ~\eqref{def:neuralnetproblem}. In particular we prove the existence of a sequence~$\{L_{\eps}\}_{\eps>0} \in \N$ and architectures~$\operatorname{arch}(\mathcal{R}_\eps)$ such that~$(\mathcal{P}^\sigma_{\mathcal{R}_\eps})$ admits an optimal solution~$(\theta^*_\eps, \mathbf{y}^*_\eps) \in \mathcal{R}_\eps \times L^\infty_\mu(Y_0;W_\infty) $ with~$(\F^\sigma_{\theta^*_\eps},\mathbf{y}^*_\eps)$ approximating minimizers of~\eqref{def:contproblem} as~$\eps \rightarrow 0$. We proceed in three steps. First, in Section
  \ref{subsec:stabbynet}
  we establish the existence of admissible points to~\eqref{def:neuralnetproblem} for suitable architectures based on Theorem~\ref{thm:network} and nonlinear perturbation results. Next, in Section~\ref{subsec:existenceneural} the existence of optimal solutions to~\eqref{def:neuralnetproblem} is argued. Finally their convergence behavior as~$\eps \rightarrow 0$ is studied in Section~\ref{sec:convergence}.


\subsection{Stabilization by deep neural network feedback laws}
  \label{subsec:stabbynet}
In this subsection well-posedness results for \eqref{eq:cloloop} with the feedback operator $\mathcal{F}^*$
replaced by a neural network based approximation are established.  For this purpose we
require an additional assumption concerning the
linearization of \eqref{eq:cloloop}.

\begin{assumption} \label{ass:controlloflinearized}
There exists a constant $C>0$, such that for arbitrary~$y_0 \in
\mathcal{N} ( Y_0)$ and  every~$\delta y_0 \in \R^n$ and~$\delta v \in
L^2(I;\R^n)$ there exists a unique function~$\delta y \in W_\infty$
satisfying
\begin{align*}
\dot{\delta y}= D \mathbf{f}(\mathbf{y}^*(y_0))\delta
y+\mathcal{B}D\mathcal{F}^*(\mathbf{y}^*(y_0))\delta y +\delta v, \quad
\delta y(0)=\delta y_0
\end{align*}
and
\begin{align*}
\wnorm{\delta y} \leq C(\|\delta v\|_{L^2(I;\R^n)}+|\delta y_0|).
\end{align*}
\end{assumption}

%
\tcr{Sufficient conditions for Assumption~\ref{ass:controlloflinearized} to hold are discussed in Appendix~\ref{app:smoothnessofvalue}}.
For~$\eps>0$ let~$\theta_\eps \in \mathcal{R}_\eps$ be the neural network obtained from Theorem~\ref{thm:network}.
Subsequently we consider the following family of closed loop systems
\begin{align} \label{eq:neuralnetstate}
\dot{y}_\eps= \mathbf{f}(y_\epsilon)+ \mathcal{B}
\F^\sigma_{\theta_\epsilon} (y_\epsilon), \quad y_\epsilon(0)=y_0
\end{align}
where~$\F^\sigma_{\theta_\eps}$ with~$\eps>0$ is induced
by~$F^\sigma_{\theta_\eps}$. Before we turn to the main theorem of this section we establish two
preparatory results.
\begin{lemma} \label{lem:lipschitzoffeed}
Let~$ \eps >0 $ as well as~$y_1,~y_2  \in \mathcal{Y}_{ad}$ be
given. Then there holds
\begin{align*}
\|\F^*(y_1)-\mathcal{F}^\sigma_{\theta_\eps}(y_1)-(\F^*(y_2)-\F^\sigma_{\theta_\eps}(y_2))\|_{L^2(I;\R^m)}\leq \eps \wnorm{y_1-y_2}.
\end{align*}
\end{lemma}
\begin{proof}
Set~$h=y_2-y_1 \in W_\infty$.
By applying the mean value theorem we obtain
\begin{align}\label{eq:aux1}
\notag
\|\F^*(y_1)-&\F^\sigma_{\theta_\eps}(y_1)-\mathcal{F}^*(y_2)+\F^\sigma_{\theta_\eps}(y_2)\|_{L^2(I;\R^m)}
\\[1.5ex]
& \leq
\sup_{s\in
[0,1]}\|D\mathcal{F}^*(y_1+sh)-D\F^\sigma_{\theta_\eps}(y_1+sh)\|_{\mathcal{L}(W_\infty,L^2(I;\R^m))}
  \wnorm{h}.
\end{align}
We note that~$y_1+sh \in \mathcal{Y}_{ad}$ for all~$s \in [0,1]$. Using
Theorem \ref{thm:network} we estimate for every~$s \in [0,1]$
and~$\delta y \in W_\infty$
\begin{align*}
\|(D\mathcal{F}^*(y_1+sh)&-D\F^\sigma_{\theta_\eps}(y_1+sh))\delta
y\|_{L^2(I;\R^m)}\\&=
\sqrt{\int_0^\infty~|(DF^*(y_1(t)+sh(t))-DF^\sigma_{\theta_\eps}(y_1(t)+sh(t)))\delta
y(t)|^2\mathrm{d} t}
\\&\leq
\sqrt{\int_0^\infty~|DF^*(y_1(t)+sh(t))-DF^\sigma_{\theta_\eps}(y_1(t)+sh(t))|^2_{\R^{m
\times n}}|\delta y(t)|^2\mathrm{d} t}
\\&\leq \eps \|\delta y\|_{L^2(I;\R^n)} \leq \eps \wnorm{\delta y}
\end{align*}
where we utilized~$|y_1(t)+sh(t)|\leq 2 \widehat{M}$ for all~$t\in I$.
Combining this estimate with \eqref{eq:aux1} yields the result.
\end{proof}

\begin{coroll}\label{lem:calmness}
Let~$\eps>0 $ and~$y \in \mathcal{Y}_{ad}$ be given. Then there holds
\begin{align*}
\|\mathcal{F}^*(y)-\mathcal{F}^\sigma_{\theta_\eps}(y)\|_{L^2(I;\R^m)}
\leq \eps \wnorm{y}.
\end{align*}
\end{coroll}
\begin{proof}
The statement immediately follows by applying
Lemma~\ref{lem:lipschitzoffeed} with~$y_1=y$ and~$y_2=0$ as well as
using~$\F^*(0)=\F^\sigma_{\theta_\eps}(0)=0$.
\end{proof}

We are now in the position to prove the following theorem.
\begin{theorem} \label{thm:existenceneuralnetwork}
Let Assumptions~\ref{ass:feedbacklaw}-\ref{ass:controlloflinearized} hold.
There exists~$\eps_1>0$ such that~\eqref{eq:neuralnetstate}
admits a unique solution~$y_\eps=\mathbf{y}_\eps (y_0)$
in~$\mathcal{Y}_{ad}$ for all~$y_0 \in Y_0$ and~$0 <\eps <\eps_1$.
Moreover there holds
\begin{align*}
\wnorm{\mathbf{y}^*(y_0)-\mathbf{y}_\eps(y_0)} \leq c\eps
\end{align*}
with a constant~$c>0$ independent of~$y_0$ and~$\eps$. The mapping
\begin{align}\label{eq:aux2a}
\mathbf{y}_\eps \colon Y_0 \to W_\infty, \quad y_0 \mapsto
\mathbf{y}_\eps(y_0)
\end{align}
fulfills
\begin{align} \label{eq:apriorineuralnet}
\wnorm{\mathbf{y}_\eps(y_0)} \leq \frac{3}{2} M |y_0| \quad \forall y_0
\in Y_0.
\end{align}
In particular we have~$\mathbf{y}_\eps\in L^\infty_\mu(Y_0;W_\infty)$ with~$\Linfbochnorm{\mathbf{y}_\eps} \leq \frac{3}{2} M_{Y_0}$.
\end{theorem}
\begin{proof}
The proof is based on a classical fixed-point argument. Let~$y_0 \in
Y_0$ be arbitrary but fixed. Define the set
\begin{align*}
\mathcal{M}= \left\{\, y \in W_\infty \;|\;\wnorm{y} \leq \frac{3}{2}
M|y_0|\,\right\} \subset \mathcal{Y}_{ad}.
\end{align*}
On~$\mathcal{M}$ we consider the mapping~$\mathcal{Z} \colon \mathcal{M}
\to W_\infty$ where~$z=\mathcal{Z}(y)\in \mathcal{Y}_{ad}$ is the unique
solution of
\begin{align} \label{eq:auxeqfixpoint}
\dot{z}= \mathbf{f}(z)+ \mathcal{B}{\F}^*(z)+
\mathcal{B}\F^\sigma_{\theta_\epsilon}(y)-\mathcal{B} \mathcal{F}^*(y),
\quad z(0)=y_0.
\end{align}
To see that this mapping is well-posed we first note that the
perturbation function~$v=\mathcal{B}\F^\sigma_{\theta_\epsilon}(y)-
\mathcal{B} \mathcal{F}^*(y)$ fulfills
\begin{align*}
\|v\|_{L^2(I;\R^n)} &\leq \|{B}\|
\|\mathcal{F}^*(y)-\F^\sigma_{\theta_\epsilon}(y)\|_{L^2(0,\infty;\R^m)}\\[1.4ex]
& \leq \eps \|{B}\| \wnorm{y} \leq \frac{3}{2} M \eps \|{B}\||y_0| \leq
\frac{3}{2} M_{Y_0} \eps \|B\|.
\end{align*}
Note that the right hand side of the final inequality is independent
of~$y \in \mathcal{M}$ and~$y_0\in Y_0$. Thus by choosing~$\eps>0$ small
enough we may invoke Theorems~\ref{thm:existspert}
and~\ref{thm:aprioriperturbed}, respectively, to get the existence of a
unique solution~$z\in \mathcal{Y}_{ad}$ to~\eqref{eq:auxeqfixpoint} with
\begin{align*}
\wnorm{z}& \leq M|y_0|+c\|v\|_{L^2(I;\R^n)} \\[1.4ex]
& \leq M |y_0|+c \frac{3}{2} M \eps \|\mathcal{B}\|_{\mathcal{L}(L^2(I;
\R^m),L^2(I; \R^n))} |y_0|  \leq \frac{3}{2} M|y_0|
\end{align*}
where~\tcr{$c>0$ is the constant from Theorem~\ref{thm:aprioriperturbed}}.
 From this we conclude~$\mathcal{Z}(\mathcal{M})\subset \mathcal{M}$ for
all~$y_0 \in Y_0$ and~$\eps >0$ small enough. It remains to prove
that~$\mathcal{Z}$ is a contraction. To this end let~$y_1,~y_2 \in
\mathcal{M}$ be given. Applying
Corollary~\ref{coroll:locallipschitzofstate} yields
\begin{align*}
\wnorm{\mathcal{Z}(y_1)-\mathcal{Z}(y_2)} &\leq c
\|\mathcal{F}^*(y_1)-\F^\sigma_{\theta_\epsilon}(y_1)-\mathcal{F}^*(y_2)+\F^\sigma_{\theta_\epsilon}(y_2)\|_{L^2(0,\infty;\R^m)}
\\& \leq c \eps \wnorm{y_1-y_2}
\end{align*}
with a constant~$c>0$ independent of~$y_1,~y_2 \in \mathcal{M}$ as well
as of~$y_0 \in Y_0$.
Note that we also utilized the Lipschitz result of
Lemma~\ref{lem:lipschitzoffeed} in the final inequality.
Choosing~$\eps>0$ small enough we conclude that~$\mathcal{Z}$ admits a
unique fixpoint~$y_\eps=\mathcal{Z}(y_\eps) \in W_\infty$
on~$\mathcal{M}$. Clearly, the function~$ \mathbf{y}_\eps(y_0):=y_\eps$
satisfies~\eqref{eq:neuralnetstate}, \eqref{eq:apriorineuralnet} as well as
\begin{align*}
\wnorm{\mathbf{y}_\eps(y_0)-\mathbf{y}^*(y_0)}&=\wnorm{\mathcal{Z}(\mathbf{y}_\eps(y_0))-\mathcal{Z}(0)}\leq
c\eps \wnorm{y_\eps}
\leq c \eps \frac{3}{2} M_{Y_0},
\end{align*}
where the constant~$c>0$ is independent of~$y_0$. Last, we mention that
the solution to equation~\eqref{eq:neuralnetstate} is unique due to Gronwall's lemma.
It remains to proof the $\mu$-measurability of the
mapping~$\mathbf{y}_\eps$ given in \eqref{eq:aux2a}  for all~$\eps>0$
small enough. Let~$y_0 \in Y_0$ as well as an arbitrary
sequence~$\{y^k_0\}_{k\in \N} \subset Y_0$ with~$y^k_0 \rightarrow y_0$
be given. By construction we have
\begin{align*}
a(\F^\sigma_{\theta_\epsilon})(\mathbf{y}_\eps(y^k_0),y^k_0,\phi)=0
\quad \forall \phi \in W_\infty
\end{align*}
and all~$k \in \N$. Since~$\{\mathbf{y}_\eps(y^k_0)\}_{k\in\N}\subset
\mathcal{Y}_{ad}$ we can extract a subsequence, denoted by the same
symbol, with~$\mathbf{y}_\eps(y^k_0) \rightharpoonup \bar{y} \in
\mathcal{Y}_{ad}$.
Passing to the limit as~$k \rightarrow \infty$ we conclude
\begin{align*}
a(\F^\sigma_{\theta_\epsilon})(\bar{y},y_0,\phi)=0 \quad \forall \phi
\in W_\infty.
\end{align*}
Since solutions to~\eqref{eq:neuralnetstate} are unique
in~$\mathcal{Y}_{ad}$ we get~$\bar{y}=\mathbf{y}_\eps(y_0)$. Due to
the arbitrary choice of the sequence~$\{y^k_0\}_{k \in\N}$, the weak
continuity of~$\mathbf{y}_\eps$ follows. We finish the proof by noting
that weakly continuous functions are~$\mu$-measurable.
\end{proof}
\subsection{Existence of optimal neural network feedback laws} \label{subsec:existenceneural}
Let us now return to the study of problem~\eqref{def:neuralnetproblem}.
For further reference define the sets of pairs fulfilling the equality constraints in~\eqref{def:neuralnetproblem} as
\begin{align*}
\mathcal{N}^{\mathcal{R}}_{eq}:= \left\{\,(\theta, \mathbf{y})\in \mathcal{R} \times \Linfboch\;|\;\abd(\F^\sigma_{\theta})(\mathbf{y}, \Phi)=0 \quad \forall \Phi \in \Lzwoboch\right\}
\end{align*}
as well as those fulfilling the bound constraints by
\begin{align*}
\mathcal{N}^{\mathcal{R}}_{const}:= \left\{\,(\theta, \mathbf{y})\in \mathcal{R} \times \Linfboch\;|\;\theta \in \mathcal{R}_{ad},~\Linfbochnorm{\mathbf{y}}\leq 2M_{Y_0}\,\right\}.
\end{align*}
Finally we get the set of admissible pairs to~\eqref{def:neuralnetproblem} as
\begin{align*}
\mathcal{N}^\mathcal{R}_{ad}=\mathcal{N}^{\mathcal{R}}_{eq} \cap \mathcal{N}^{\mathcal{R}}_{const} \subset \mathcal{R} \times \Lzwoboch.
\end{align*}
The product space is endowed with the norm
\begin{align*}
\|(\theta, \mathbf{y})\|_{\mathcal{R} \times \Lzwoboch}= \sqrt{\|\theta\|^2_{\mathcal{R}}+\Lzwobochnorm{\mathbf{y}}^2}.
\end{align*}
In this section we shall prove the following two existence results.

\begin{theorem} \label{thm:existence}
Let Assumption~\ref{ass:feedbacklaw} and~\ref{ass:smoothofact} hold and
assume that~$\mathcal{N}^{\mathcal{R}}_{ad}\neq \emptyset$. Then there exists at least one global minimizer~$(\theta^*,\mathbf{y}_{\theta^*}) \in \mathcal{R}_{ad}\times \Linfboch$ to~\eqref{def:neuralnetproblem}.
\end{theorem}

\begin{theorem} \label{thm:condfornonempty}
Let Assumption~\ref{ass:feedbacklaw}-\ref{ass:controlloflinearized} hold and denote  by~$\operatorname{arch}(\mathcal{R}_\eps)$,~$\eps>0$, the family of architectures from Theorem~\ref{thm:network}. Then there exists~$\eps_1>0$ such that~$(\mathcal{P}^\sigma_{\mathcal{R}_\eps})$ admits a global minimizer for every~$0 <\eps<\eps_1$.
\end{theorem}

We require several technical results.
\begin{lemma} \label{lem:contofnetworks}
Let~$L \in \N $ and an architecture~$\operatorname{arch}(\mathcal{R}) \in \N^{L+1}$ be given. Consider neural networks~$\theta_k,~\theta \in \mathcal{R}$,~$k\in\N$, with~$\theta_k \rightarrow  \theta$ in~$\mathcal{R}$. Denote by~$\F^\sigma_{\theta_k}$ and~$\F^\sigma_\theta$ the Nemitsky operators induced by the corresponding realizations. Then there holds
\begin{align*}
\lim_{k \rightarrow \infty} \sup_{\Linfbochnorm{\mathbf{y}}\leq 2M_{Y_0}}\|\F^\sigma_{\theta_k}(\mathbf{y})-\F^\sigma_{\theta}(\mathbf{y})\|^2_{L^2_\mu(Y_0;L^2(I;\R^m))}=0.
\end{align*}
\end{lemma}
\begin{proof}
We proceed in several steps. First let~$N_1,N_2,N_3 \in \N$,~$g^k_1,$ let $g_1 \in \mathcal{C}^1(\R^{N_1};\R^{N_2})$ as well as~$g^k_2,~g_2 \in \mathcal{C}^1(\R^{N_2};\R^{N_3})$,~$k\in \N$, with
\begin{align*}
\lim_{k \rightarrow \infty} \lbrack \|g^k_1-g_1\|_{\mathcal{C}^1(K_1;\R^{N_2})}+\|g^k_2-g_2\|_{\mathcal{C}^1(K_2;\R^{N_3})} \rbrack=0
\end{align*}
for all compact sets~$K_1 \subset \R^{N_1}$ and~$K_2 \subset \R^{N_2}$, respectively, be given. Set~$g^k=g^k_2 \circ g^k_1$ and~$g= g_2 \circ g_1$. By a compactness argument it follows that
\begin{align} \label{eq:aux6}
\lim_{k\rightarrow \infty}\|g^k-g \|_{\mathcal{C}^1(K_1;\R^{N_3})}=0.
\end{align}
Now let us return to the proof of the stated result. We estimate
\begin{align*}
|F^\sigma_{\theta_k}(x)- F^\sigma_{\theta}(x)|&= |F^\sigma_{\theta_k}(x)- F^\sigma_{\theta}(x)-F^\sigma_{\theta_k}(0)+F^\sigma_{\theta}(0)| \\ &\leq
\sup_{|\widehat{x}|\leq 2\widehat{M}} \|DF^\sigma_{\theta_k}(\widehat{x})- DF^\sigma_{\theta}(\widehat{x})\| \, |x|
\end{align*}
for every~$x\in \R^n $,~$|x|\leq 2 \widehat{M}$. Utilizing~Assumption~\ref{ass:smoothofact} as well as~$\theta_k \rightarrow \theta$ in~$\mathcal{R}$ it is straightforward to show
\begin{align*}
\lim_{k \rightarrow \infty} \|f^\sigma_{i,\theta}-f^\sigma_{i,\theta_k}\|_{\mathcal{C}^1( K_{i-1}; \R^{N_i})}=0, \quad \text{for}~ i=1, \dots,L
\end{align*}
and arbitrary compact sets~$K_{i-1} \subset \R^{N_{i-1}}$.
We apply~\eqref{eq:aux6} repeatedly to conclude
\begin{align*}
&\lim_{k \rightarrow \infty}\sup_{|\widehat{x}|\leq 2 \widehat{M}} \|DF^\sigma_{\theta_k}(\widehat{x})- DF^\sigma_{\theta}(\widehat{x})\|\\&=\lim_{k \rightarrow \infty}\sup_{|\widehat{x}|\leq 2 \widehat{M}} \|D(f_{L,\theta}^\sigma \circ \dots \circ f_{2,\theta}^\sigma \circ f_{1,\theta}^\sigma)(\widehat{x})-D(f_{L,\theta_k}^\sigma \circ \dots \circ f_{2,\theta_k}^\sigma \circ f_{1,\theta_k}^\sigma)(\widehat{x})\|=0.
\end{align*}
For abbreviation define the null sequence
\begin{align*}
m_k:=\sup_{|\widehat{x}|\leq 2\widehat{M}} \|DF^\sigma_{\theta_k}(\widehat{x})- DF^\sigma_{\theta}(\widehat{x})\| \quad \forall k\in\N.
\end{align*}
Next consider an arbitrary~$y \in \mathcal{Y}_{ad}$. We have~$\|y\|_{\mathcal{C}_b(I;\R^n)} \leq 2 \widehat{M}$ and therefore
\begin{align*}
\|\F^\sigma_{\theta_k}(y)-\F^\sigma_{\theta}(y)\|^2_{L^2(I;\R^m)}= \int_{I} |\F^\sigma_{\theta_k}(y(t))-\F^\sigma_{\theta}(y(t))|^2~\de t
 \leq m^2_k \int_I |y(t)|^2 ~\de t  \leq 4 m^2_k M^2_{Y_0} .
\end{align*}

Finally let~$\mathbf{y}\in \Linfboch$ with~$\Linfbochnorm{\mathbf{y}}\leq 2M_{Y_0}$ be given.
Note the existence of a~$\mu$-zero set~$O \in \mathcal{A}$ such that~$\mathbf{y}(y_0)\in \mathcal{Y}_{ad}$ for~$y_0 \in Y_0 \setminus O$. Thus we can estimate
\begin{align*}
\|\F^\sigma_{\theta_k}(\mathbf{y})-\F^\sigma_{\theta}(\mathbf{y})\|^2_{L^2_\mu(Y_0;L^2(I;\R^m))}&= \int_{Y_0} \|\F^\sigma_{\theta_k}(\mathbf{y}(y_0))-\F^\sigma_{\theta}(\mathbf{y}(y_0))\|^2_{L^2(I;\R^m)}~\de \mu(y_0)
\\& \leq m^2_k \int_{Y_0} \|\mathbf{y}(y_0)\|^2_{L^2(I;\R^n)}~\de \mu(y_0)
\leq 4 m_k^2 M^2_{Y_0}.
\end{align*}
Since the right hand side tends to~$0$ as~$k \rightarrow \infty$ independently of~$\mathbf{y}\in \Linfboch$, this completes the proof.
\end{proof}

\begin{remark} \label{rem:convneuralnetpoint}
Along the lines of the previous proof we also conclude that
\begin{align*}
\lim_{k \rightarrow \infty} \sup_{y \in \mathcal{Y}_{ad}} \|\F^\sigma_{\theta_k}(y)-\F^\sigma_{\theta}(y) \|_{L^2(I;\R^m)} =0.
\end{align*}
\end{remark}
\begin{lemma} \label{lem:convofensembleforms}
Let~$\{\theta_k\}_{k\in\N} \subset \mathcal{R}$ be a given sequence of neural networks satisfying~$\theta_k \rightarrow \theta \in \mathcal{R}$. For every~$k\in\N$ assume that there exists an ensemble solution~$\mathbf{y}_k \in \Linfboch$ fulfilling
\begin{align*}
\abd(\F^\sigma_{\theta_k})(\mathbf{y}_k, \Phi)=0 \quad \forall \Phi \in \Lzwoboch,~\Linfbochnorm{\mathbf{y}_k} \leq 2 M_{Y_0}.
\end{align*}
Then there exists~$\mathbf{y}\in \Lzwoboch$ with~$\mathbf{y}_k \rightharpoonup \mathbf{y}$ in~$L^2_\mu(Y_0;W_\infty)$ as well as
\begin{align*}
\abd(\F^\sigma_{\theta})(\mathbf{y}, \Phi)=0 \quad \forall \Phi \in \Lzwoboch,~\Linfbochnorm{\mathbf{y}} \leq 2 M_{Y_0}.
\end{align*}
\end{lemma}
\begin{proof}
Since~$\{\mathbf{y}_k\}_{k \in \N}$ is bounded in~$\Linfboch$ there exists a subsequence, denoted by the same symbol, as well as an element~$\mathbf{y}\in \Lzwoboch$ with~$\mathbf{y}_k \rightharpoonup \mathbf{y}$ in~$\Lzwoboch$. Invoking Proposition~\ref{prop:parametrizedstate}, let~$O \in \mathcal{A}$ denote a~$\mu$ zero set such that
\begin{align*}
a(\F^\sigma_{\theta_k})(\mathbf{y}_k(y_0),y_0,\phi)=0 \quad \forall \phi \in W_\infty,~\wnorm{\mathbf{y}_k(y_0)} \leq 2M_{Y_0}
\end{align*}
for all~$k \in \N$ and all~$y_0 \in Y_0 \setminus O$. Fix an arbitrary~$y_0 \in Y_0 \setminus O$. Note that the sequence~$\{\mathbf{y}_k (y_0)\}_{k\in\N}$ is bounded in~$W_\infty$. Let~$\{\mathbf{y}_{k_i}(y_0)\}_{i \in \N}$ denote an arbitrary but fixed subsequence with~$\mathbf{y}_{k_i}(y_0) \rightharpoonup \widehat{\mathbf{y}}(y_0)$ for some~$\widehat{\mathbf{y}}(y_0)\in W_\infty$ as~$i \rightarrow \infty$. Due to the weak lower semicontinuity of the norm on~$W_\infty$ we immediately get~$\widehat{\mathbf{y}}(y_0) \in \mathcal{Y}_{ad}$. Next fix~$\phi \in W_\infty$. Utilizing the weak convergence of~$\{\mathbf{y}_{k_i}(y_0)\}_{i \in \N}$ yields
\begin{align*}
\lim_{i \rightarrow \infty}|\scalarl{\dot{\mathbf{y}}_{k_i}(y_0)}{\phi}-\scalarl{\dot{\widehat{\mathbf{y}}}(y_0)}{\phi}|=0
\end{align*}
and
\begin{align*}
\lim_{i \rightarrow \infty}|\scalarrn{\mathbf{y}_{k_i}(y_0)(0)}{\phi(0)}-\scalarrn{\widehat{\mathbf{y}}(y_0)(0)}{\phi(0)}|=0.
\end{align*}
Furthermore, applying Corollary \ref{cor21a} (with $F$ replaced by $f$) , we conclude
\begin{align*}
\lim_{i \rightarrow \infty} |\scalarl{\mathbf{f}(\mathbf{y}_{k_i}(y_0))}{\phi}-\scalarl{\mathbf{f}(\widehat{\mathbf{y}}(y_0))}{\phi}|=0.
\end{align*}
Last we split
\begin{align*}
\mathcal{F}^\sigma_{\theta_{k_i}}(\mathbf{y}_{k_i}(y_0))-\mathcal{F}^\sigma_{\theta}(\widehat{\mathbf{y}}(y_0))=
\mathcal{F}^\sigma_{\theta_{k_i}}(\mathbf{y}_{k_i}(y_0))-\mathcal{F}^\sigma_{\theta}(\mathbf{y}_{k_i}(y_0))+\mathcal{F}^\sigma_{\theta}(\mathbf{y}_{k_i}(y_0))-\mathcal{F}^\sigma_{\theta}(\widehat{\mathbf{y}}(y_0))
\end{align*}
to get
\begin{align*}
\lim_{i \rightarrow \infty} |\scalarl{\mathcal{B}\mathcal{F}^\sigma_{\theta_{k_i}}(\mathbf{y}_{k_i}(y_0))}{\phi}-\scalarl{\mathcal{B}\mathcal{F}^\sigma_{\theta}(\widehat{\mathbf{y}}(y_0))}{\phi}|=0.
\end{align*}
Here we again used the weak convergence of~$\{\mathbf{y}_{k_i}(y_0)\}_{i \in \N}$, Corollary \ref{cor21a}, as well as Remark~\ref{rem:convneuralnetpoint} and~$\theta_{k_i} \rightarrow \theta$ in~$\mathcal{R}$. Recalling that~$\phi \in W_\infty$ was chosen arbitrarily we arrive at
\begin{align*}
0= \lim_{i \rightarrow \infty} a(\F^\sigma_{\theta_{k_i}})(\mathbf{y}_{k_i}(y_0),y_0,\phi)= a(\F^\sigma_{\theta})(\widehat{\mathbf{y}}(y_0),y_0,\phi) \quad \forall \phi \in W_\infty.
\end{align*}
Thus~$\widehat{\mathbf{y}}(y_0) \in \mathcal{Y}_{ad}$ is the unique solution to the closed loop system~\eqref{eq:closedloopstrong} associated to~$\mathcal{F}^\sigma_\theta \in \mathcal{H}$ and~$y_0 \in Y_0 \setminus O$ cf. Proposition~\ref{prop:equivalofweakandstrong}. Combining this conclusion with the arbitrary choice of the subsequence~$\{\mathbf{y}_{k_i}\}_{i\in\N}$  we get the weak convergence of the whole sequence i.e.~$\mathbf{y}_k(y_0) \rightharpoonup \widehat{\mathbf{y}}(y_0)$. Last note that the previous argument can be repeated for every~$y_0 \in Y_0 \setminus O$ i.e. there exists a family~$\{\widehat{y}(y_0)\}_{y_0 \in Y_0 \setminus O} \subset \mathcal{Y}_{ad}$ with
\begin{align*}
\mathbf{y}_k(y_0) \rightharpoonup \widehat{\mathbf{y}}(y_0), \quad a(\F^\sigma_{\theta})(\widehat{\mathbf{y}}(y_0),y_0,\phi)=0 \quad \forall \phi \in W_\infty
\end{align*}
and every~$y_0 \in Y_0 \setminus O$.

Summarizing our previous findings we therefore have~$\Linfbochnorm{\mathbf{y}_k}\leq 2 M_{Y_0}$,~$\mathbf{y}_k \rightharpoonup \mathbf{y}$ in~$\Lzwoboch$ as well as~$\mathbf{y}_k(y_0)\rightharpoonup \widehat{\mathbf{y}}(y_0)$ in~$W_\infty$ for~$\mu$-a.e.~$y_0 \in Y_0$. In particular, the requirements of Lemma~\ref{lem:measoflimit} are met  implying that~$\mathbf{y} \in \Linfboch$ with~$\Linfbochnorm{\mathbf{y}}\leq 2M_{Y_0}$ and~$\mathbf{y}(y_0)=\widehat{\mathbf{y}}(y_0)$ for~$\mu$-a.e.~$y_0 \in Y_0$. Moreover it is evident that~$\mathbf{y}$ is an ensemble solution to~\eqref{eq:parametrizedstate} given~$\F^\sigma_\theta \in \mathcal{H}$ i.e.
\begin{align*}
\abd(\F^\sigma_{\theta})(\mathbf{y}, \Phi)=0 \quad \forall \Phi \in \Lzwoboch,
\end{align*}
see Proposition~\ref{prop:parametrizedstate}.

Finally note that the weakly convergent subsequence of~$\{\mathbf{y}_k\}_{k\in\N}$ was chosen arbitrary in the beginning and ensemble solutions are unique, see again Proposition~\ref{prop:parametrizedstate}. Thus we obtain~$\mathbf{y}_k \rightharpoonup \mathbf{{y}}$ in~$\Lzwoboch$ for the whole sequence. This completes the proof.
\end{proof}
We make the following observations.
\begin{lemma} \label{lem:weakclosedness}
There holds:
\begin{itemize}
\setlength{\itemsep}{1em}
\item \textit{Weak closedness}: The set~$\mathcal{N}^\mathcal{R}_{ad}$ is weakly closed in~$\mathcal{R}\times\Lzwoboch$.
\item \textit{Radial unboundedness}: For every sequence~$\{(\theta_k, \mathbf{y}_k)\}_{k \in \N} \subset \mathcal{N}^{\mathcal{R}}_{ad}$ we have
\begin{align*}
\|(\theta_k, \mathbf{y}_k)\|_{\mathcal{R} \times \Lzwoboch} \rightarrow \infty \Rightarrow j( \mathbf{y}_k, \mathcal{F}^\sigma_{\theta_k})+ \mathcal{G}_{\mathcal{R}}(\theta_k) \rightarrow \infty.
\end{align*}
\item \textit{Lower semicontinuity}: For every sequence~$\{(\theta_k, \mathbf{y}_k)\}_{k \in \N} \subset \mathcal{N}^{\mathcal{R}}_{ad}$ we have
\begin{align*}
(\theta_k, \mathbf{y}_k) \rightharpoonup (\theta, \mathbf{y}) \Rightarrow j( \mathbf{y}, \mathcal{F}^\sigma_\theta)+ \mathcal{G}_{\mathcal{R}}(\theta) \leq \liminf_{k \rightarrow \infty} \lbrack j( \mathbf{y}_k, \mathcal{F}^\sigma_{\theta_k})+ \mathcal{G}_{\mathcal{R}}(\theta_k)\rbrack.
\end{align*}
\end{itemize}
\end{lemma}
\begin{proof}
In the following~$\{(\theta_k, \mathbf{y}_k)\}_{k \in \N}\subset \mathcal{N}^{\mathcal{R}}_{ad}$ always denotes an arbitrary sequence. Let us prove the first statement. Assume that
\begin{align*}
\{(\theta_k, \mathbf{y}_k)\}_{k \in \N} \rightharpoonup \{(\theta, \mathbf{y})\}_{k \in \N} \quad \text{in}~\mathcal{R} \times \Lzwoboch.
\end{align*}
Clearly, this is equivalent to~$\theta_k \rightarrow \theta$ in~$\mathcal{R}$ and~$\mathbf{y}_k \rightharpoonup \mathbf{y}$ in~$\Lzwoboch$. By definition of~$\mathcal{R}_{ad}$ we immediately get~$\theta \in \mathcal{R}_{ad}$.
Applying Lemma~\ref{lem:convofensembleforms} we further conclude
\begin{align*}
\abd(\F^\sigma_{\theta})(\mathbf{y}, \Phi)=0 \quad \forall \Phi \in \Lzwoboch,~\Linfbochnorm{\mathbf{y}} \leq 2 M_{Y_0},
\end{align*}
i.e.~$(\theta, \mathbf{y})\in \mathcal{N}_{ad}$. Thus~$\mathcal{N}^{\mathcal{R}}_{ad}$ is weakly closed.
By Remark \ref{rem:convneuralnetpoint} and \eqref{eq:kk2} of Lemma \ref{lem:measoflimit} we conclude that
\begin{align*}
\F^\sigma_{\theta_k}(\mathbf{y}_k) \rightharpoonup \F^\sigma_{\theta}(\mathbf{y}) \quad \text{in}~\Lzwobochm.
\end{align*}
Next we introduce the linear and continuous mapping
\begin{align*}
\mathbf{Q} \colon \Linfboch \to \Linfboch, \quad \text{where}\quad (\mathbf{Q}\mathbf{y})(y_0)(t)=Q\mathbf{y}(y_0)(t)
\end{align*}
for~$\mu$-a.e.~$y_0 \in Y_0$ and~$t \in I$.
We rewrite~$j( \mathbf{y}_k, \mathcal{F}^\sigma_{\theta_k})$ as
\begin{align*}
j( \mathbf{y}_k,\F^\sigma_{\theta_k})= \frac{1}{2} \Lzwobochnormn{\mathbf{Q}\mathbf{y}_k}^2+ \frac{\beta}{2} \Lzwobochnormm{\mathcal{F}^\sigma_{\theta_k}(\mathbf{y}_k)}^2.
\end{align*}
Since the squared norm is convex and continuous it is weakly lower semicontinuous on~$\Lzwoboch$. As a consequence we get
\begin{align*}
\liminf_{k \rightarrow \infty} j( \mathbf{y}_k,\F^\sigma_{\theta_k}) \geq j( \mathbf{y},\F^\sigma_{\theta}).
\end{align*}
Together with~$\mathcal{G}_\mathcal{R}(\theta_k)\rightarrow \mathcal{G}_\mathcal{R}(\theta) $ this yields the third statement.

It remains to comment on the radial unboundedness of the objective functional on~$\mathcal{N}^{\mathcal{R}}_{ad}$. Assume now that we have~$\|(\theta_k, \mathbf{y}_k)\|_{\mathcal{R} \times \Lzwoboch} \rightarrow \infty$. This is equivalent to
\begin{align*}
\lim_{k\rightarrow \infty} \mathcal{G}_{\mathcal{R}}(\theta_k) / \alpha_{\mathcal{R}}= \lim_{k \rightarrow \infty} \sum^L_{i=2} \|W^k_i\|^2=+\infty.
\end{align*}
We finish the proof noting that
\begin{align*}
0 \leq \alpha_{\mathcal{R}} \, \frac{\mathcal{G}_{\mathcal{R}}(\theta_k)} { \alpha_{\mathcal{R}}} \leq j( \mathbf{y}_k,\F^\sigma_{\theta_k})+ \mathcal{G}_{\mathcal{R}}(\theta_k) \rightarrow \infty, \text{ for } k \to \infty.
\end{align*}
\end{proof}
Based on the results of the previous lemma the existence of at least one solution to~\eqref{def:neuralnetproblem} is imminent if the set of admissible pairs is not empty. This verifies the assertion of Theorem \ref{thm:existence}.

\begin{proof}[Proof of Theorem \ref{thm:condfornonempty}]
For every~$\eps>0$ denote the set of admissible pairs to~$(\mathcal{P}^\sigma_{\mathcal{R}_\eps})$ by~$\mathcal{N}^\eps_{ad}$.
We show that~$\mathcal{N}^\eps_{ad}$ is nonemtpy for all~$\eps>0$ small enough. The existence of a global minimizer  then follows from Theorem~\ref{thm:existence}. Let~$\theta_\eps$,~$\eps>0$, be the family of neural networks from Theorem~\ref{thm:network}. According to Theorem~\ref{thm:existenceneuralnetwork} there exist~$\eps_1>0$ and functions~$\mathbf{y}_\eps \in \Linfboch$ with
\begin{align*}
\abd(\F^\sigma_{\theta_\eps})(\mathbf{y}_\eps, \Phi)=0 \quad \forall \Phi \in \Lzwoboch,~\Linfbochnorm{\mathbf{y}_\eps}\leq 2 M_{Y_0}
\end{align*}
for every~$0<\eps <\eps_1$. Since~$\theta_\eps \in \mathcal{R}_{\eps,ad}$,~$\eps>0$, by construction we conclude that~$(\theta_\eps, \mathbf{y}_\eps)\in \mathcal{N}^\eps_{ad}$,~$0<\eps<\eps_1$.
\end{proof}
\begin{remark} \label{rem:alphaiszero}
To close this subsection we again point to the explicit form of the objective functional~$j$ in~\eqref{def:neuralnetproblem} which is given by
\begin{align*}
j(\mathbf{y},\mathcal{F})= \frac{1}{2}\int_{Y_0} \, \int_I (|Q \mathbf{y}(y_0)(t)|^2 +
\beta|F(\mathbf{y}(y_0)(t))|^2 dt \, d\mu(y_0).
\end{align*}
In particular recall that the parameter~$\beta$ is strictly positive. \tcr{This is often required, in the absence of further control constraints, to show } existence of minimizers to the open loop problem~\eqref{def:openloopproblem} for~$y_0 \in Y_0$. Since we parametrize admissible feedback laws by neural networks whose weights are additionally regularized this is not necessary to guarantee existence in the approximating problem~\eqref{def:neuralnetproblem}. In fact all results derived in this section remain valid for~$\beta=0$.
\end{remark}

%
%

\subsection{Convergence \& a priori estimates} \label{sec:convergence}
This subsection addresses the convergence of feedback controls  induced by realizations of optimal solutions to~$(\mathcal{P}^\sigma_{\mathcal{R}_\eps})$ as~$\eps \rightarrow 0$. For this purpose we again point to the additional regularization term~$\mathcal{G}_{\mathcal{R}_\eps}$ which was added to the objective functional in~\eqref{def:neuralnetproblem} in order to ensure the existence of global minimizers. Since no similar term appears in the original problem~\eqref{def:contproblem} it should vanish at a certain rate as~$\eps>0$ goes to zero. This will be achieved by choosing the regularization parameter~$\alpha_{\mathcal{R}_\eps}$ in dependence on the admissible set of neural networks. For preparation we make the following observation.
\begin{coroll} \label{corr:nonconstant}
Let~$\theta_\eps=(W^\eps_1,b^\eps_1, \dots, W^\eps_{L_\eps}, b^\eps_{L_\eps})$,~$\eps>0$, denote the family of neural networks from Theorem~\ref{thm:network}. Then there exists~$0< \eps_2$ such that $\sum^{L_\eps}_{i=2} \|W^\eps_i\|^2>0$ for all~$0 < \eps <\eps_2$.
\end{coroll}
\begin{proof}
Assume that~there exists a sequence~$\{\eps_{k}\}_{k\in \N}\subset \R_+ \setminus 0$ with~$\eps_k \rightarrow 0$ as well as~$\sum^{L_{\eps_k}}_{i=2} \|W^{\eps_k}_i\|^2=0$ for all~$k\in \N$. Then by definition the realizations associated to~$\theta_{\eps_k}$ are constant i.e.
\begin{align*}
F^\sigma_{\theta_{\eps_k}}(x)=F^\sigma_{\theta_{\eps_k}}(0) \quad \forall x \in \bar{B}_{2 \widehat{M}}(0)
\end{align*}
and all~$k \in \N$. Moreover due to Theorem~\ref{thm:network} we have
\begin{align*}
\lim_{k \rightarrow \infty}\sup_{x \in \bar{B}_{2\widehat{M}}(0)} |F^*(x)-F^\sigma_{\theta_{\eps_k}}(x)|=0.
\end{align*}
This gives a contradiction since~$F^*$ is not constant on~$\bar{B}_{2\widehat{M}}(0)$.
\end{proof}
The preceding corollary justifies the following assumption on the regularization parameter.
\begin{assumption} \label{ass:regparameter}
Let Assumption~\ref{ass:feedbacklaw}-\ref{ass:controlloflinearized} hold and denote by~$\operatorname{arch}(\mathcal{R}_\eps)$ and~$\theta_\eps \in \mathcal{R}_\eps$,~$\eps>0,$ the families of architectures and neural networks according to Theorem~\ref{thm:network}. The family of regularization parameters~$\{\alpha_{\mathcal{R}_\eps}\}_{\eps>0}$ is chosen such that
\begin{align*}
0<\alpha_{\mathcal{R}_\eps} \leq \frac{\eps}{2\sum^{L_\eps}_{i=2} \|W^\eps_i\|^2}, \quad 0 <\eps <\eps_2.
\end{align*}
\end{assumption}
We first derive an a priori estimate for the optimal objective function values.
\begin{theorem} \label{thm:aprioriest}
Let Assumption~\ref{ass:feedbacklaw}-\ref{ass:regparameter} hold. For~$0<\eps<\eps_2$ let~$(\theta^*_\eps, \mathbf{y}_{\theta^*_\eps})$ be a minimizing pair to~$(\mathcal{P}^\sigma_{\theta_\eps})$. Then we have
\begin{align*}
0\leq j(\mathbf{y}_{\theta^*_\eps}, \F^\sigma_{\theta^*_\eps})+\mathcal{G}_{\mathcal{R}_\eps}(\theta^*_\eps)-j( \mathbf{y}^*,\F^*)\leq c \eps
\end{align*}
for some constant~$c>0$ independent of~$\eps$. In particular we have
\begin{align*}
 j(\mathbf{y}_{\theta^*_\eps}, \F^\sigma_{\theta^*_\eps}) \rightarrow j(\mathbf{y}^*, \F^*)
\end{align*}
as~$\eps \rightarrow 0$.
\end{theorem}
\begin{proof}
First recall that the considered approximation scheme is conforming i.e.~$\mathcal{H}^\sigma_{\mathcal{R}_\eps}\subset \mathcal{H}$ for every~$\eps>0$. Thus for the family~$\{\theta_\eps\}_{\eps>0}$ from Theorem~\ref{thm:network} and~$\eps>0$ small enough we get
\begin{align} \label{eq:aux2}
0 &\leq j(\mathbf{y}_{\theta^*_\eps}, \F^\sigma_{\theta^*_\eps})-j( \mathbf{y}^*,\F^*) \leq j(\mathbf{y}_{\theta^*_\eps}, \F^\sigma_{\theta^*_\eps})+\mathcal{G}_{\mathcal{R}_\eps}(\theta^*_\eps)-j( \mathbf{y}^*,\F^*) \\
&\leq j(\mathbf{y}_{\theta_\eps}, \F^\sigma_{\theta_\eps})+\mathcal{G}_{\mathcal{R}_\eps}(\theta_\eps)-j( \mathbf{y}^*,\F^*) \notag
\end{align}
due to the optimality of~$(\F^*,\mathbf{y}^*)$ for~\eqref{def:contproblem}, the optimality of~$(\theta^*_\eps, \mathbf{y}_{\theta^*_\eps})$ for~$(\mathcal{P}^\sigma_{\mathcal{R}_\eps})$ as well as~$\mathcal{G}_{\mathcal{R}_\eps}(\theta^*_\eps)\geq 0$. Let~$0<\eps<\eps_2$ be given. We further estimate the difference on the right hand side of the last inequality. For this purpose we proceed similarly to the proof of Lemma~\ref{lem:weakclosedness} and rewrite
\begin{align*}
j&(\mathbf{y}_{\theta_\eps}, \F^\sigma_{\theta_\eps})-j( \mathbf{y}^*,\F^*)/2 =D_1+D_2
\end{align*}
where
\begin{align*}
 & D_1=\Lzwobochnormn{\mathbf{Q}\mathbf{y}_\eps}^2- \Lzwobochnormn{\mathbf{Q}\mathbf{y}^*}^2, \\
 & D_2=\beta( \Lzwobochnormm{\mathcal{F}^\sigma_{\theta_\eps}(\mathbf{y}_\eps)}^2 - \Lzwobochnormm{\F^*(\mathbf{y}^*)}^2).
\end{align*}
The first term is estimated by
\begin{align*}
&|\Lzwobochnormn{\mathbf{Q}\mathbf{y}_\eps}^2- \Lzwobochnormn{\mathbf{Q}\mathbf{y}^*}^2|\\& \leq|\Linfbochnorm{\mathbf{Q}\mathbf{y}_\eps}+ \Linfbochnorm{\mathbf{Q}\mathbf{y}^*}||\Lzwobochnormn{\mathbf{Q}\mathbf{y}_\eps}- \Lzwobochnormn{\mathbf{Q}\mathbf{y}^*}| \\
& \leq  4 \|Q\|^2 M_{Y_0} \Lzwobochnormn{\mathbf{y}_\eps-\mathbf{y}^*}
 \leq
4 \|Q\|^2 M_{Y_0}  \Linfbochnorm{\mathbf{y}_\eps-\mathbf{y}^*}
 \leq 4 \|Q\|^2 M_{Y_0} c \eps,
\end{align*}
where we used
\begin{align*}
\max\{\,\Linfbochnorm{\mathbf{y}_\eps},\Linfbochnorm{\mathbf{y}^*}\,\} \leq 2 M_{Y_0}
\end{align*}
in the second inequality and the convergence result of Theorem~\ref{thm:existenceneuralnetwork} in the last one. To bound the second term we first point to
\begin{align*}
\Lzwobochnormm{\mathcal{F}^\sigma_{\theta_\eps}(\mathbf{y}_\eps)&-\F^*(\mathbf{y}^*)}\\& \leq \Lzwobochnormm{\mathcal{F}^\sigma_{\theta_\eps}(\mathbf{y}_\eps)-\F^*(\mathbf{y}_\eps)}+\Lzwobochnormm{\F^*(\mathbf{y}_\eps)-\F^*(\mathbf{y}^*)} \\& \leq \eps \Linfbochnorm{\mathbf{y}_\eps} + L_{F^*,2\widehat{M}} \Linfbochnorm{\mathbf{y}_\eps-\mathbf{y}^*}
\leq 2 M_{Y_0} \eps+ c\eps
.
\end{align*}
Here we utilized again Theorem~\ref{thm:existenceneuralnetwork}
as well as Corollary~\ref{lem:calmness} and the Lipschitz continuity of~$F^*$ on~$\bar{B}_{2 \widehat{M}}(0)$. Similar to the first term we estimate
\begin{align*}
&|\Lzwobochnormm{\mathcal{F}^\sigma_{\theta_\eps}(\mathbf{y}_\eps)}^2 - \Lzwobochnormm{\F^*(\mathbf{y}^*)}^2|\\
& \leq (\Lzwobochnormm{\mathcal{F}^\sigma_{\theta_\eps}(\mathbf{y}_\eps)}+\Lzwobochnormm{\F^*(\mathbf{y}^*)}) \Lzwobochnormm{\mathcal{F}^\sigma_{\theta_\eps}(\mathbf{y}_\eps)\!-\!\F^*(\mathbf{y}^*)} \\
& \leq c \eps
\end{align*}
where we used the last estimate as well as
\begin{align*}
\sup_{\eps>0}\{\Lzwobochnormm{\mathcal{F}^\sigma_{\theta_\eps}(\mathbf{y}_\eps)}\}< \infty.
\end{align*}

Summarizing the previous observations we conclude
\begin{align*}
0\leq j(\mathbf{y}_{\theta^*_\eps}, \F^\sigma_{\theta^*_\eps})+\mathcal{G}_{\mathcal{R}_\eps}(\theta^*_\eps)-j( \mathbf{y}^*,\F^*) \leq c \eps+ \mathcal{G}_{\mathcal{R}_\eps}(\theta_\eps)
\end{align*}
for some constant~$c>0$ independent of~$\eps$.
The first part of the proof is now finished by noting that~$\mathcal{G}_{\mathcal{R}_\eps}(\theta_\eps)\leq \eps$ due to Assumption~\ref{ass:regparameter}. From~\eqref{eq:aux2} we finally get
\begin{align*}
0 \leq j(\mathbf{y}_{\theta^*_\eps}, \F^\sigma_{\theta^*_\eps}) - j(\mathbf{y}^*, \F^*) \leq c \eps
\end{align*}
yielding~$j(\mathbf{y}_{\theta^*_\eps}, \F^\sigma_{\theta^*_\eps}) \rightarrow j(\mathbf{y}^*, \F^*)$ as~$\eps \rightarrow 0$.
\end{proof}
The following convergence result is a direct consequence of the a priori estimate on the optimal objective functional values.
\begin{theorem}\label{thm:convergeofopt}
Let Assumption~\ref{ass:feedbacklaw}-\ref{ass:regparameter} hold. Denote by~$\{\eps_k\}_{k \in \N} \subset (0, \eps_1)$ an arbitrary null sequence. For each~$k\in \N$ let~$(\theta^*_{\eps_k}, \mathbf{y}_{\theta^*_{\eps_k}}) \in \mathcal{R}_\eps \times \Linfboch$ be a minimizing pair to~$(\mathcal{P}^\sigma_{\eps_k})$. The sequence~$\{\,(\mathbf{y}_{\theta^*_{\eps_k}},\F^\sigma_{\theta^*_{\eps_k}}(\mathbf{y}_{\theta^*_{\eps_k}}))\,\}_{k \in \N}$ admits at least one weak accumulation point~$(\widehat{\mathbf{y}}, \widehat{\mathbf{u}})$ in~$\Lzwoboch \times \Lzwobochm$. Each such point fulfills~$\Linfbochnorm{\widehat{\mathbf{y}}}\leq 2 M_{Y_0}$ as well as
\begin{align*}
\dot{\widehat{\mathbf{y}}}(y_0)= \mathbf{f}(\widehat{\mathbf{y}}(y_0))+ \mathcal{B}\widehat{\mathbf{u}}(y_0),~\widehat{\mathbf{y}}(y_0)(0)=y_0, \quad (\widehat{\mathbf{y}}(y_0), \widehat{\mathbf{u}}(y_0)) \in \argmin \eqref{def:openloopproblem}
\end{align*}
for~$\mu$-a.e.~$y_0 \in Y_0$.

If~$Q$ is invertible every weakly convergent subsequence of~$\{\,(\mathbf{y}_{\theta^*_{\eps_k}},\F^\sigma_{\theta^*_{\eps_k}}(\mathbf{y}_{\theta^*_{\eps_k}}))\,\}_{k \in \N}$ in~$\Lzwoboch \times \Lzwobochm$ also converges strongly.
\end{theorem}
\begin{proof}
First note that~$\Linfbochnorm{\mathbf{y}_{\theta^*_{\eps_k}}}\leq 2 M_{Y_0}$ as well as
\begin{align*}
\Lzwobochnormm{\F^\sigma_{\theta^*_{\eps_k}}(\mathbf{y}_{\theta^*_{\eps_k}})} \leq \frac{j(\mathbf{y}_{\theta^*_{\eps_k}},\F^\sigma_{\theta^*_{\eps_k}}(\mathbf{y}^*_{\theta^*_{\eps_k}}))}{\beta}
\end{align*}
for every~$k \in \N$. Since the right handside of this inequality converges as~$k \rightarrow \infty$, the sequence~$(\mathbf{y}_{\theta^*_{\eps_k}}, \F^\sigma_{\theta^*_{\eps_k}}(\mathbf{y}_{\theta^*_{\eps_k}}))$ is bounded in~$\Linfboch \times \Lzwobochm$ and thus admits a weakly convergent subsequence in~$\Lzwoboch \times \Lzwobochm$ denoted by the same symbol. Let~$(\widehat{\mathbf{y}}, \widehat{\mathbf{u}}) \in \Lzwoboch \times \Lzwobochm$ be its weak limit. Due to the  boundedness of~$\{\mathbf{y}_{\theta^*_{\eps_k}}\}_{k\in\N}$ in~$\Linfboch$ there holds~$\widehat{\mathbf{y}} \in \Linfboch$ with~$\Linfbochnorm{\widehat{\mathbf{y}}}\leq 2M_{Y_0}$.

By arguing similarly as in the proof of Lemma~\ref{lem:convofensembleforms} we conclude
\begin{align*}
0&=\lim_{k \rightarrow \infty} \abd(\F^\sigma_{\theta_{\eps_k}})( \mathbf{y}_{\theta^*_{\eps_k}} , \Phi) \\
&= \int_{Y_0} \lbrack \scalarl{\dot{\widehat{\mathbf{y}}}(y_0)-\mathbf{f}(\widehat{\mathbf{y}}(y_0))-\mathcal{B}\widehat{\mathbf{u}}(y_0)}
{\Phi(y_0)}+\scalarrn{\widehat{\mathbf{y}}(y_0)(0)\!-\!y_0}{\Phi(y_0)(0)} \rbrack~ \de \mu(y_0),
\end{align*}
for every~$\Phi \in \Lzwoboch$. Repeating the arguments of Propositions~\ref{prop:equivalofweakandstrong} and~\ref{prop:parametrizedstate} this implies
\begin{align*}
\scalarl{\dot{\widehat{\mathbf{y}}}(y_0)-\mathbf{f}(\widehat{\mathbf{y}}(y_0))-\mathcal{B}\widehat{\mathbf{u}}(y_0)}{\phi}+\scalarrn{\widehat{\mathbf{y}}(y_0)(0)-y_0}{\phi(0)}=0 \quad \forall \phi \in W_\infty
\end{align*}
and thus
\begin{align*}
\dot{\widehat{\mathbf{y}}}(y_0)= \mathbf{f}(\widehat{\mathbf{y}}(y_0))+ \mathcal{B} \widehat{\mathbf{u}}(y_0),~\widehat{\mathbf{y}}(y_0)(0)=y_0
\end{align*}
for~$\mu$-a.e~$y_0 \in Y_0$. Due to the weak lower semicontinuity of the norms on~$\Lzwobochn$ and~$\Lzwobochm$, respectively, we arrive at
\begin{align*}
\int_{Y_0} J(\widehat{\mathbf{y}}(y_0), \widehat{\mathbf{u}}(y_0))~\de \mu(y_0) \leq \liminf_{k \rightarrow \infty} j(\mathbf{y}_{\theta^*_{\eps_k}}, \F^\sigma_{\theta^*_{\eps_k}})= j(\mathbf{y}^*,\F^*)= \int_{Y_0} V(y_0)~\de \mu(y_0).
\end{align*}
According to the previous discussions the pair~$(\widehat{\mathbf{y}}(y_0), \widehat{\mathbf{u}}(y_0)) \in W_\infty \in L^2(I;\R^m)$ is admissible for~\eqref{def:openloopproblem} for~$\mu$-a.e.~$y_0 \in Y_0$. Therefore we have~$V(y_0)\leq J(\widehat{\mathbf{y}}(y_0), \widehat{\mathbf{u}}(y_0))$ for~$\mu$-a.e.~$y_0 \in Y_0$ and
\begin{align*}
\int_{Y_0} V(y_0)~\de \mu(y_0) \leq \int_{Y_0} J(\widehat{\mathbf{y}}(y_0), \widehat{\mathbf{u}}(y_0))~\de \mu(y_0).
\end{align*}
Thus we finally get~$J(\widehat{\mathbf{y}}(y_0), \widehat{\mathbf{u}}(y_0))=V(y_0)$ as well as
\begin{align*}
(\widehat{\mathbf{y}}(y_0), \widehat{\mathbf{u}}(y_0)) \in \argmin \eqref{def:openloopproblem}
\end{align*}
for~$\mu$-a.e.~$y_0 \in Y_0$
by definition of the value function.

Last if $Q$ is invertible the mapping
\begin{align*}
\widehat{\jmath} \colon \Lzwobochn \times \Lzwobochm \to \R
\end{align*}
with
\begin{align*}
\widehat{\jmath}(\mathbf{y}, \mathbf{u})= \sqrt{\frac{1}{2}\Lzwobochnormn{\mathbf{Q}\mathbf{y}}^2+\frac{\beta}{2} \Lzwobochnormm{\mathbf{u}}^2}
\end{align*}
defines a norm on~$\Lzwobochn \times \Lzwobochm$ which is induced by an inner product and equivalent to the canonical norm
\begin{align*}
\|(\mathbf{y}, \mathbf{u})\|_{\Lzwobochn \times \Lzwobochm}= \sqrt{\Lzwobochnormn{\mathbf{y}}^2+\Lzwobochnormm{\mathbf{u}}^2}.
\end{align*}
Clearly~$\Lzwobochn \times \Lzwobochm$ is complete with respect to~$\widehat{\jmath}$ making it a Hilbert space. Moreover, from Theorem~\ref{thm:aprioriest} we get
\begin{align*}
\widehat{\jmath}(\mathbf{y}_{\theta^*_{\eps_k}},\F^\sigma_{\theta^*_{\eps_k}}(\mathbf{y}^*_{\theta^*_{\eps_k}}))^2= j(\mathbf{y}_{\theta^*_{\eps_k}},\F^\sigma_{\theta^*_{\eps_k}}(\mathbf{y}^*_{\theta^*_{\eps_k}})) \rightarrow j(\widehat{\mathbf{y}}, \widehat{\mathbf{u}})= \widehat{\jmath}(\widehat{\mathbf{y}}, \widehat{\mathbf{u}})^2.
\end{align*}
Utilizing the weak convergence of the sequence as well as convergence of the norms of its elements we conclude
\begin{align*}
\lim_{k \rightarrow \infty}\widehat{\jmath}(\mathbf{y}_{\theta^*_{\eps_k}}-\widehat{\mathbf{y}}, \F^\sigma_{\theta^*_{\eps_k}}(\mathbf{y}_{\theta^*_{\eps_k}})-\widehat{\mathbf{u}})=0.
\end{align*}
Finally the equivalence of~$\widehat{\jmath}$ and the canoncial norm yields
\begin{align*}
\lim_{k \rightarrow \infty}\|(\mathbf{y}_{\theta^*_{\eps_k}}, \F^\sigma_{\theta^*_{\eps_k}}(\mathbf{y}_{\theta^*_{\eps_k}}))-(\widehat{\mathbf{y}},\widehat{\mathbf{u}})\|_{\Lzwobochn \times \Lzwobochm}=0.
\end{align*}
Last we again apply Lemma~\ref{lem:Nemitsky} to conclude
\begin{align*}
\|\dot{\mathbf{y}}_{\theta^*_{\eps_k}} &- \dot{\widehat{\mathbf{y}}}\|_{L^2_\mu (Y_0; L^2(I;\R^n))} \\& \leq c (\|\mathbf{f}(\mathbf{y}_{\theta^*_{\eps_k}}) - \mathbf{f}(\widehat{\mathbf{y}})\|_{L^2_\mu (Y_0; L^2(I;\R^n))}+\|\F^\sigma_{\theta^*_{\eps_k}}(\mathbf{y}_{\theta^*_{\eps_k}})-\widehat{\mathbf{u}}\|_{L^2_\mu (Y_0; L^2(I;\R^m))})
\\ & \leq
c (\|\mathbf{y}_{\theta^*_{\eps_k}} - \widehat{\mathbf{y}}\|_{L^2_\mu (Y_0; L^2(I;\R^n))}+\|\F^\sigma_{\theta^*_{\eps_k}}(\mathbf{y}_{\theta^*_{\eps_k}})-\widehat{\mathbf{u}}\|_{L^2_\mu (Y_0; L^2(I;\R^m))}).
\end{align*}
Thus~$\mathbf{y}_{\theta^*_{\eps_k}} \rightarrow \widehat{\mathbf{y}}$ in~$\Lzwoboch$.
The proof is finished noting that the weakly convergent subsequence was chosen arbitrarily in the beginning.
\end{proof}
\begin{remark} \label{rem:problem}
Note that the previous theorem does not ensure the existence of a Nemitsky operator~$\widehat{\F} \in \mathcal{H}$ such that
\begin{align} \label{eq:feddbacklimit}
\dot{\widehat{\mathbf{y}}}(y_0)= \mathbf{f}(\widehat{\mathbf{y}}(y_0))+ \mathcal{B} \widehat{\F}(\widehat{\mathbf{y}}(y_0)),~\widehat{\mathbf{y}}(y_0)(0)=y_0,~\widehat{\mathbf{u}}(y_0)= \widehat{\F}(\widehat{\mathbf{y}}(y_0))
\end{align}
for~$\mu$-a.e.~$y_0 \in Y_0$. In practice, however, we can only solve instances of the approximating problem~$(\mathcal{P}^\sigma_{\mathcal{R}_\eps})$ for small but nonzero~$\eps>0$ yielding controls in feedback form. Nevertheless, sufficient conditions on the considered neural networks which ensure the existence of an operator~$\widehat{\F}$ fulfilling~\eqref{eq:feddbacklimit} would be of great theoretical interest. We postpone a thorough discussion of this question to future research.
\end{remark}
\section{Practical realization} \label{sec:practicalrealization}
We address selected topics for the computational realization of the learning problem. First the
infinite time horizon in~$(\mathcal{P}^\sigma_{\mathcal{R}_\eps})$ is replaced by a problem with finite time horizon. Together with its first order optimality  conditions it is discussed in Section \ref{subsec:firstorder}. These first order conditions are subsequently
used within a gradient based algorithm. Combined  with a
discretization approach this is   described in Section~\ref{subsec:discandalg}.

\subsection{Finite time horizon problem} \label{subsec:firstorder}

We fix~$\eps>0$ and~$T>0$,  and set~$I_T=[0,T]$. For every~$y_0 \in Y_0$ and~$\theta \in \mathcal{R}_\eps$ we consider the finite time horizon learning  problem given by
\begin{equation} \label{def:neuralnetproblemfinite} \tag{$\mathcal{P}^\sigma_{{\mathcal{R}_\eps},T}$}
\left\{
\begin{aligned}
\quad & \min_{\substack{\theta \in \mathcal{R}_{\eps,ad}, \\ \mathbf{y} \in L^\infty_\mu(Y_0;W_T)}} j_T(\mathbf{y},\F^\sigma_\theta)+ \mathcal{G}_{\mathcal{R}_\eps}(\theta)
\\
&s.t. \quad  \abd_T(\F^\sigma_\theta)(\mathbf{y}, \Phi) =0 \quad \forall \Phi \in L^2_\mu(Y_0;W_T),~\|\mathbf{y}\|_{L^\infty_\mu(Y_0;W_T)}\leq 2 M_{Y_0},
\end{aligned}
\right.
\end{equation}
where
\begin{align*}
j_T(\mathbf{y},\mathcal{F})= \frac{1}{2}\int_{Y_0} \, \int_{I_T} (|Q \mathbf{y}(y_0)(t)|^2 +
\beta|F(\mathbf{y}(y_0)(t))|^2 dt \, d\mu(y_0),
\end{align*}
and $W_T$ and $\abd_T(\F^\sigma_\theta)$ denote  the canonical  restrictions of $W_\infty$ and $\abd$ to $I_T$.

\tcr{It is straightforward to argue existence for a minimizing pair~$(\theta^*_{\eps,T},\mathbf{y}_{\theta^*_{\eps,T}})$ of~\eqref{def:neuralnetproblemfinite} if~$(\mathcal{P}_{\mathcal{R}\eps})$ admits a minimizer.
We formally derive the following necessary first order conditions, a detailed proof, however, goes beyond the scope of this paper. For simplicity let us focus on the case of inactive constraints i.e. we assume that
\begin{align*}
\theta^*_{\eps,T} \in \operatorname{int} \mathcal{R}_{ad,\eps}, \quad \|\mathbf{y}_{\theta^*_{\eps,T}}\|_{L^\infty_\mu(Y_0;W_T)}< 2 M_{Y_0}.
\end{align*}
}

Then there exists~$\mathbf{p}_{\theta^*_{\eps,T}} \in L^\infty_\mu(Y_0;W_T)$ such that
\begin{align}
(\dot{\mathbf{y}}_{\theta^*_{\eps,T}}-\mathbf{f}(\mathbf{y}_{\theta^*_{\eps,T}})-\mathcal{B}\mathcal{F}^\sigma_{\theta^*_{\eps,T}}(\mathbf{y}_{\theta^*_{\eps,T}}),\Phi)_{L^2_\mu(Y_0;L^2(I;\R^n))}&+
(\mathbf{y}_{\theta^*_{\eps,T}}(0)-\mathbf{y}_0,\Phi(0))_{L^2_\mu(Y_0;\R^n)} \notag \\
&=0 \quad
 \forall \Phi \in L^2_\mu(Y_0; W_T), \label{eq:stateeqoptimal}
\end{align}
\begin{align}
(\dot{\mathbf{p}}_{\theta^*_{\eps,T}}&+D\mathbf{f}(\mathbf{y}_{\theta^*_{\eps,T}})^*\mathbf{p}_{\theta^*_{\eps,T}}+D_{x}\mathcal{F}^\sigma_{\theta^*_{\eps,T}}(\mathbf{y}_{\theta^*_{\eps,T}})^*\mathcal{B}^*\mathbf{p}_{\theta^*_{\eps,T}},\Phi)_{L^2_\mu(Y_0;L^2(I;\R^n))}-(\mathbf{p}_{\theta^*_{\eps,T}}(T),\Phi(T))_{L^2_\mu(Y_0;\R^n)} \notag \\
&=-(\mathcal{Q}^*\mathcal{Q}\mathbf{y}_{\theta^*_{\eps,T}}+\beta D_x \F^\sigma_{\theta^*_{\eps,T}}(\mathbf{y}_{\theta^*_{\eps,T}})^*\F^\sigma_{\theta^*_{\eps,T}}(\mathbf{y}_{\theta^*_{\eps,T}}),\Phi)_{L^2_\mu(Y_0;L^2(I;\R^n))}
 \quad \forall \Phi \in L^2_\mu(Y_0; W_T), \label{eq:adjointeqtwoopt}
\end{align}
as well as
\begin{align}
\left( \widehat{\theta}_{\eps,T}+S'_{\eps,T}(\theta^*_{\eps,T})^*\nu_{\theta^*_{\eps,T}} +\nabla{G}_{\mathcal{R}_\eps}(\theta^*_{\eps,T} ), \theta-\theta^*_{\eps,T}\right)_{\mathcal{R}_\eps}\geq 0 \quad \forall \theta \in \mathcal{R}_{ad,\eps}\label{eq:gradcondopt}
\end{align}
where
\begin{align*}
\nabla {G}_{\mathcal{R}_\eps}(\theta^*_{\eps,T})=2\alpha_{\mathcal{R}_\eps} (0,0,W^*_2,0,\cdots, W^*_L,0) \in \mathcal{R}_\eps
\end{align*}
and
\begin{align*}
 \widehat{\theta}_{\eps,T}:=\int_{Y_0} \int^T_0 \left \lbrack D_\theta \F^\sigma_{\theta^*_{\eps,T}}(\mathbf{y}_{\theta^*_{\eps,T}})^* \mathcal{B}^* \mathbf{p}_{\theta^*_{\eps,T}}+\beta D_\theta \F^\sigma_{\theta^*_{\eps,T}}(\mathbf{y}_{\theta^*_{\eps,T}})^* \F^\sigma_{\theta^*_{\eps,T}}(\mathbf{y}_{\theta^*_{\eps,T}})\right \rbrack~\de t~\de \mu(y_0).
\end{align*}


\subsection{Discretization \& Algorithmic Solution}
\label{subsec:discandalg}
In order to compute a minimizer of \eqref{def:neuralnetproblemfinite} we also have to approximate the integrals with respect to the probability  measure~$\mu$ as well as the closed loop system and the adjoint equation. While a thorough discussion of different discretization schemes including a rigorous analysis of their convergence is postponed to future research we nevertheless briefly outline our chosen approach. First replace the measure~$\mu$ by a suitable convex combination of Dirac Delta functions
\begin{align*}
\mu_N= \sum^N_{i=1} \omega_i \delta_{y^i_0}, \quad \sum^N_{i=1} \omega_i=1,~\omega_i>0
\end{align*}
for ~$i=1,\dots,N$, supported on a finite set~$\widehat{Y}_0 \subset Y_0$. Without loss of generality we may assume~$0 \in \widehat{Y}_0$.

The new ensemble objective functional is then given by
\begin{align*}
j_{N,T}(\mathbf{y},\mathcal{F})= \frac{1}{2} \sum^N_{i=1} \omega_i \left\lbrack\int_{I_T} (|Q \mathbf{y}(\tcr{y^i_0})(t)|^2 +
\beta|F(\mathbf{y}(y^i_0)(t))|^2 dt \right \rbrack.
\end{align*}
Due to the finite cardinality of~$\widehat{Y}_0$ we identify
\begin{align*}
L^2_{\mu_N}(\widehat{Y}_0;W_T)\simeq (W^N_T,\|\cdot\|_{W_T,2}), \quad \|\mathbf{y}\|_{W_T,2}= \sqrt{\sum^N_{i=1} \omega_i \|\mathbf{y}_i\|^2_{W_T}},
\end{align*}
as well as
\begin{align*}
L^\infty_{\mu_N}(\widehat{Y}_0;W_T)\simeq (W^N_T,\|\cdot\|_{W_T,\infty}), \quad \|\mathbf{y}\|_{W_T,\infty}= \max_{i=1,\dots,N} \|\mathbf{y}_i\|_{W_T,\infty}.
\end{align*}
Further there holds~$L^\infty_{\mu_N}(\widehat{Y}_0;W_T)\simeq \mathcal{C}_b(\widehat{Y}_0;W_T).$
This leads to the  learning problem which we treat numerically
\begin{align} \label{def:problemdiscstoch}
\min_{\substack{\theta \in \mathcal{R}_{\eps,ad}, \\ \mathbf{y} \in W^N_T }}\left \lbrack\frac{1}{2} \sum^N_{i=1} \omega_i \left\lbrack\int_{I_T} (|Q \mathbf{y}_i(t)|^2 +
\beta|F^\sigma_\theta(\mathbf{y}_i(t))|^2 dt \right \rbrack+ \mathcal{G}_{\mathcal{R}_\eps}(\theta) \right \rbrack,
\end{align}
subject to the constraints
\begin{align*}
\dot{\mathbf{y}}_i= \mathbf{f}(\mathbf{y}_i)+ \F^\sigma_\theta(\mathbf{y}_i),\quad \mathbf{y}_i(0)=y^i_0, \quad \|\mathbf{y}_i\|_{W_T} \leq 2M_{Y_0}, \quad i=1,\dots,N.
\end{align*}

Again, we put ourselves in the situation when there exists a minimizing pair~$(\theta^*,\mathbf{y}^*)\in \mathcal{R}_{ad,\eps}\times W^N_T$ for~\eqref{def:problemdiscstoch} with
\begin{align*}
\theta^* \in \operatorname{int} \mathcal{R}_{ad,\eps}, \quad \|\mathbf{y}^*_i\|_{W_T}< 2M_{Y_0}, \quad i=1,\dots,N,
\end{align*}
satisfying
\begin{align} \label{eq:stateass}
\dot{\mathbf{y}}_i= \mathbf{f}(\mathbf{y}_i)+ \F^\sigma_\theta(\mathbf{y}_i),\quad \mathbf{y}_i(0)=y^i_0, \quad i=1,\dots,N.
\end{align}
These assumptions together with~\eqref{eq:stateeqoptimal}-\eqref{eq:gradcondopt} imply the existence of~$\mathbf{p}^*\in W^T_N$ such that the triple~$(\theta^*,\mathbf{y}^*,\mathbf{p}^*)$ is a solution to
\begin{align}
\dot{\mathbf{y}}_i= \mathbf{f}(\mathbf{y}_i)+ \F^\sigma_\theta(\mathbf{y}_i),~ \mathbf{y}_i(0)=y^i_0 \label{eq:stateeqdisc} \\
-\dot{\mathbf{p}}_i= D \mathbf{f}(\mathbf{y}_i)^*\mathbf{p}_i+ D_x \F^\sigma_\theta(\mathbf{y}_i)^*\mathcal{B}^*\mathbf{p}_i+\mathcal{Q}^*\mathcal{Q}\mathbf{y}_i+\beta D_x \F^\sigma_{\theta}(\mathbf{y}_i)^*\F^\sigma_{\theta}(\mathbf{y}_i),~ \mathbf{p}_i(T)=0
\end{align}
for all~$i=1,\dots,N$ and
\begin{align} \label{eq:condeqdisc}
\nabla{G}_{\mathcal{R}_\eps}(\theta )+ \sum^N_{i=1} \omega_i \int^T_0 \left \lbrack D_\theta \F^\sigma_{\theta}(\mathbf{y}_i)^* \mathcal{B}^* \mathbf{p}_i+\beta D_\theta \F^\sigma_{\theta}(\mathbf{y}_i)^* \F^\sigma_{\theta}(\mathbf{y}_i)\right \rbrack~\de t =0.
\end{align}
In practice we compute a triple $(\theta,\mathbf{y},\mathbf{p})$ satisfying~\eqref{eq:stateeqdisc}-\eqref{eq:condeqdisc} by applying a gradient descent method to the reduced problem
\begin{align*}
\min_{\theta  \in \mathcal{R}_\eps} \left \lbrack\frac{1}{2} \sum^N_{i=1} \omega_i \left\lbrack\int_{I_T} (|Q S(\theta)_i(t)|^2 +
\beta|F^\sigma_\theta(S(\theta)_i(t))|^2 dt \right \rbrack+ \mathcal{G}_{\mathcal{R}_\eps}(\theta) \right \rbrack,
\end{align*}
where~$S$ denotes the operator mapping neural networks to the vector of solutions to system~\eqref{eq:stateass}. The procedure is summarized in Algorithm~\ref{alg:grad}.
\begin{algorithm}
\begin{algorithmic}
\STATE 1. Let~$\theta^1 \in \mathcal{R}_\eps$.
\WHILE {not converged}
\STATE 2. For~$i=1,\dots,N$ solve
\begin{align*}
\dot{\mathbf{y}}^k_i= \mathbf{f}(\mathbf{y}^k_i)+ \F^\sigma_{\theta^k}(\mathbf{y}^k_i),\quad \mathbf{y}^k_i(0)=y^i_0.
\end{align*}
as well as
\begin{align*}
-\dot{\mathbf{p}}^k_i= D \mathbf{f}(\mathbf{y}^k_i)^*\mathbf{p}^k_i+ D_x \F^\sigma_{\theta^k}(\mathbf{y}^k_i)^*\mathcal{B}^*\mathbf{p}^k_i+\mathcal{Q}^*\mathcal{Q}\mathbf{y}^k_i+\beta D_x \F^\sigma_{\theta^k}(\mathbf{y}^k_i)^*\F^\sigma_{\theta}(\mathbf{y}^k_i)
\end{align*}
with~$\mathbf{p}^k_i(T)=0$.
\STATE 4. Compute
\begin{align*}
\widehat{\theta}^k=\nabla{G}_{\mathcal{R}_\eps}(\theta^k )+ \sum^N_{i=1} \omega_i \int^T_0 \left \lbrack D_\theta \F^\sigma_{\theta^k}(\mathbf{y}^k_i)^* \mathcal{B}^* \mathbf{p}^k_i+\beta D_\theta \F^\sigma_{\theta^k}(\mathbf{y}^k_i)^* \F^\sigma_{\theta^k}(\mathbf{y}^k_i)\right \rbrack~\de t.
\end{align*}
\STATE 5. Choose~$s^k>0$ and set~$\theta^{k+1}=\theta^k-s^k \widehat{\theta}^k$.
\ENDWHILE
\end{algorithmic}
\caption{Gradient descent feedback learning}
\label{alg:grad}
\end{algorithm}
Concerning the choice of the stepsize~$s^k$ in step 5. we note that each evaluation of the reduced objective functional requires the solution of~$N$ nonlinear ordinary differential equations. Hence implicit stepsize rules such as Armijo backtracking are infeasible for the problem. In our implementation we select~$s^k$ according to the Barzilai-Borwein stepsize rule. Concretely this method is based on
\begin{align*}
s^1_{\text{BB}}:= \frac{(\mathcal{S}_{k-1},\mathcal{E}_{k-1})_{\mathcal{R}_\eps}}{(\mathcal{S}_{k-1},\mathcal{S}_{k-1})_{\mathcal{R}_\eps}} \quad \text{and}\quad s^2_{\text{BB}}:= \frac{(\mathcal{E}_{k-1},\mathcal{E}_{k-1})_{\mathcal{R}_\eps}}{(\mathcal{S}_{k-1},\mathcal{E}_{k-1})_{\mathcal{R}_\eps}},
\end{align*}
with~$\mathcal{S}_{k-1}=\theta^k-\theta^{k-1}$ and~$\mathcal{E}_{k-1}=\widehat{\theta}^k-\widehat{\theta}^{k-1}$. Subsequently we either choose
\begin{align*}
s^k= \max \left\{s_{\text{min}},\min \left\{s^1_{\text{BB}},s_{\text{max}}\right\}\right\} \quad \text{or} \quad s^k= \max \left\{s_{\text{min}},\min \left\{s^2_{\text{BB}},s_{\text{max}}\right\}\right\}.
\end{align*}
where~$s_{\text{min}},s_{\text{max}}>0$ are fixed constants independent of~$k\in\N$.
At the same moment we have to mention that the computation of the reduced gradient in steps 2.-4. requires the solution of~$N$ nonlinear closed loop systems and~$N$ linear adjoint equations per iteration. While only a moderate number~$N$ of initial conditions is considered in our numerical experiments we propose a stochastic version of Algorithm~\ref{alg:grad} and/or the use of inexact gradients for large~$N$.

\section{Numerical Examples} \label{sec:examples}
In the last section we report on three examples which illustrate the
practical applicability of the proposed approach. It is based on
\eqref{eq:stateass}-\eqref{eq:condeqdisc} where the state dynamics are
replaced by a continuous Galerkin approximation of order one and all
arising temporal integrals are treated by the trapezoidal rule. This
corresponds to a Crank-Nicolson scheme for the closed loop system.
Accordingly we derive first order necessary optimality conditions for
the discretized learning problem. Finally a neural network feedback law
is obtained from Algorithm~\ref{alg:grad}.

The finite time horizon~$T$, the number of layers~$L$ as well as the
activation function~$\sigma$ vary between the considered examples and
are chosen based on numerical testing. \tcr{Our experiences with the learning problem suggest that the time horizon~$T$ has to be chosen sufficiently large to ensure successful stabilization}. We set~$N_i=n$,~$i=1,\dots,L-1$
and slightly depart from the presentation in Section~\ref{sec:neuralnet}
by considering layers of the form
\begin{align*}
f^\sigma_{i,\theta}(x)= \sigma(W_ix+b_i)+x \quad \forall x \in
\R^n,~i=1,\dots,L-1
\end{align*}
including an additional identity mapping. Such skip or shortcut
connections,~\cite{W18,He15}, prevent the Jacobians of the individual
layers from vanishing for small weights and thus greatly help in the
training of deep neural networks.

Special attention is paid to comparisons between effects  of the neural
network feedbacks and the feedback laws obtained by the linear-quadratic regulator~(LQR) applied to the linearized
system. It is given by~$\F_{\text{LQR}}(y)(t)=-(1/\beta)B^\top \Pi
y(t)$, c.f.~\cite[Chapter~6.2]{CurZ95}, where~$\Pi \in \R^{n \times n}$ denotes the unique
positive-definite solution to the algebraic Riccati equation
\begin{align} \label{eq:algebraicric}
A^\top \Pi+ \Pi A-(1/\beta)\Pi B^\top B \Pi+Q=0  \quad \text{where}
\quad A= Df(0).
\end{align}
Comparisons are also carried out to  the power series expansion type
controller (PSE)
\begin{align*}
\F_{\text{PSE}}(y)(t)=-(1/\beta) B^\top(\Pi y(t)-(A^\top-(1/\beta)\Pi
BB^\top \Pi)^{-1}\Pi f_l(y(t))
\end{align*}
with~$f_l$ describing the nonlinear part of the state dynamics. Such
feedback laws can be related to formal Taylor approximations of first
and second order for the value function around zero. For further
reference see e.g.~\cite{Ray10}.

While we do not aim here at a quantitative comparison among different
feedback methodologies with respect
to  e.g. approximation  capabilities and computational effort, we
nonetheless point out that, at least
for the examples  considered in the following, even with only a moderate
number of initial conditions for
the learning  problem and with early stopping after few iterations of
~Algorithm~\ref{alg:grad} the NN-based feedback gain provides competitive results.
There are cases when it is successful,
while the Riccati- or the PSE-controllers fail.  We point out that the
computations  in step 2. can be
parallelized leaving the overall  computation time in the range of
minutes. Finally the application
of NN-feedback laws only require the  computation of matrix-vector
products as well as the application
of the nonlinear activation function  which can be efficiently realized.
All computations in this
section have been conducted in Matlab 2019a.
\subsection*{LC-circuit}
\begin{wrapfigure}{r}{4.5cm}
  \includegraphics[width=4cm]{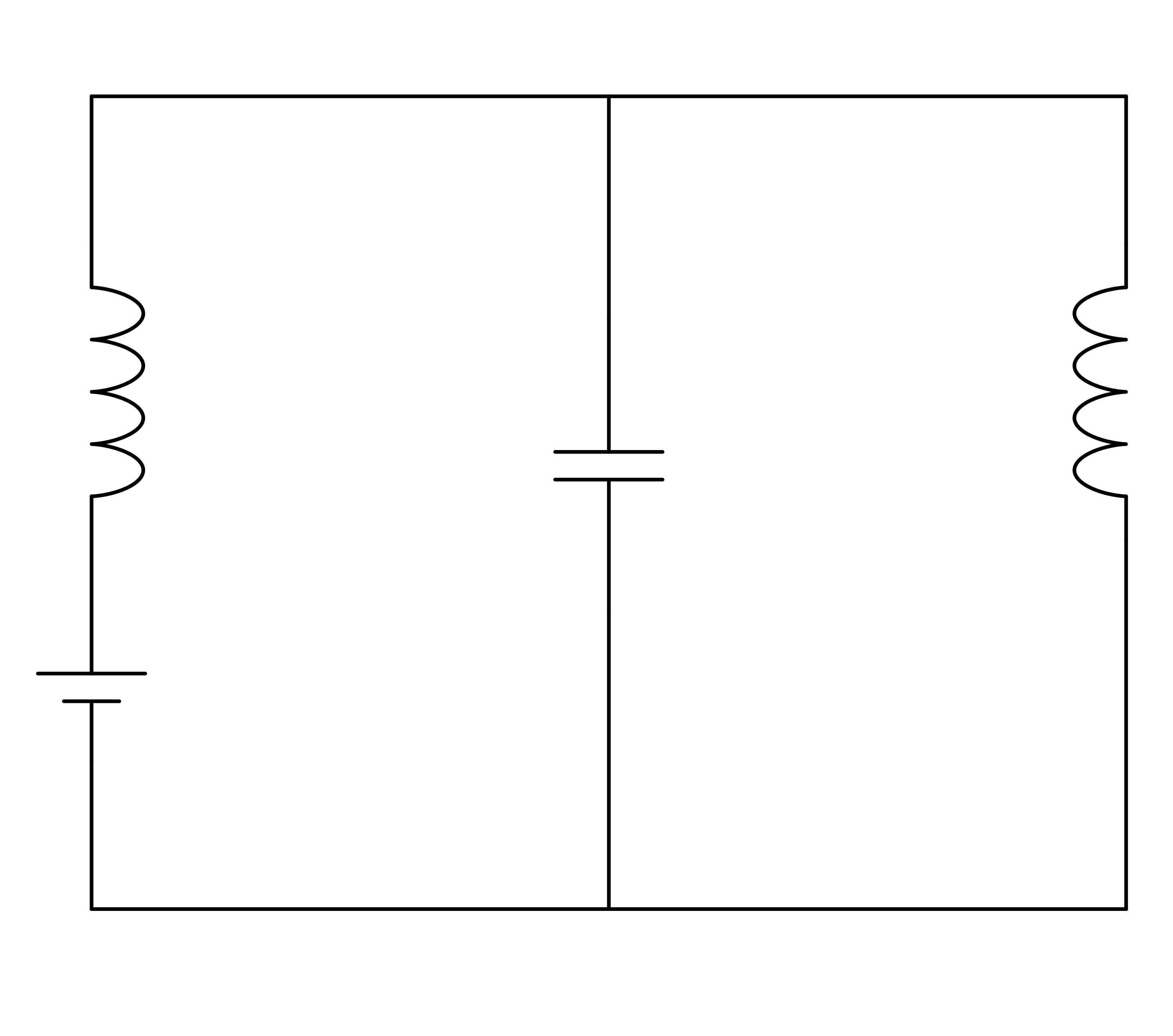}
  \caption{LC-circuit}
  \label{fig:circuit}
\end{wrapfigure}
As a first example consider a LC-circuit consisting of two inductors of unit inductivity and a capacitator of unit capacity see~\cite[Example 4.2.1]{vdS96}. The setup is schematized in Figure \ref{fig:circuit}.
The magnetic fluxes in the left and right hand side inductor are described by time-dependent functions~$\phi_1$ and~$\phi_2$, respectively. Let~$q$ denote the charge stored on the capacitor at any given time, and set~$y=(\phi_1,\phi_2,q)$. We assume the linearity of all involved elements. As a consequence the combined magnetic and electric energy in the circuit at time~$t\in [0,\infty)$ is given by~$\frac{1}{2}\,|y(t)|^2$.
The system can be influenced by applying a voltage~$u$ through the generator at the lower left side. According to Kirchhoff's laws the vector-valued function~$y$ thus follows
\begin{align} \label{eq:circuit}
\dot{y}=\left( \begin{array}{rrr}
0 & 1 & -1 \\ -1 & 0 & 0 \\ 1 &0 &0
\end{array}\right)y+ \left(\begin{array}{r} 0 \\ 1 \\ 0
\end{array} \right) u, \quad y(0)=y_0,
\end{align}
where~$y_0 \in \R^3$ is the initial distribution of magnetic fluxes and charge. We stress that in the absence of a voltage~$u$ the system is norm-constant i.e.~$|y(t)|=|y_0|$ for all~$t \in [0,\infty)$.

Given~$y_0$ our aim is to drive the energy stored in the circuit to zero by applying a voltage~$u$ in feedback form. We set~$Q=\operatorname{Id}\in \R^{3 \times 3}$, and~$\beta=0.1$.

Since ~\eqref{eq:circuit} is linear the value function of the associated open loop optimal control problem is~$V(y_0)= 1/2\, y^\top_0 \Pi y_0$, where~$\Pi$ denotes the solution of the corresponding Riccati equation~\eqref{eq:algebraicric}.
Thus as described in Section \ref{sec3.2} an optimal feedback law is given by~$F^*(y_0)=-1/\beta\, B^\top \Pi y_0$, and  ~$W^*=-1/\beta\, B^\top \Pi \in \R^{1\times 3}$ is the Riccati feedback gain.

\begin{figure}[htb]
\centering
\begin{subfigure}[t]{.45\linewidth}
\scalebox{.48}
{\input{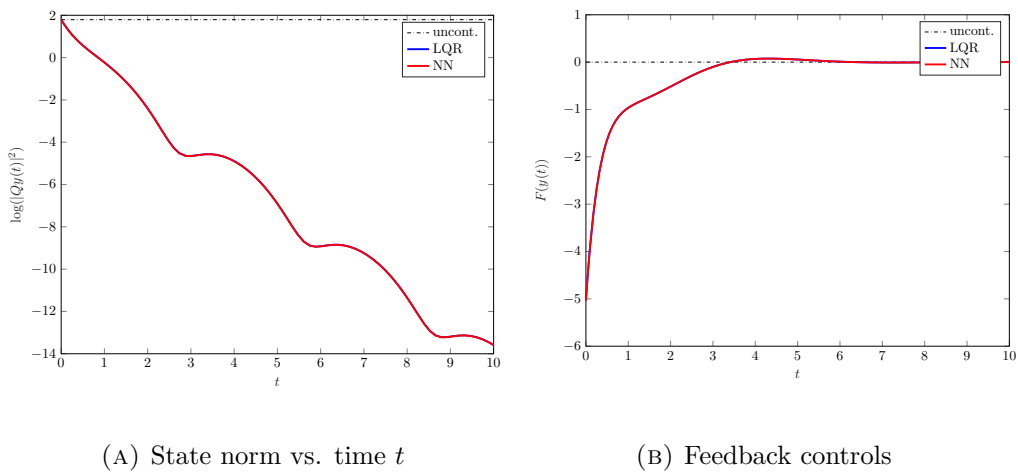}}\\
\caption{State norm vs. time~$t$}
\end{subfigure}
\begin{subfigure}[t]{.45\linewidth}
\scalebox{.47}
{
%
%
\begin{tikzpicture}

\begin{axis}[%
width=4.667in,
height=3.681in,
at={(0.783in,0.497in)},
scale only axis,
xmin=0,
xmax=10,
xlabel style={font=\color{white!15!black}},
xlabel={$t$},
ymin=-6,
ymax=1,
ylabel style={font=\color{white!15!black}},
ylabel={$F(y(t))$},
axis background/.style={fill=white},
legend style={legend cell align=left, align=left, draw=white!15!black}
]
\addplot [color=black, dashdotted]
  table[row sep=crcr]{%
0	0\\
0.0248878506744702	0\\
0.0497757013489404	0\\
0.0746635520234106	0\\
0.0995514026978808	0\\
0.175510214935839	0\\
0.251469027173797	0\\
0.327427839411756	0\\
0.403386651649714	0\\
0.477435269009474	0\\
0.551483886369234	0\\
0.625532503728994	0\\
0.699581121088754	0\\
0.779311267623463	0\\
0.859041414158173	0\\
0.938771560692882	0\\
1.01850170722759	0\\
1.10641188560017	0\\
1.19432206397276	0\\
1.28223224234534	0\\
1.37014242071792	0\\
1.46904048798509	0\\
1.56793855525226	0\\
1.66683662251943	0\\
1.7657346897866	0\\
1.87810700343565	0\\
1.9904793170847	0\\
2.10285163073375	0\\
2.2152239443828	0\\
2.34418160209877	0\\
2.47313925981474	0\\
2.60209691753072	0\\
2.73105457524669	0\\
2.88123327939415	0\\
3.03141198354162	0\\
3.18159068768908	0\\
3.33176939183655	0\\
3.4711392031884	0\\
3.61050901454025	0\\
3.74987882589209	0\\
3.88924863724394	0\\
4.03298086382336	0\\
4.17671309040277	0\\
4.32044531698218	0\\
4.46417754356159	0\\
4.59031387906046	0\\
4.71645021455933	0\\
4.8425865500582	0\\
4.96872288555707	0\\
5.09485922105594	0\\
5.2209955565548	0\\
5.34713189205367	0\\
5.47326822755254	0\\
5.61147512536495	0\\
5.74968202317736	0\\
5.88788892098977	0\\
6.02609581880218	0\\
6.15313159799401	0\\
6.28016737718584	0\\
6.40720315637766	0\\
6.53423893556949	0\\
6.66765474366559	0\\
6.8010705517617	0\\
6.9344863598578	0\\
7.0679021679539	0\\
7.2127678944546	0\\
7.3576336209553	0\\
7.502499347456	0\\
7.6473650739567	0\\
7.80102707138918	0\\
7.95468906882167	0\\
8.10835106625416	0\\
8.26201306368664	0\\
8.3942628778449	0\\
8.52651269200315	0\\
8.6587625061614	0\\
8.79101232031965	0\\
8.92735632518717	0\\
9.06370033005468	0\\
9.2000443349222	0\\
9.33638833978971	0\\
9.47394749331553	0\\
9.61150664684135	0\\
9.74906580036717	0\\
9.88662495389299	0\\
9.91496871541974	0\\
9.9433124769465	0\\
9.97165623847325	0\\
10	0\\
};
\addlegendentry{uncont.}

\addplot [color=blue, line width=1.5pt]
  table[row sep=crcr]{%
0	-5.04177735408707\\
0.0248592261418808	-4.71873032970943\\
0.0497184522837617	-4.41940495235828\\
0.0745776784256425	-4.14216532168564\\
0.0994369045675233	-3.88547379867058\\
0.175292390405325	-3.21373685450406\\
0.251147876243126	-2.68549594685266\\
0.327003362080928	-2.27291524285154\\
0.402858847918729	-1.95073549483392\\
0.476808291115508	-1.70404924266157\\
0.550757734312288	-1.51119582393957\\
0.624707177509068	-1.36101674605378\\
0.698656620705847	-1.2434638716043\\
0.778286677959331	-1.14403216044749\\
0.857916735212815	-1.06640201119103\\
0.937546792466298	-1.00500024863774\\
1.01717684971978	-0.954790910382244\\
1.10498530428415	-0.907698064537202\\
1.19279375884853	-0.866574759863523\\
1.2806022134129	-0.829092303909617\\
1.36841066797727	-0.793098920228463\\
1.4672077038794	-0.752490634758724\\
1.56600473978152	-0.711309413622071\\
1.66480177568365	-0.669206802354664\\
1.76359881158577	-0.625751412912283\\
1.87587468667859	-0.574559496937753\\
1.98815056177142	-0.522283942800639\\
2.10042643686424	-0.469651993268385\\
2.21270231195706	-0.417118701000704\\
2.34156760995539	-0.357590343952696\\
2.47043290795373	-0.300114159226548\\
2.59929820595206	-0.245637077168747\\
2.72816350395039	-0.194721966600427\\
2.87834667191079	-0.140455502989212\\
3.02852983987119	-0.0923948083349644\\
3.17871300783159	-0.0509009255635523\\
3.32889617579198	-0.015960897914314\\
3.46826550692587	0.0107736577084834\\
3.60763483805976	0.0322409096224876\\
3.74700416919364	0.048802260583428\\
3.88637350032753	0.0608778313082875\\
4.03010063744857	0.0691237479089516\\
4.1738277745696	0.0736660907290103\\
4.31755491169064	0.0750801708813196\\
4.46128204881167	0.0739004478013612\\
4.58736555096708	0.0711342096707989\\
4.71344905312248	0.067113656990602\\
4.83953255527788	0.0621510027769727\\
4.96561605743329	0.0565195219180176\\
5.09169955958869	0.050462162570263\\
5.21778306174409	0.0442015090895693\\
5.3438665638995	0.0379248782164737\\
5.4699500660549	0.0317846415575121\\
5.60812830480862	0.025357865082601\\
5.74630654356235	0.0193721671060509\\
5.88448478231607	0.0139197455490183\\
6.02266302106979	0.0090601463510879\\
6.1496168212826	0.00514550752097478\\
6.27657062149541	0.0017645457900641\\
6.40352442170822	-0.00109082187302281\\
6.53047822192104	-0.00343825908778382\\
6.66383729550371	-0.00538635133426488\\
6.79719636908638	-0.00684682631104196\\
6.93055544266906	-0.00786772755338323\\
7.06391451625173	-0.00849757589030588\\
7.20874869855508	-0.00879676307943821\\
7.35358288085844	-0.00876326178530044\\
7.49841706316179	-0.008462969430629\\
7.64325124546514	-0.00795421744316623\\
7.79735595621741	-0.00724375699776195\\
7.95146066696968	-0.00642036618727091\\
8.10556537772196	-0.00553709093661978\\
8.25967008847423	-0.00463589070245238\\
8.3918784425142	-0.00387445775312833\\
8.52408679655417	-0.00314385673267625\\
8.65629515059414	-0.00245810805256256\\
8.78850350463411	-0.00182740468163742\\
8.92482128416352	-0.00124212575101956\\
9.06113906369293	-0.00072726449783633\\
9.19745684322235	-0.000284481765860271\\
9.33377462275176	8.66474647021044e-05\\
9.4712838167466	0.0003907496147323\\
9.60879301074144	0.00062867411830416\\
9.74630220473628	0.000805822059119367\\
9.88381139873112	0.00092795205362763\\
9.91285854904834	0.000947291608017235\\
9.94190569936556	0.000964528603652967\\
9.97095284968278	0.00097972631982707\\
10	0.000992948126620473\\
};
\addlegendentry{LQR}

\addplot [color=red, line width=1.5pt]
  table[row sep=crcr]{%
0	-5.03712874103067\\
0.0248878506744702	-4.71445764546426\\
0.0497757013489404	-4.41547796118483\\
0.0746635520234106	-4.13855613375254\\
0.0995514026978808	-3.88215670093031\\
0.175510214935839	-3.21104849115364\\
0.251469027173797	-2.68331526437949\\
0.327427839411756	-2.27114038897843\\
0.403386651649714	-1.94928337926416\\
0.477435269009474	-1.70285085540237\\
0.551483886369234	-1.51019241846296\\
0.625532503728994	-1.36015801579793\\
0.699581121088754	-1.24270783235071\\
0.779311267623463	-1.14335529996693\\
0.859041414158173	-1.06576645277809\\
0.938771560692882	-1.00437487177502\\
1.01850170722759	-0.954151365755318\\
1.10641188560017	-0.90702597760372\\
1.19432206397276	-0.865851755264614\\
1.28223224234534	-0.828305442559919\\
1.37014242071792	-0.792240718916114\\
1.46904048798509	-0.751554524778332\\
1.56793855525226	-0.710295945144406\\
1.66683662251943	-0.668120260555828\\
1.7657346897866	-0.62459978251749\\
1.87810700343565	-0.57335567529491\\
1.9904793170847	-0.521043730937882\\
2.10285163073375	-0.468392062286373\\
2.2152239443828	-0.415856578287575\\
2.34418160209877	-0.356355762428216\\
2.47313925981474	-0.298927526880597\\
2.60209691753072	-0.244516215299515\\
2.73105457524669	-0.193682361016177\\
2.88123327939415	-0.139564860145323\\
3.03141198354162	-0.0916553571328904\\
3.18159068768908	-0.0503092320771033\\
3.33176939183655	-0.0155099331812316\\
3.4711392031884	0.0111035860701712\\
3.61050901454025	0.0324604374752406\\
3.74987882589209	0.0489234120202235\\
3.88924863724394	0.0609133663956365\\
4.03298086382336	0.0690848409738168\\
4.17671309040277	0.0735662784036067\\
4.32044531698218	0.0749325095690625\\
4.46417754356159	0.0737172417375966\\
4.59031387906046	0.0709278409344962\\
4.71645021455933	0.0668911863455427\\
4.8425865500582	0.0619189676865107\\
4.96872288555707	0.0562838704172667\\
5.09485922105594	0.0502281968391055\\
5.2209955565548	0.043973820711029\\
5.34713189205367	0.0377073412992289\\
5.47326822755254	0.0315804376188236\\
5.61147512536495	0.0251724447818358\\
5.74968202317736	0.0192074505694652\\
5.88788892098977	0.0137767937581084\\
6.02609581880218	0.00893929527880312\\
6.15313159799401	0.00504313598355663\\
6.28016737718584	0.00168044731903138\\
6.40720315637766	-0.00115721787406147\\
6.53423893556949	-0.00348782100900916\\
6.66765474366559	-0.00541906166449056\\
6.8010705517617	-0.00686434727923344\\
6.9344863598578	-0.00787185219464031\\
7.0679021679539	-0.00849016275660531\\
7.2127678944546	-0.00877891545460283\\
7.3576336209553	-0.0087372572219472\\
7.502499347456	-0.00843093978620568\\
7.6473650739567	-0.00791811960342631\\
7.80102707138918	-0.00720767103888177\\
7.95468906882167	-0.00638650928562289\\
8.10835106625416	-0.00550693142486046\\
8.26201306368664	-0.00461032034876417\\
8.3942628778449	-0.00385079512618439\\
8.52651269200315	-0.003122371107583\\
8.6587625061614	-0.00243897931096906\\
8.79101232031965	-0.0018107340975206\\
8.92735632518717	-0.00122808826917995\\
9.06370033005468	-0.000715806105115457\\
9.2000443349222	-0.000275490075767162\\
9.33638833978971	9.33297486465086e-05\\
9.47394749331553	0.000395341848239251\\
9.61150664684135	0.000631374945445514\\
9.74906580036717	0.000806848697312563\\
9.88662495389299	0.000927531620135752\\
9.91496871541974	0.000946151818945031\\
9.9433124769465	0.000962779801176569\\
9.97165623847325	0.000977474264629229\\
10	0.000990293985757101\\
};
\addlegendentry{NN}

\end{axis}

\begin{axis}[%
width=6.022in,
height=4.516in,
at={(0in,0in)},
scale only axis,
xmin=0,
xmax=1,
ymin=0,
ymax=1,
axis line style={draw=none},
ticks=none,
axis x line*=bottom,
axis y line*=left,
legend style={legend cell align=left, align=left, draw=white!15!black}
]
\end{axis}
\end{tikzpicture}%
}\\
\caption{Feedback controls}
\end{subfigure}
\caption{Results for~$y_0=(-1,2,1)^\top$.}
\label{fig:circuitres}
\end{figure}

Clearly in this case it suffices to consider linear feedback laws in the learning problem i.e.~$L$ is set to one and we only optimize~$W_L\in \R^{1\times 3}$. The finite time horizon is chosen as~$T=20$ and the training set of initial conditions solely contains the first canonical basis vector~$(1,0,0)^\top$. Algorithm~\ref{alg:grad} stops with
\begin{align*}
W_L=-(3.567,4.137,0.331)
\end{align*}
after~$38$ iterations with gradient norm smaller than~$10^{-6}$.
As a comparison the true Riccati feedback gain is computed using the icare routine in Matlab. This results in
\begin{align*}
W^*=-(3.571,4.140,0.332).
\end{align*}
Remarkably the learning based approach with a single initial condition yields a good approximation to the true Riccati feedback feedback gain. The optimal and neural network feedback controls for~$y_0=(-1,2,1)^\top$ as well as the evolution of the norm of the closed loop states  are depicted in Figure~\ref{fig:circuitres}.

\subsection*{Van-der-Pol like system}
The second example addresses the stabilization of a nonlinear oscillator whose dynamics are described by
\begin{align}\label{eq:nonlinearosc}
\ddot{Y}(t)=1.5\,(1-Y^2(t)) \dot{Y}(t)-Y(t)+0.8\, Y^3(t)+u, \quad (Y(0),\dot{Y}(0))=y_0,
\end{align}
for some~$y_0 \in \R^2$.
This corresponds to the governing equation of a Van-der-Pol oscillator with an added cubic nonlinearity which destabilizes the system. We equivalently rewrite~\eqref{eq:nonlinearosc} as a two-dimensional system
\begin{align*}
\dot{y}= \left( \begin{array}{rr}
0 & 1 \\ -1 & 1.5
\end{array}\right)y+ \left( \begin{array}{c}
0 \\ -1.5\,y^2_1 y_2-y_1+0.8\,y^3_1
\end{array} \right)+ \left(\begin{array}{r} 0 \\ 1
\end{array} \right) u, \quad y(0)=y_0,
\end{align*}
where~$y=(Y,\dot{Y})\in W_\infty$. We set $Q=(1,0)^\top(1,0)$.
Further we fix
$$ \beta=0.001,\,  T=3, \text{ and }\sigma(x)=\max(x^{1.01},0).$$
A set of~$20$ initial conditions  is obtained by sampling from a uniform distribution on~$Y_0=[-10,10]^2$ . Subsequently it is equally split into a set~$\widehat{Y}_0=\{y^i_0\}^{10}_{i=1}$ used as training data in the learning problem with equal weights~$w_i=1/10$,~$i=1,\dots,10$, and~its complement on which we validate the computed results. In order to stabilize the system we trained neural networks with~$L=3$ and~$5$ layers, respectively. Since both lead to very similar results we only report on those for the smaller network  to avoid ambiguities in the following.

\begin{figure}[htb]
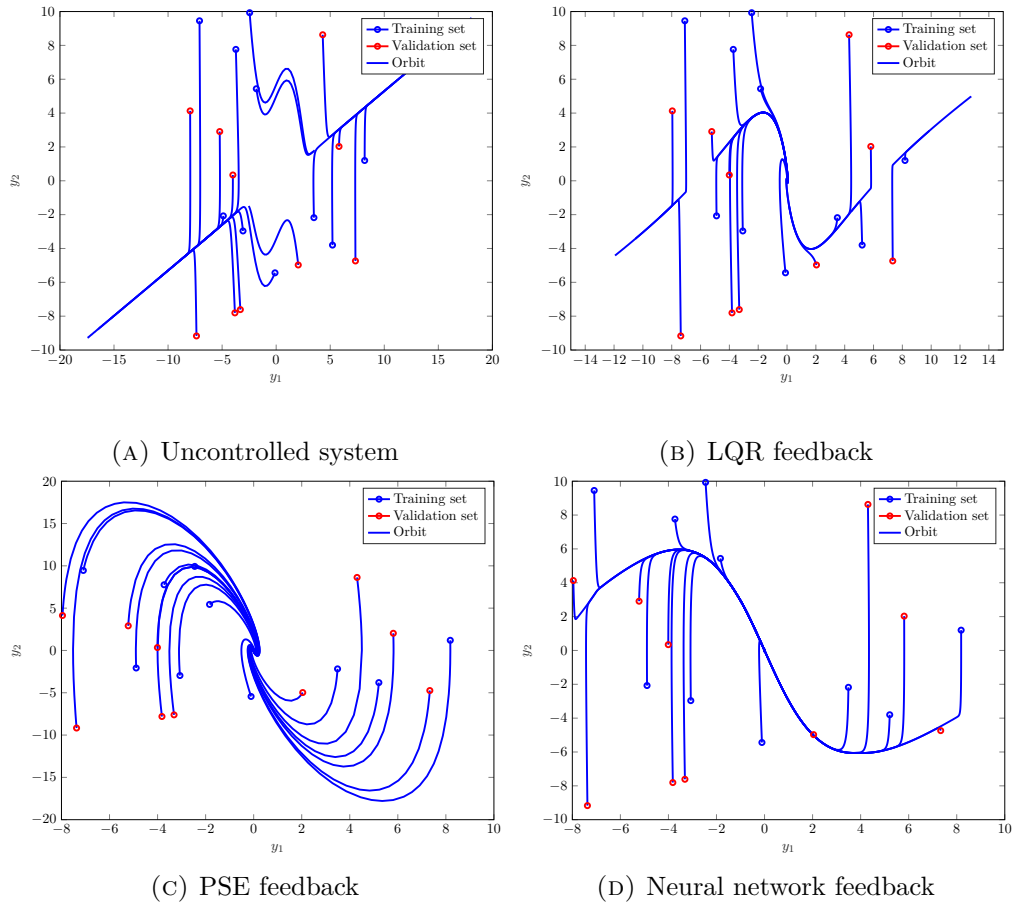

\centering
\begin{subfigure}[t]{.45\linewidth}
\scalebox{.48}
{\input{orbituncont.tex}}\\
\caption{Uncontrolled system}
\end{subfigure}
\begin{subfigure}[t]{.45\linewidth}
\scalebox{.48}
{\input{orbitric.tex}}\\
\caption{LQR feedback}
\end{subfigure}
\begin{subfigure}[t]{.45\linewidth}
\centering
\scalebox{.48}
{\input{orbitpse.tex}}\\
\caption{PSE feedback}
\end{subfigure}
\begin{subfigure}[t]{.45\linewidth}
\centering
\scalebox{.48}
{\input{orbitnn.tex}}\\
\caption{Neural network feedback}
\end{subfigure}
\caption{Orbits of the closed loop systems}
\label{fig:orbits}
\end{figure}

We illustrate the effect of the different feedback controllers on the system. For this purpose the orbits of the associated closed loop dynamics originating from the test and validation sets are plotted in Figure~\ref{fig:orbits}. In Table~\ref{tab:vander} the state and feedback control norms as well as open loop objective functional values, approximated by a suitably large time horizon, are summarized for several initial conditions in the validation set.

As we can see from the results presented  in Figure~\ref{fig:orbits} as well as Table~\ref{tab:vander} the uncontrolled trajectories diverge rapidly. This is due to the presence of the cubic nonlinearity.  In contrast the linearization based feedback stabilizes the system in most of the considered test cases but fails for initial conditions with relatively large~$y_1$ component. Last, while both nonlinear feedback laws succeed for all sampled initial conditions the computed trajectories remarkably differ, especially in the transient phase. In fact the neural network feedback controller first drives the state almost in parallel to the~$y_2$-axis towards a curve crossing the origin. On this lower dimensional manifold the trajectories are controlled to the origin. We point out that none of the initial conditions considered in the learning problem (indicated by blue dots in the Figures) lie directly on this curve. Again, in view of the Bellman principle, this highlights the influence of the whole trajectory, and not just the initial condition, on the learned feedback controller cf. Corollary~\ref{cor:largerset}.

 The initial changes in the~$y_2$ component of the PSE closed-loop state differ from the neural network dynamics both in magnitude and, in some cases, in sign. In particular note that the neural network states do not leave~$Y_0$ while this occurs for the PSE dynamics if the first component of the initial condition is too large or small.

  Comparing the computed results in Table~\ref{tab:vander} we observe that the LQR and PSE feedback controllers lead to smaller objective functional values if the first component of~$y_0$ is relatively close to~zero. This can be expected since both feedback laws are constructed based on a polynomial expansion of the value function around the origin. Moreover, in contrast to LQR, the PSE controller efficiently stabilizes for initial conditions with relatively large/small first component. However, we observe substantially larger control norms in comparison to the neural network feedback and thus also larger open loop objective functional values.
\begin{table}[!htb]
\begin{minipage}{.48\linewidth}
    \centering

    \medskip
\begin{tabular}{ c c c c }

\multicolumn{4}{c}{$y_0=(-7.37,-9.17)$}\\
\toprule
$\mathcal{F}$ &  $\|Qy\|_{L^2}$ & $\|\F(y)\|_{L^2}$ &  $J(y,\mathcal{F}(y))$  \\
\midrule
uncont. & $+\infty$ & 0 & $+\infty$\\
LQR & $+\infty$ & $+\infty$ & $+\infty$\\
PSE & 3.7 & 494 & 129\\
NN  & 5.76 & 379 & 88.5 \\
\bottomrule
\end{tabular}
\end{minipage}\hfill
\begin{minipage}{.48\linewidth}
    \centering

    \medskip
\begin{tabular}{ c c c c }

\multicolumn{4}{c}{$y_0=(2.04,-4.97)$}\\
\toprule
$\mathcal{F}$ &  $\|Qy\|_{L^2}$ & $\|\F(y)\|_{L^2}$ &  $J(y,\mathcal{F}(y))$  \\
\midrule
uncont. & $+\infty$ & 0 & $+\infty$\\
LQR & 0.82 & 7.85 & 0.37\\
PSE & 0.72 & 14.7 & 0.37\\
NN  & 0.87 & 6.66 & 0.4 \\
\bottomrule
\end{tabular}
\end{minipage}\hfill
\begin{minipage}{.48\linewidth}
    \centering

    \medskip
\begin{tabular}{ c c c c }

\multicolumn{4}{c}{$y_0=(5.81,2.03)$}\\
\toprule
$\mathcal{F}$ &  $\|Qy\|_{L^2}$ & $\|\F(y)\|_{L^2}$ &  $J(y,\mathcal{F}(y))$  \\
\midrule
uncont. & $+\infty$ & 0 & $+\infty$\\
LQR & 7.32 &  205  &  47.8\\
PSE  & 2.7 &  240  &  32.4\\
NN  & 3.46 &  200  &  25.88  \\
\bottomrule
\end{tabular}
\end{minipage}\hfill
\begin{minipage}{.48\linewidth}
    \centering

    \medskip
\begin{tabular}{ c c c c }

\multicolumn{4}{c}{$y_0=(-3.31,-7.61)$}\\
\toprule
$\mathcal{F}$ &  $\|Qy\|_{L^2}$ & $\|\F(y)\|_{L^2}$ &  $J(y,\mathcal{F}(y))$  \\
\midrule
uncont. & $+\infty$ & 0 & $+\infty$\\
LQR  & 2.26 &  52.4   &  3.91\\
PSE  & 1.75 &  68.7 &  3.89\\
NN  & 1.71 &  81.0  &  4.74  \\
\bottomrule
\end{tabular}
\end{minipage}\hfill
\caption{Results for~different~$y_0$}
\label{tab:vander}
\end{table}

\subsection*{Viscous Burgers' equation}
The final example is dedicated to the stabilization of a one dimensional Burgers'-like equation given by
\begin{align*}
\partial_t \mathcal{Y}(x,t)&= 0.2 \partial_{xx} \mathcal{Y}(x,t)+ \mathcal{Y}(x,t) \partial_x Y(x,t)+ \delta \mathcal{Y}^p+ \chi_{\omega}(x) u(t) \\
 \mathcal{Y}(-1,t)&= \mathcal{Y}(1,t)=0 \\
\mathcal{Y}(x,0)&= \mathcal{Y}_0,
\end{align*}
for all~$(x,t) \in (-1,1)\times I$,~$\delta \geq 0$,~$p\in \{1,3\}$ and initial datum~$\mathcal{Y}_0$. The time-dependent control signal~$u\in L^2(I)$ acts on the subset~$\omega=(-0.5,-0.2)$ of the spatial domain, with~$\chi_\omega$ denoting its characteristic function.

To fit this setting into the perspective of the  manuscript we approximate the state by~$\mathcal{Y}_N (x,t)=\sum_{j=0}^N y_j(t)  \phi_j(x)$ where~$N\in\N$,~$x_j=\cos(j\pi /N)$,~$j=0,\dots,N$, are the Chebyshev nodes, and~$\{\phi_j\}^N_{j=0}$ is the set of Lagrange basis polynomials associated to~$\{x_j\}^N_{j=0}$. The time-dependent coefficient function~$y_j$ corresponds to an approximation of the state~$\mathcal{Y}$ at the~$j$-th collocation point~$x_j$. In more detail, we set~$y_0=y_N=0$, and require~$y=(y_1, \dots,y_{N-1})^\top$ to fulfill
\begin{align} \label{eq:burgersapprox}
\dot{y}(t)= D_1 y(t)+ y(t) * D_2 y(t)+\delta \,y^p+ Bu  , \quad
\end{align}
where the matrices~$D_1,D_2 \in \R^{(N-1)\times(N-1)}$ are pseudospectral approximations of the first and second derivative with respect to~$x$, $y^p$ denotes the coordinate-wise $p$-th power of $y$,  and
$y_0=(\mathcal{Y}_0(x_1), \dots, \mathcal{Y}_0(x_{N-1}))$.
The control operator is given by~$B=(\chi(x_1), \dots,\chi(x_{N-1}))$, and the symbol~$*$ denotes the Hadamard product between two vectors.
We approximate the~$L^2$-norm of~$\mathcal{Y}$ by
\begin{align*}
\int_I \int^1_{-1} |\mathcal{Y}(x,t)|^2~\mathrm{d}x~\mathrm{d}t \approx \int_I |Qy(t)|^2~\mathrm{d}t
\end{align*}
 where~$Q$ is a diagonal matrix containing the square roots of the Clenshaw-Curtis quadrature weights. Further we fix
$$ \beta=0.1,\,  T=20,\, \sigma(x)=\ln(1+\exp(x)), \text{ and } N=14.$$
Our network-based feedback laws for~\eqref{eq:burgersapprox}  are obtained based on a~four layer neural network on the set~$Y_0=[-3,3]^{13}$. A set of~$40$ random training initial conditions is randomly sampled in~$Y_0$. We apply equal weights~$w_i=1/40$ and first train a neural network controller for~$\delta=2$,~$p=1$.
As in the previous example, initial conditions that were not part of the training set are used to compare the neural network feedback with the LQR and PSE controllers.

Turning to the validation of the computed results, Table~\ref{tab:linreac} summarizes the norm of the states and feedback controls as well as the open loop objective function value approximated with a finite time horizon of~$T=50$. The entries~"$+\infty^*$" indicate that the  closed-loop state does not converge to zero while~"$+\infty$" marks finite-time blowups. Additionally we plot the temporal evolution the logarithm of the norm of the state, as well as the feedback controls for two particular initial conditions in Figures~\ref{fig:lincond3} and~\ref{fig:lincond4}.

While all considered feedback laws stabilize~initial conditions~$\mathcal{Y}^3_0$ and~$\mathcal{Y}^1_0
$ increasing the magnitude of the latter one, cf.~$\mathcal{Y}^2_0$, causes blowups of the LQR and PSE closed-loop states. This emphasizes the local nature of these controllers. Moreover the nonlinear PSE feedback is unable to stabilize~$\mathcal{Y}^3_0$ and leads to the largest open loop objective functional values for all initial conditions. On the contrary the neural network feedback stabilizes the dynamical system at an exponential rate in all test cases and admits the lowest open loop objective functional value among the considered controllers. Such a behavior can be expected since the neural network controller is found as an approximate solution to a minimization problem involving the nonlinear closed-loop system. We also point to the difference in magnitude between the LQR optimal controls and the neural network feedback which is a well known drawback of linearization based feedback laws. It is also reflected in the smaller norm of the nonlinear feedback controls compared to the linear ones. Last we highlight the remarkably different transient behavior of the neural network feedback in Figure~\ref{fig:lincond3} especially in comparison to the PSE controller.

\begin{table}[!htb]
\begin{minipage}{.48\linewidth}
    \centering

    \medskip
\begin{tabular}{ c c c c }

\multicolumn{4}{c}{$\mathcal{Y}^1_0(x)=\cos(2\pi x)cos(\pi x)+0.5$}\\
\toprule
$\mathcal{F}$ &  $\|Qy\|_{L^2}$ & $\|\F(y)\|_{L^2}$ &  $J(y,\mathcal{F}(y))$  \\
\midrule
uncont. & $+\infty^*$ & 0 & $+\infty^*$\\
LQR & 1.0 & 5.28 & 1.9\\
PSE & 1.38 & 4.56 & 1.99\\
NN  & 1.12 & 4.86 & 1.81 \\
\bottomrule
\end{tabular}
\end{minipage}\hfill
\begin{minipage}{.48\linewidth}
    \centering

    \medskip
\begin{tabular}{ c c c c }
\multicolumn{4}{c}{$\mathcal{Y}^2_0(x)=\cos(2\pi x)cos(\pi x)+1.5$}\\
\toprule
$\mathcal{F}$ &  $\|Qy\|_{L^2}$ & $\|\F(y)\|_{L^2}$ &  $J(y,\mathcal{F}(y))$  \\
\midrule
uncont. & $+\infty^*$ & 0 & $+\infty^*$\\
LQR & $+\infty$& $+\infty$& $+\infty$ \\
PSE & $+\infty$& $+\infty$& $+\infty$\\
NN  & 2.34 & 11.6 & 9.47\\
\bottomrule
\end{tabular}
\end{minipage}\hfill
\begin{minipage}{.48\linewidth}
    \centering

    \medskip

\begin{tabular}{ c c c c }
\multicolumn{4}{c}{$\mathcal{Y}^3_0(x)=-2\operatorname{sign}(x)$}\\
\toprule
$\mathcal{F}$ &  $\|Qy\|_{L^2}$ & $\|\F(y)\|_{L^2}$ &  $J(y,\mathcal{F}(y))$  \\
\midrule
uncont. & $+\infty^*$ & $0$ & $+\infty^*$\\
LQR & 2.67 & 5.5 & 5.08  \\
PSE & 3.49 & 7.04 & 8.57\\
NN  & 2.5 & 1.36 & 3.23\\
\bottomrule
\end{tabular}
\end{minipage}\hfill
\begin{minipage}{.48\linewidth}
    \centering

    \medskip

\begin{tabular}{ c c c c }
\multicolumn{4}{c}{$\mathcal{Y}^4_0(x)=2.5 (x-1)^2 (x+1)^2$}\\
\toprule
$\mathcal{F}$ &  $\|Qy\|_{L^2}$ & $\|\F(y)\|_{L^2}$ &  $J(y,\mathcal{F}(y))$  \\
\midrule
uncont. & $+\infty^*$ & 0 & $+\infty^*$\\
LQR & 1.85 & 13.1 & 10.3 \\
PSE & $+\infty$ & $+\infty$ & $+\infty$\\
NN  & 1.93 & 11.9 & 8.94 \\
\bottomrule
\end{tabular}
\end{minipage}
\caption{Results for~$\delta=2$,~$p=1$.}
\label{tab:linreac}
\end{table}

\begin{figure}[htb]
\centering
\begin{subfigure}[t]{.45\linewidth}
\scalebox{.48}
{\input{costlincont3.tex}}\\
\caption{Running cost vs. time~$t$}
\end{subfigure}
\begin{subfigure}[t]{.45\linewidth}
\scalebox{.48}
{
%
%
\begin{tikzpicture}

\begin{axis}[%
width=4.667in,
height=3.681in,
at={(0.783in,0.497in)},
scale only axis,
xmin=0,
xmax=15,
xlabel style={font=\color{white!15!black}},
xlabel={$t$},
ymin=-20,
ymax=5,
ylabel style={font=\color{white!15!black}},
ylabel={$\log(|Qy(t)|^2)$},
axis background/.style={fill=white},
legend style={legend cell align=left, align=left, draw=white!15!black}
]
\addplot [color=black, dashdotted]
  table[row sep=crcr]{%
0	2.07430014217942\\
0.007764162435695	1.99280895809891\\
0.0208605086034037	1.94497357313253\\
0.0386724908874641	1.91343560819034\\
0.065982185619891	1.88532412519834\\
0.134853234485716	1.82105653449592\\
0.172112308942598	1.7705586473965\\
0.213581305997659	1.69845016641222\\
0.266508729444698	1.58734903861207\\
0.383187883648747	1.31513972077753\\
0.477473365505906	1.10802355215123\\
0.558831719735185	0.949505405032241\\
0.639020684207814	0.811648333289783\\
0.720647569187825	0.688103284539594\\
0.804728692756603	0.576194508360127\\
0.891806407729483	0.474425345178895\\
0.982004175919627	0.382105230363564\\
1.07494352177869	0.299167820103076\\
1.17022266004444	0.225573519993143\\
1.26698777353367	0.161611037118679\\
1.36436793839345	0.107469503510085\\
1.46159233934821	0.0631978226846801\\
1.55766207360868	0.0288386404205632\\
1.65240939706396	0.0040540382852452\\
1.74539777797629	-0.0113967855965829\\
1.8368404682119	-0.0178590335303248\\
1.9268869071456	-0.0155685187416825\\
2.01640193730342	-0.00458420410217464\\
2.10560945432345	0.0151696140240869\\
2.19503317297617	0.0439141477879659\\
2.28586120551348	0.0823162210409372\\
2.37873650699892	0.131106255618429\\
2.4747424059939	0.191433239645498\\
2.57530226129024	0.264953591036813\\
2.68255170442376	0.3542310191597\\
2.79913156704712	0.462694297418812\\
2.93046053385729	0.596902749462867\\
3.08983996890843	0.772550851400737\\
3.3806038224208	1.10874905796883\\
3.58508073560128	1.33826587989916\\
3.72938515640221	1.48784845959117\\
3.85442568691514	1.60553690523413\\
3.96986416055293	1.70270340023737\\
4.08047231565781	1.78477393548363\\
4.18922568980054	1.85489653282691\\
4.29883487782778	1.91539576377391\\
4.41144293825095	1.96773470239734\\
4.52962342618861	2.01313771219489\\
4.65604080734652	2.05241636288961\\
4.79416148064333	2.08622198859738\\
4.94843437207393	2.11500468238807\\
5.12528988458754	2.13911066340613\\
5.334613531345	2.15879045932889\\
5.59309216826716	2.17422708666389\\
5.93241619795384	2.1855585121956\\
6.42634248290805	2.19293746204341\\
7.32013362537621	2.1966129475167\\
10.6268545297292	2.19735947912748\\
15	2.1973605171166\\
};
\addlegendentry{uncont.}

\addplot [color=blue, line width=1.5pt]
  table[row sep=crcr]{%
0	2.07430014217942\\
0.00670162130014873	2.02100549725326\\
0.0150392499156702	2.01053737171062\\
0.0245530215941834	2.01547994017825\\
0.0417974000634764	2.0403657110417\\
0.122191455535873	2.17098823940199\\
0.147130394424657	2.19252166153143\\
0.169037901058942	2.2006712462571\\
0.189941363650398	2.19889319224836\\
0.212551388555902	2.18702789257283\\
0.240456085827752	2.15988087062359\\
0.273343838777201	2.11374483906414\\
0.316558711915864	2.03700493262108\\
0.382391357009183	1.90110189279422\\
0.670895392049697	1.28787007918762\\
0.786305352125886	1.06976103977641\\
0.902612952873145	0.867495156138574\\
1.02562385198657	0.670448796174128\\
1.15597588213118	0.477653110417023\\
1.29722888500671	0.284010062406473\\
1.45129908604723	0.0872970066229861\\
1.62432963412025	-0.119625761300057\\
1.82135592326662	-0.341828743910387\\
2.05872903018985	-0.596391573847619\\
2.37380703326487	-0.921082679633869\\
3.55158018809662	-2.12342542860353\\
3.89283122699069	-2.49005425026101\\
4.21533136661758	-2.84886180106074\\
4.52990189737373	-3.21143507230151\\
4.84291672946687	-3.58501984202281\\
5.16144222004936	-3.97832835103121\\
5.48882066497983	-4.39604430587381\\
5.82970091832749	-4.84479079117717\\
6.18731117446895	-5.32959883424675\\
6.56767015694616	-5.85950251650609\\
6.97801650161724	-6.44565731323095\\
7.42947808026117	-7.10529141197839\\
7.93804041768088	-7.8634596569073\\
8.52332084426188	-8.75135025666607\\
9.22433206443939	-9.83043874704277\\
10.114690572765	-11.2170386271421\\
11.3539763737645	-13.1635959249266\\
13.44500367365	-16.4655737842219\\
15	-18.9255884782786\\
};
\addlegendentry{LQR}

\addplot [color=green, line width=1.5pt]
  table[row sep=crcr]{%
0	2.07430014217942\\
0.00670170321618713	2.01315424722737\\
0.0150392430506265	1.9985308193178\\
0.0245533009624541	2.00207549069208\\
0.0812393164403566	2.05642807769143\\
0.0954373377896083	2.05192897086168\\
0.112018829003976	2.03228383550104\\
0.129163850479419	1.99561317706062\\
0.152418390455857	1.92329031760763\\
0.27108253022217	1.51021670113416\\
0.325863809612422	1.37498681067638\\
0.370695750926393	1.28386038303541\\
0.407201689768906	1.22678813263188\\
0.440959368694424	1.18878478716525\\
0.476199315474197	1.16253827414886\\
0.515593778280623	1.14504810800763\\
0.572872001772417	1.13085477773612\\
0.703041466694451	1.10063699401595\\
0.822787634872519	1.06208849839832\\
1.09049503436995	0.972272592922206\\
1.63077957750563	0.809837141641868\\
2.32087236858069	0.588312291963517\\
3.1449431169146	0.323877141945911\\
3.49954149410244	0.219302750614407\\
3.78793497315879	0.142834727656846\\
4.04576587279305	0.0830106581983436\\
4.29609235107076	0.0335437864502843\\
4.55897684035358	-0.00966616743559179\\
4.89826359155634	-0.056149027944473\\
5.33554041242896	-0.11693003045904\\
5.54530718789612	-0.155055271745411\\
5.71807310791466	-0.194794669960018\\
5.87503667765612	-0.239484457796483\\
6.02099087834408	-0.289838433104402\\
6.16067901926098	-0.347097030676156\\
6.29467426625742	-0.411223792417577\\
6.42469780734037	-0.482731322716726\\
6.55431649167847	-0.563634289165567\\
6.68486713235666	-0.655217876431182\\
6.81529334981874	-0.757095720078917\\
6.94780485699414	-0.871359229520248\\
7.08288196769978	-0.998985126022697\\
7.22296633072473	-1.14306818063666\\
7.36907900163064	-1.30572258601286\\
7.52150272052583	-1.48829280377849\\
7.6813284613671	-1.69303163872192\\
7.85234234489195	-1.92600280348776\\
8.03510175565824	-2.18924902761895\\
8.23621477369433	-2.49374460914216\\
8.46356290113	-2.85346150084662\\
8.73138134466239	-3.29346842638446\\
9.0715000552907	-3.86931248195311\\
9.64111116012161	-4.85293927525334\\
10.4321911498155	-6.21540289581336\\
10.995484721173	-7.16812277680448\\
11.5759638184843	-8.13322593402822\\
12.2279789633944	-9.20062847162297\\
13.0169457200315	-10.4755971480429\\
14.0657939394602	-12.1537213583731\\
15	-13.6399516455854\\
};
\addlegendentry{PSE}

\addplot [color=red, line width=1.5pt]
  table[row sep=crcr]{%
0	2.07430014217942\\
0.0074559525330784	1.99884062579611\\
0.0193251223581079	1.95839049129526\\
0.0347650897767124	1.93428716995631\\
0.0602945486200461	1.91419192656673\\
0.122256206758721	1.87336524085414\\
0.155041880213769	1.83847672772337\\
0.190190821854504	1.78764588455945\\
0.233690034078137	1.70814891770226\\
0.29754978833369	1.57122518109257\\
0.429717750592562	1.28565431758764\\
0.501324736805426	1.15245017143205\\
0.571195797079053	1.03974229330575\\
0.642945626771414	0.939600486797747\\
0.720144779639767	0.846161714275691\\
0.806761145076717	0.754806737358411\\
0.906870896142088	0.661982425574003\\
1.02638699176849	0.563148255391361\\
1.17777637359958	0.449353647908438\\
1.3908332243676	0.300562708364932\\
1.75639195600223	0.0570768091832221\\
2.43857042386368	-0.397480426780096\\
2.81805731498022	-0.661129831516842\\
3.15839038643408	-0.907696041281323\\
3.47642520053222	-1.14842075717995\\
3.77983674657101	-1.3885692436465\\
4.07253224775628	-1.63088095898213\\
4.356955122206	-1.87709820059752\\
4.63854768295823	-2.13190098469476\\
4.91677131566229	-2.39490167473474\\
5.19601687418132	-2.67045641866529\\
5.47447379696409	-2.95702636735168\\
5.75671769867131	-3.25963660568299\\
6.04186259510469	-3.57778652409662\\
6.32927781876887	-3.91101847944432\\
6.62559945783596	-4.26756061811168\\
6.92676513017752	-4.64311969346763\\
7.23716231328371	-5.04366777382385\\
7.55699990709465	-5.47015397915993\\
7.88854347193925	-5.92630637604275\\
8.23308959061328	-6.41470088426468\\
8.5920569898874	-6.93815569315001\\
8.96691514750597	-7.49963307006498\\
9.36128073992272	-8.10549331663739\\
9.77680566893383	-8.7593275940623\\
10.2153710043639	-9.4651401400177\\
10.678732640119	-10.2267837427717\\
11.1697530581907	-11.0500151040603\\
11.6933169506923	-11.944235273723\\
12.2507126926826	-12.9129271948188\\
12.8477980016552	-13.9676560553581\\
13.4751996889804	-15.0929280386392\\
14.1527738738526	-16.3256859070542\\
14.9182439129245	-17.737958432217\\
15	-17.8899487616417\\
};
\addlegendentry{NN}

\end{axis}
\end{tikzpicture}%
}\\
\caption{State norm vs. time~$t$}
\end{subfigure}
\begin{subfigure}[t]{.45\linewidth}
\centering
\scalebox{.48}
{\input{ulincont3.tex}}\\
\caption{Feedback controls}
\end{subfigure}
\caption{$\delta=2$,~$p=1$,~$\mathcal{Y}_0=-2 \operatorname{sign}(x)$}
\label{fig:lincond3}
\end{figure}
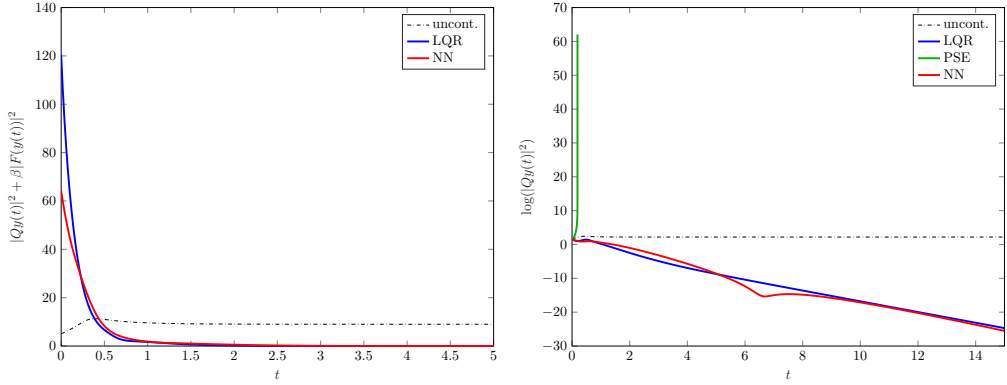
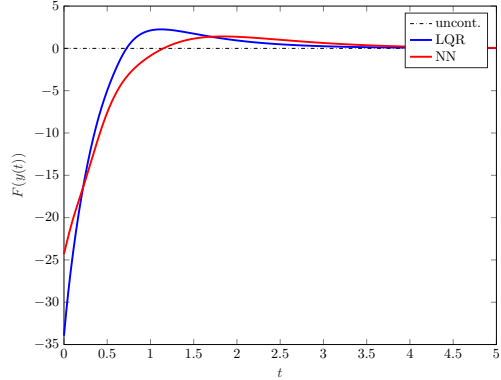
\begin{figure}[htb]
\centering
\begin{subfigure}[t]{.45\linewidth}
\scalebox{.48}
{
%
%
\begin{tikzpicture}

\begin{axis}[%
width=4.667in,
height=3.681in,
at={(0.783in,0.497in)},
scale only axis,
xmin=0,
xmax=5,
xlabel style={font=\color{white!15!black}},
xlabel={$t$},
ymin=0,
ymax=140,
ylabel style={font=\color{white!15!black}},
ylabel={$|Qy(t)|^2+\beta|F(y(t))|^2$},
axis background/.style={fill=white},
legend style={legend cell align=left, align=left, draw=white!15!black}
]
\addplot [color=black, dashdotted]
  table[row sep=crcr]{%
0	5.07936507936508\\
0.027415065779504	5.48904368119995\\
0.056899315240404	5.97191911878231\\
0.0879511696196627	6.52696797279962\\
0.12114246552588	7.17129292640425\\
0.159040293449452	7.96348885792154\\
0.25900426170233	10.0874564698495\\
0.279916552354058	10.4518558652271\\
0.297445989961135	10.7112656716889\\
0.312800100197826	10.8992863858659\\
0.326551130169113	11.0352394134445\\
0.339010682790002	11.1322089569216\\
0.350332230144481	11.1996405491269\\
0.360588592001488	11.2448227518434\\
0.369835073613588	11.2735968605196\\
0.378078371258592	11.2904889052978\\
0.385328456431555	11.2991038783689\\
0.391817561279375	11.3022774714018\\
0.398103662283397	11.3016046144751\\
0.404792398146581	11.2971884592686\\
0.41249460156563	11.2878354694496\\
0.421667856544634	11.2714300367502\\
0.432827968126597	11.2448796257497\\
0.446513349550425	11.2042968903776\\
0.463822801872892	11.1434376299297\\
0.4873175871704	11.0496563977811\\
0.528374114492731	10.8721227113924\\
0.5869428163727	10.6215022144473\\
0.624675995737144	10.4725005570869\\
0.65929024542314	10.3469646336454\\
0.692501084892143	10.2368399053184\\
0.725096376455557	10.1382663519821\\
0.757603858490922	10.0487571119865\\
0.790259250497821	9.96698135281519\\
0.823138142306266	9.89213821045366\\
0.85641560559228	9.82328857207045\\
0.890261826622499	9.75964630785507\\
0.924850410733431	9.70054801987582\\
0.960364721591132	9.64543137131219\\
0.996695865657353	9.59423144975462\\
1.03402495731494	9.54647375787925\\
1.07255261564507	9.50175460618821\\
1.1123315605581	9.4599034269064\\
1.15359680785815	9.42059831537033\\
1.19642581372689	9.38372282798193\\
1.24090166594044	9.34916792121531\\
1.28732010824836	9.3166948214766\\
1.33580736119778	9.28623284762536\\
1.38672870279943	9.25759523543738\\
1.44027357036688	9.23074189668684\\
1.49690429683111	9.20553017753323\\
1.5569050874819	9.18194453725963\\
1.62060756874725	9.15996761237825\\
1.68839729753507	9.1395814229119\\
1.76134480087892	9.12061766770338\\
1.83978675291709	9.1031619750394\\
1.92515507026759	9.0870883818408\\
2.01813560516952	9.072478800664\\
2.12078482917312	9.05924017057124\\
2.23490086570365	9.04740434437814\\
2.36334734781819	9.03695586228633\\
2.51002889473799	9.02789033515816\\
2.68114324626446	9.0201830842734\\
2.88669967481442	9.01381214187433\\
3.14210795708518	9.00879361463523\\
3.48019881879674	9.00508562060762\\
3.97598935292351	9.00266337785296\\
4.87880550444448	9.00146239049697\\
5	9.00141121127081\\
};
\addlegendentry{uncont.}

\addplot [color=blue, line width=1.5pt]
  table[row sep=crcr]{%
0	120.377319048264\\
0.0169637275151189	107.291069919764\\
0.0357434072654144	94.8515113568984\\
0.0567627715029317	82.9164642982501\\
0.0784530366271099	72.3426596009596\\
0.0983334046170228	63.915080896129\\
0.119460653380969	56.0641328894469\\
0.140900103761666	49.0884475137475\\
0.162298160382491	42.9839613632902\\
0.184355411286305	37.4725715357088\\
0.206592435249689	32.6217077485758\\
0.228076821617165	28.5309571743841\\
0.24940911585918	24.9859237206374\\
0.270073719373173	21.9919618285109\\
0.291668069870198	19.2777910171936\\
0.312334176191229	17.0347486283998\\
0.33283411929564	15.1139700224946\\
0.35185792779626	13.5714751765364\\
0.371220658008014	12.210006229318\\
0.391754315708965	10.9657989916554\\
0.41054378970388	9.98147641762225\\
0.430346724072024	9.07817365079009\\
0.451732810340729	8.22961783141521\\
0.474891811042255	7.4278393042993\\
0.499959780832143	6.66240169179666\\
0.52887817323483	5.87211161721227\\
0.56073908703344	5.08118401266516\\
0.591512675366658	4.38374987525972\\
0.617951814546828	3.84371937821861\\
0.640072288488383	3.44399440002842\\
0.659799653483134	3.13397013574381\\
0.677516159228759	2.89540038829539\\
0.693160970958175	2.71608063171637\\
0.710194078826646	2.55257746972359\\
0.725603140083848	2.43056832239964\\
0.741227955519406	2.32875514515779\\
0.756051158027944	2.24939378974115\\
0.771792975827438	2.18008054219459\\
0.786984822242559	2.12472564759518\\
0.8045276068883	2.07138277958087\\
0.823911962903338	2.02170627654988\\
0.848165792472045	1.96801104018853\\
0.883620705171879	1.89728982370805\\
0.957551574285375	1.75075730057802\\
1.00410785984333	1.65154224666477\\
1.05938181321225	1.52696773235769\\
1.16381757401639	1.2827653426988\\
1.2399598556763	1.10846559425396\\
1.29661923083381	0.98566811218538\\
1.34425934367441	0.888617362707478\\
1.38965401759623	0.802130280746027\\
1.43325803399293	0.724889809936329\\
1.47660556175694	0.653885328325785\\
1.51713292625737	0.592682488552384\\
1.55779623822686	0.536183781224295\\
1.59886914044053	0.48391719074877\\
1.63894577691391	0.437340498712061\\
1.68043003400759	0.393458234797507\\
1.72087307253645	0.354638112313381\\
1.7620621344833	0.318828056321536\\
1.80477065776975	0.28535402150797\\
1.84785707635611	0.255034762376141\\
1.89266974868816	0.226840556981486\\
1.93913596732396	0.200850073514616\\
1.98787161095632	0.176766065266051\\
2.03918534117118	0.154526461424425\\
2.09305735565891	0.134205035702365\\
2.15027471672316	0.11557859470021\\
2.21293763080109	0.0981863105964322\\
2.27918204739602	0.0827033927804308\\
2.3508481417073	0.0687663855915162\\
2.43015686295489	0.0561502159395246\\
2.51895446078727	0.0448467593653561\\
2.61916972336911	0.034903180405351\\
2.73597074357454	0.0261741696890283\\
2.87319731441372	0.018784239755874\\
3.03981998858954	0.0126784108419713\\
3.25067195220552	0.00783078995942788\\
3.53656700792756	0.00418861229299239\\
3.96756717147818	0.00172456089799766\\
4.76278714921625	0.000385985484527396\\
5	0.000253434964363919\\
};
\addlegendentry{LQR}

\addplot [color=red, line width=1.5pt]
  table[row sep=crcr]{%
0	64.1759461203391\\
0.0184771428873063	59.9044237700135\\
0.0457170935585225	54.2346475859262\\
0.0715681637595083	49.4158276190396\\
0.0966378377681281	45.2661808670201\\
0.120870363168905	41.6938192911296\\
0.149539044739768	37.8959265956779\\
0.1816638249324	34.01478985682\\
0.220889070572383	29.6044552305825\\
0.264512392183917	25.0058602042226\\
0.301522563336789	21.3770537173252\\
0.334355740234642	18.4222565371604\\
0.361816217018955	16.1764204739457\\
0.388040553163904	14.2410981416086\\
0.413762709628472	12.5477263641221\\
0.436853147238608	11.1990551142441\\
0.4600958503692	9.99863962914108\\
0.483564285615842	8.93669882130546\\
0.506226650359409	8.04236098093132\\
0.528512554709096	7.27588523304043\\
0.549823233973939	6.63619895952375\\
0.571449301099264	6.06874161025232\\
0.593551969374417	5.56254298706966\\
0.615945072075945	5.11476929306494\\
0.637578565455769	4.73494847330871\\
0.66001500020235	4.38711640412059\\
0.684065379545927	4.0578690504346\\
0.708910802052145	3.75698242519543\\
0.734686526319123	3.47958699586364\\
0.760126795057843	3.23485508926581\\
0.785047429308477	3.0188304197863\\
0.812265525291124	2.8061346831259\\
0.838640719227186	2.6205065382257\\
0.867015078557117	2.44100341136266\\
0.894283074154814	2.28651798620972\\
0.920097011313061	2.15508626983227\\
0.946781099067934	2.03326338702452\\
0.973595694020474	1.92389418505624\\
0.999138740931429	1.83083915477616\\
1.02640853416537	1.74224200632048\\
1.05195103092004	1.66825312141644\\
1.07861101988854	1.59920791848334\\
1.10542951256205	1.53714214386568\\
1.13351910749546	1.47897373119706\\
1.16226198022541	1.4256191893185\\
1.19357213938986	1.37341023409829\\
1.22712060697492	1.32304649996783\\
1.26321187083184	1.27395096142794\\
1.3032100071575	1.22422240965597\\
1.3490120458404	1.17164484638636\\
1.40765584182276	1.10881588276783\\
1.48108312771186	1.03416525437055\\
1.61320722472387	0.904519733778272\\
1.77242693850717	0.752290950780818\\
1.8674203456238	0.664836714135276\\
1.95132343019087	0.591040720729666\\
2.02797705448721	0.52725024186617\\
2.09641938963399	0.473690426807423\\
2.16297642349171	0.424941939893387\\
2.23382674909524	0.37686586636174\\
2.29974481962383	0.335711630856679\\
2.36135012809315	0.30035351749234\\
2.42537471480951	0.266727768493851\\
2.48924853001131	0.236239808957464\\
2.55536867252412	0.20774289667456\\
2.61954422633882	0.182881083872999\\
2.69064352099063	0.158343628881852\\
2.76158761125234	0.136754395804331\\
2.83820531247203	0.116384924990598\\
2.91654374486373	0.0983936079161936\\
3.00242687694417	0.0815795153693557\\
3.09376080333756	0.0666041522066223\\
3.19350079959553	0.0531609966863584\\
3.3008921299741	0.041523623645972\\
3.42118917290028	0.0313269177779887\\
3.55651491648082	0.0226788722152946\\
3.71411908627411	0.0154485027884164\\
3.90353070234042	0.00963623011460868\\
4.14002945916154	0.00526071232329173\\
4.4629538348141	0.00223557756149262\\
4.97420995854725	0.000533630061013923\\
5	0.000494981180239051\\
};
\addlegendentry{NN}

\end{axis}
\end{tikzpicture}%
}\\
\caption{Running cost vs. time~$t$}
\end{subfigure}
\begin{subfigure}[t]{.45\linewidth}
\scalebox{.48}
{
%
%
\begin{tikzpicture}

\begin{axis}[%
width=4.667in,
height=3.681in,
at={(0.783in,0.497in)},
scale only axis,
xmin=0,
xmax=15,
xlabel style={font=\color{white!15!black}},
xlabel={$t$},
ymin=-30,
ymax=70,
ylabel style={font=\color{white!15!black}},
ylabel={$\log(|Qy(t)|^2)$},
axis background/.style={fill=white},
legend style={legend cell align=left, align=left, draw=white!15!black}
]
\addplot [color=black, dashdotted]
  table[row sep=crcr]{%
0	1.62518626940224\\
0.220413727057737	2.23050941002925\\
0.263998932376438	2.32035619616855\\
0.298133706218579	2.37215704119858\\
0.328080878803492	2.40229203953352\\
0.356855063711617	2.41859331239291\\
0.38720602352789	2.42484219158457\\
0.42305526331425	2.42201112035116\\
0.472581594926433	2.40781447477748\\
0.616187266242981	2.35184676684659\\
0.726022802539735	2.31605367842048\\
0.840728492246143	2.2879714468914\\
0.974665903680384	2.26433119502907\\
1.13959270826811	2.24426704116924\\
1.35005990426609	2.22763408995515\\
1.62884485383693	2.21455509913919\\
2.02076476633015	2.20520425038989\\
2.64500015621853	2.1996216076218\\
4.04320675610736	2.19750044064581\\
15	2.19736051728935\\
};
\addlegendentry{uncont.}

\addplot [color=blue, line width=1.5pt]
  table[row sep=crcr]{%
0	1.62518626940224\\
0.027620023978514	1.43082247801699\\
0.0567627715029353	1.28417749571264\\
0.0860310000966216	1.17962438046112\\
0.113762759459398	1.11049099727122\\
0.138257937597334	1.06956385411864\\
0.16036048299717	1.047135406724\\
0.181806543705367	1.03758268937796\\
0.20232554521759	1.03897866009762\\
0.223808798103146	1.05063275204758\\
0.249409115859184	1.07659896119275\\
0.280452596854488	1.12236987933782\\
0.326658800230838	1.20852348468716\\
0.405646309098454	1.35612693557808\\
0.44033015085628	1.40249058555494\\
0.466369866940095	1.4242331170673\\
0.488525726717835	1.43190582820737\\
0.507418297793794	1.4295424157276\\
0.527018272118536	1.4176851511989\\
0.548108848456547	1.39370347617449\\
0.571545222191233	1.35332872848051\\
0.599141482649447	1.28852048359932\\
0.637681779971043	1.17386880788833\\
0.762431802923047	0.783447535631968\\
0.832878082721113	0.606223009096173\\
0.994271066041968	0.205529996021024\\
1.1714937926485	-0.264595704372592\\
2.05673996846576	-2.63496445474841\\
2.30994326827832	-3.27851192271783\\
2.54680576470417	-3.85805072323449\\
2.7755806130941	-4.39601762942644\\
3.00199932111327	-4.90740410191542\\
3.23047510622363	-5.40303551447985\\
3.46351420311966	-5.8888647838724\\
3.7072931282093	-6.37775222633969\\
3.96324162014096	-6.87225220081468\\
4.23620275183726	-7.38141478736692\\
4.53466719242365	-7.92027060790765\\
4.86503694967437	-8.49928010413437\\
5.24249974083738	-9.14363491375972\\
5.68655207352325	-9.88466737312867\\
6.23101059099803	-10.776418141745\\
6.94240728712835	-11.9247307982696\\
7.96654664292643	-13.5607793443806\\
9.74603686621248	-16.3856538832461\\
14.0562478572611	-23.2101376003013\\
14.1696454367244	-23.3900789562048\\
14.685849239428	-24.2057651384998\\
14.7013301504165	-24.2316584854537\\
14.9908079975187	-24.6879632306033\\
14.9956500621871	-24.6972110200161\\
15	-24.7030185914873\\
};
\addlegendentry{LQR}

\addplot [color=green, line width=1.5pt]
  table[row sep=crcr]{%
0	1.62518626940224\\
0.0215885496664328	1.74777504941456\\
0.0433151968979786	1.92129658128176\\
0.0641317441289289	2.15351589552004\\
0.0862085376067796	2.49904623041804\\
0.107699280673891	2.97276084848353\\
0.128685672120106	3.62227757229198\\
0.148829858735347	4.50948152552005\\
0.166483300828368	5.68684607833633\\
0.179557321101704	7.22159519019607\\
0.186086039180928	8.75017218625634\\
0.188932466090058	12.9631809890727\\
0.189235248910698	33.5227508680903\\
0.189235250094072	62.0824659350554\\
};
\addlegendentry{PSE}

\addplot [color=red, line width=1.5pt]
  table[row sep=crcr]{%
0	1.62518626940224\\
0.0325313581039808	1.4780875381803\\
0.0765197668495681	1.31894529624148\\
0.128119869560013	1.16322520204484\\
0.174614088093314	1.04629862179944\\
0.213481849346632	0.968809931637605\\
0.246780987097964	0.919348636413037\\
0.276887417156026	0.888553191668883\\
0.305963931580635	0.870786682690248\\
0.334355740234642	0.863582395655062\\
0.366968447937641	0.865543801170173\\
0.403555238360301	0.877333566465886\\
0.464975842197589	0.908560721430327\\
0.534210952267884	0.941271654909258\\
0.577260980520389	0.951845143903221\\
0.615945072075952	0.951897428022079\\
0.652477742511959	0.942721297382601\\
0.690764159780041	0.923328595685952\\
0.732364805157239	0.891651536466526\\
0.782430470379794	0.841281433516169\\
0.849779262589948	0.759344325095885\\
1.32084521802652	0.143109790682921\\
1.44023316161549	-0.0312205151690783\\
1.56505406937489	-0.228677523264345\\
1.70186975067375	-0.461172028241407\\
1.85425256138774	-0.737039918829893\\
2.02277451986904	-1.05974479651374\\
2.21035051149865	-1.43747828978265\\
2.41499440394821	-1.86864873379876\\
2.63948708214611	-2.36136357917725\\
2.8811490690937	-2.91201238155113\\
3.13720222445923	-3.51622864015808\\
3.40304416183145	-4.16496925417995\\
3.67095862082984	-4.84071481228058\\
3.93514013858026	-5.5295054006257\\
4.19131626336196	-6.22054796473162\\
4.43595207235908	-6.90439360373994\\
4.6652734554176	-7.57007553820509\\
4.87854829552244	-8.21462016677868\\
5.07485784341496	-8.8341029963829\\
5.25450871467593	-9.42795793183866\\
5.42129058987926	-10.007710096799\\
5.57505565365318	-10.572307462052\\
5.71685608848991	-11.1247901957431\\
5.84643579806639	-11.6627179258165\\
5.96829415769279	-12.2038400590681\\
6.0841092451551	-12.7557785361685\\
6.20172051797546	-13.3578056583555\\
6.35817212261903	-14.2086732660142\\
6.45871310114702	-14.7353638023131\\
6.51396609143203	-14.9830037226868\\
6.55641157717177	-15.1369800773499\\
6.58952581207085	-15.2294774403119\\
6.61813980235209	-15.2880342360106\\
6.64309297249018	-15.3227818907761\\
6.66701737701058	-15.3424415441714\\
6.68992064258289	-15.3497479612578\\
6.71355369573347	-15.3468512957264\\
6.7396238608104	-15.333242813\\
6.76933928734123	-15.3071449507388\\
6.81037481664816	-15.2580520658445\\
6.87506885292767	-15.1646712001536\\
7.00470884639032	-14.9763577491403\\
7.07482074726161	-14.8905916380245\\
7.14015457355641	-14.8241051348659\\
7.20318222907952	-14.7722566040993\\
7.26516406357751	-14.7324257539228\\
7.32637487683766	-14.7031436166951\\
7.38821701388989	-14.6828423270686\\
7.45245924579982	-14.6707320516059\\
7.51935771367465	-14.6669155970973\\
7.5897532028501	-14.6716181330175\\
7.66569724721431	-14.6856802034875\\
7.74495029820428	-14.7092421987539\\
7.83210490854199	-14.744476961216\\
7.92463192766004	-14.7913581761699\\
8.02502825533769	-14.8519825847777\\
8.13563581251439	-14.9291279376544\\
8.25548703413944	-15.02345471895\\
8.38703713846893	-15.1382204270726\\
8.52979667512003	-15.2743394086369\\
8.68472265848453	-15.4338980633836\\
8.85446825542589	-15.6209669232433\\
9.03902436516256	-15.8369222963467\\
9.24833754950859	-16.095497625811\\
9.47286814111269	-16.3867174235059\\
9.72018331644886	-16.7217092042292\\
9.98086112999171	-17.0886370410495\\
10.2503027686016	-17.4807490304656\\
10.6029322979565	-18.0108572370991\\
10.9256973249627	-18.5104726035562\\
11.3398872757289	-19.1688363994177\\
11.7698223072295	-19.8697232587506\\
12.2212944702585	-20.6221622565403\\
12.7363965232308	-21.4982626178354\\
13.4395380643961	-22.7199785421839\\
13.5856316855102	-22.9772514131333\\
14.5202008295181	-24.645762361445\\
14.5644084803163	-24.72560539733\\
14.6488685894384	-24.8767954532962\\
14.653710700314	-24.8871162663441\\
14.6926779013654	-24.9555223707052\\
14.6979183507016	-24.9671866207454\\
14.7070402078885	-24.9838991185014\\
14.7872205699146	-25.1291464512026\\
14.7896416253317	-25.1242693154163\\
14.7955244084991	-25.1446366433651\\
14.831428219937	-25.2094360819068\\
14.8340484445831	-25.2055969865927\\
14.84055007705	-25.2262404165856\\
14.9158883276552	-25.3603649281991\\
14.9207304384525	-25.3718277131962\\
14.9231514938511	-25.3644837119804\\
14.9290342769761	-25.387464908951\\
14.9596976388973	-25.4391225568052\\
14.9649380881503	-25.4523282126964\\
14.9675583127768	-25.4462003626888\\
14.9697255235837	-25.4614663641504\\
14.9740599451975	-25.4692353114831\\
14.9762271560044	-25.4646445645573\\
14.9835030623915	-25.4860100701987\\
15	-25.5170210302163\\
};
\addlegendentry{NN}

\end{axis}
\end{tikzpicture}%
}\\
\caption{State norm vs. time~$t$}
\end{subfigure}
\begin{subfigure}[t]{.45\linewidth}
\centering
\scalebox{.48}
{
%
%
\begin{tikzpicture}

\begin{axis}[%
width=4.667in,
height=3.681in,
at={(0.783in,0.497in)},
scale only axis,
xmin=0,
xmax=5,
xlabel style={font=\color{white!15!black}},
xlabel={$t$},
ymin=-35,
ymax=5,
ylabel style={font=\color{white!15!black}},
ylabel={$F(y(t))$},
axis background/.style={fill=white},
legend style={legend cell align=left, align=left, draw=white!15!black}
]
\addplot [color=black, dashdotted]
  table[row sep=crcr]{%
0	0\\
5	0\\
};
\addlegendentry{uncont.}

\addplot [color=blue, line width=1.5pt]
  table[row sep=crcr]{%
0	-33.9555524132503\\
0.0234883535340913	-31.3787781679662\\
0.0492818634165175	-28.8479677750364\\
0.0754867967932924	-26.5204982742978\\
0.105734074115794	-24.0739878707376\\
0.138257937597331	-21.6760334180044\\
0.1715793467894	-19.4275599125058\\
0.204458990233647	-17.3868955588557\\
0.23763171278727	-15.4878480838546\\
0.272118566539099	-13.669963727762\\
0.304809279254769	-12.0849463741626\\
0.338806236951974	-10.5715585224681\\
0.371220658008014	-9.24940801806961\\
0.403192463330669	-8.05266154302642\\
0.436259669928354	-6.91699999922033\\
0.470965601248167	-5.82375888989766\\
0.504940212999585	-4.83740866848484\\
0.541025945229983	-3.8656490532961\\
0.579123153928094	-2.91161955059475\\
0.615864315219774	-2.05691329136512\\
0.650036513843283	-1.32421367108889\\
0.679363709194362	-0.75057100662508\\
0.705671072551624	-0.284685386906027\\
0.731591085109528	0.126497113247957\\
0.756051158027951	0.470534066778924\\
0.779967015438622	0.76686570702838\\
0.802073457602987	1.00757512729039\\
0.823911962903345	1.21637915322161\\
0.845571199724453	1.39740704017668\\
0.8663613748245	1.54894576422613\\
0.886235834017086	1.67530220738579\\
0.907258639490841	1.79106802712\\
0.927081981730915	1.88490060432475\\
0.947264478543474	1.96653307230226\\
0.966782157861722	2.03331296352393\\
0.985715357539796	2.08765694860926\\
1.00410785984332	2.13143821879623\\
1.02108168610895	2.16460335010334\\
1.0369683713188	2.18983693261265\\
1.05259983854837	2.20959163639574\\
1.06692838477592	2.22359422498558\\
1.08129832941498	2.23396940840735\\
1.09485912098307	2.24062034291438\\
1.10777666399538	2.24430961740707\\
1.12295517138162	2.24557358755278\\
1.13769848104994	2.24385454702931\\
1.15355880435393	2.23901872115971\\
1.16893505310447	2.23161466592362\\
1.18787164722588	2.21915373987032\\
1.20693978038863	2.20324574747026\\
1.228894538709	2.18122796809329\\
1.25293692277911	2.1531516744698\\
1.28120491052073	2.11564125101621\\
1.3118852781264	2.07040452955832\\
1.34964079285626	2.00969320614522\\
1.39637602413936	1.92902916115664\\
1.46185396341794	1.81003469010407\\
1.64158374169109	1.480728221432\\
1.7095980992534	1.36330510992364\\
1.77138827573038	1.26210983566713\\
1.83043524119601	1.1707709277152\\
1.88690450258207	1.08851344283479\\
1.94124555029062	1.01409041864723\\
1.99639876637937	0.943241879970614\\
2.05192781613842	0.876553773393937\\
2.10733444428489	0.814484472811081\\
2.16222971251685	0.757195836574247\\
2.21805784095454	0.703010090250473\\
2.27457565751668	0.652106536411765\\
2.3318380340155	0.604339526273733\\
2.38991146449101	0.559554722967967\\
2.44783472065787	0.518304753518606\\
2.50788980013975	0.478883126116095\\
2.57036933527628	0.441213094962222\\
2.63433130262313	0.405899511971064\\
2.70065822630395	0.372462450916835\\
2.76945173790407	0.340904903227056\\
2.84159434199605	0.310912049273476\\
2.91593680833495	0.283005227518323\\
2.99473956534804	0.256411264361887\\
3.076856621327	0.231620927922044\\
3.16425536152519	0.208142965056872\\
3.25699553900485	0.186117726603513\\
3.35656590019213	0.165360896620001\\
3.46351420311966	0.145950430596343\\
3.57955456132814	0.127781351068002\\
3.70547895825168	0.110943374131338\\
3.84230644716027	0.0954817535863768\\
3.9939168860282	0.081183484359137\\
4.16276308693495	0.0680924983231819\\
4.35250465131745	0.0562013449240837\\
4.5701976419467	0.0453998583597723\\
4.82206775232289	0.0357505070082382\\
5	0.0303343030913155\\
};
\addlegendentry{LQR}

\addplot [color=red, line width=1.5pt]
  table[row sep=crcr]{%
0	-24.3097883662063\\
0.0235015500912361	-23.2937449243605\\
0.0642131284567675	-21.6408333653295\\
0.096637837768121	-20.434837653336\\
0.130530819121393	-19.2847635445323\\
0.179376141837697	-17.7381907849722\\
0.250892188339492	-15.4637084240183\\
0.318300726142574	-13.2113649600587\\
0.398451502726218	-10.5507416610655\\
0.440719198418094	-9.24132820690971\\
0.477995197074403	-8.17161003579998\\
0.511401361733462	-7.29046166126597\\
0.542617693195677	-6.53643895045899\\
0.573444656947892	-5.8574685206326\\
0.603710948908134	-5.25223867846001\\
0.633373331199309	-4.71472764376934\\
0.662455848258723	-4.23716735036534\\
0.690764159780041	-3.81504864533637\\
0.721244132858448	-3.40270418913099\\
0.752380758170744	-3.02109592284719\\
0.782430470379794	-2.68541293633749\\
0.816494682454032	-2.33795015955956\\
0.851586845520725	-2.01058232415825\\
0.887110295669149	-1.70527298522794\\
0.922675552213295	-1.42177169840359\\
0.963812996685323	-1.11732722025932\\
1.0037968954563	-0.842696161547266\\
1.04409281507133	-0.585330523219874\\
1.08282061891317	-0.355535059404147\\
1.12120918051922	-0.144158363251211\\
1.1584323122934	0.0455337324410934\\
1.19357213938986	0.210960832580728\\
1.22457495588881	0.346288826116588\\
1.25897118863509	0.484854364083226\\
1.29240512384338	0.608242859586834\\
1.32563374391724	0.720109764321595\\
1.35410500092592	0.807878172631309\\
1.38723639443553	0.900728020440901\\
1.41681856180312	0.97563432266961\\
1.44813366754772	1.04699086230095\\
1.47477905171693	1.10165081527002\\
1.50380069013853	1.15505615344274\\
1.52972106253132	1.19760340292612\\
1.55848852612134	1.23941613757013\\
1.58504400986221	1.27327399159798\\
1.6157579491729	1.30692356824337\\
1.64237324753592	1.33162439999105\\
1.66656235943811	1.35076511423945\\
1.69064696052633	1.36678015877503\\
1.71721681430685	1.38118043453128\\
1.74794252495002	1.39370444654942\\
1.77457261944644	1.40126615344921\\
1.79877090512801	1.40573770234765\\
1.82478203346679	1.40807572015513\\
1.85185241537789	1.40799639537907\\
1.88529853278195	1.40465805889504\\
1.91783899796016	1.39825870008531\\
1.95132343019088	1.38873388594608\\
1.98655801789571	1.37579930140345\\
2.01761789474926	1.36226829203743\\
2.06103553872838	1.3402822501594\\
2.10667825231354	1.31388567106804\\
2.16297642349172	1.27749592360787\\
2.20651650698608	1.24705178168293\\
2.2627562056506	1.20540817782951\\
2.34859537065915	1.13827414972362\\
2.47348333712856	1.03679320302752\\
2.6751243315937	0.874053548467881\\
2.78183025336035	0.792140529733263\\
2.89699691631812	0.708738442381577\\
2.98695194630125	0.647629830062709\\
3.06877018669183	0.595322503905358\\
3.15065631288478	0.546154737902771\\
3.23492127528392	0.498875416024919\\
3.32099538055568	0.453993943193531\\
3.4088132295369	0.411658856087687\\
3.49685926531651	0.372570843600634\\
3.5767327399696	0.339882503754719\\
3.65902921389926	0.308821277271377\\
3.75118212585157	0.277001396155764\\
3.83797011537392	0.249706844611428\\
3.93706452055876	0.221484628539034\\
4.03480563613601	0.196479727116095\\
4.13745694084598	0.172981807984236\\
4.24839279742985	0.150469398350872\\
4.35996166265668	0.130542127597327\\
4.48339544890862	0.111314008868749\\
4.61418723687319	0.0937841426257471\\
4.75671077096386	0.0775724959969679\\
4.91212401237416	0.0628303505304224\\
5	0.0556663902343004\\
};
\addlegendentry{NN}

\end{axis}
\end{tikzpicture}%
}\\
\caption{Feedback controls}
\end{subfigure}
\caption{$\delta=2$,~$p=1$,~$\mathcal{Y}_0=2.5 (x-1)^2 (x+1)^2$}
\label{fig:lincond4}
\end{figure}


The final set of results are obtained for the choice ~$\delta=0.5$ and~$p=3$, while leaving all other specifications unaltered.

The obtained results are depicted in Table~\ref{tab:cubreac}. While the uncontrolled system is stable for~$\mathcal{Y}^1_0$ and~$\mathcal{Y}^3_0$ , the additional cubic term leads to a finite time blow up for larger initial conditions. Since the Fr\'echet derivative of the cubic nonlinearity vanishes at zero the same behavior can be observed for the linearization based feedback. In contrast, the nonlinear controllers take this destabilizing effect into account and drive the state to zero at an exponential rate for all considered initial conditions. It is worth mentioning that the neural network closed loop state is the smallest, in the~$L^2$ sense, in all test cases. This comes at the cost of larger control norms and objective functional values. A possible explanation for this behavior can be found in the early termination of Algorithm~\ref{alg:grad} after only few iterations. However this was particularly crucial in the present example since the use of explicit stepsizes in the gradient method occasionally led to iterates with finite-time blow ups on the training set. A more sophisticated and rigorous algorithmic treatment of systems with exploding behavior in the context of the proposed feedback setting is subject to future research.

\begin{table}[!htb]
\begin{minipage}{.48\linewidth}
    \centering

    \medskip
\begin{tabular}{ c c c c }

\multicolumn{4}{c}{$\mathcal{Y}^1_0(x)=\cos(2\pi x)cos(\pi x)+0.5$}\\
\toprule
$\mathcal{F}$ &  $\|Qy\|_{L^2}$ & $\|\F(y)\|_{L^2}$ &  $J(y,\mathcal{F}(y))$  \\
\midrule
uncont. & 0.98 & 0 & 0.48\\
LQR & 0.67 & 1.19 & 0.26\\
PSE & 0.59 & 1.32 & 0.26\\
NN  & 0.46 & 2.56 & 0.43 \\
\bottomrule
\end{tabular}
\end{minipage}\hfill
\begin{minipage}{.48\linewidth}
    \centering

    \medskip
\begin{tabular}{ c c c c }
\multicolumn{4}{c}{$\mathcal{Y}^2_0(x)=\cos(2\pi x)cos(\pi x)+1.5$}\\
\toprule
$\mathcal{F}$ &  $\|Qy\|_{L^2}$ & $\|\F(y)\|_{L^2}$ &  $J(y,\mathcal{F}(y))$  \\
\midrule
uncont. & $+\infty$ & 0 & $+\infty$\\
LQR & $+\infty$ & $+\infty$ & $+\infty$ \\
PSE & 1.73 & 6.79 & 3.81\\
NN  & 1.53 & 8.03 & 4.4\\
\bottomrule
\end{tabular}
\end{minipage}\hfill
\begin{minipage}{.48\linewidth}
    \centering

    \medskip

\begin{tabular}{ c c c c }
\multicolumn{4}{c}{$\mathcal{Y}^3_0(x)=-2\operatorname{sign}(x)$}\\
\toprule
$\mathcal{F}$ &  $\|Qy\|_{L^2}$ & $\|\F(y)\|_{L^2}$ &  $J(y,\mathcal{F}(y))$  \\
\midrule
uncont. & 1.85 & 0 & 1.72\\
LQR & 1.84 & 0.28 & 1.7  \\
PSE & 1.86 & 0.52 & 1.75\\
NN  & 1.82 & 1.25 & 1.74\\
\bottomrule
\end{tabular}
\end{minipage}\hfill
\begin{minipage}{.48\linewidth}
    \centering

    \medskip

\begin{tabular}{ c c c c }
\multicolumn{4}{c}{$\mathcal{Y}^4_0(x)=2.5 (x-1)^2 (x+1)^2$}\\
\toprule
$\mathcal{F}$ &  $\|Qy\|_{L^2}$ & $\|\F(y)\|_{L^2}$ &  $J(y,\mathcal{F}(y))$  \\
\midrule
uncont. & $+\infty$ & 0 & $+\infty$\\
LQR & $+\infty$ & $+\infty$ & $+\infty$ \\
PSE & 1.86 & 8.01 & 4.94 \\
NN  & 1.59 & 9.2 & 5.49 \\
\bottomrule
\end{tabular}
\end{minipage}
\caption{Results for~$\delta=0.5$,~$p=3$.}
\label{tab:cubreac}
\end{table}

\tcr{
\section{Conclusion}
A feedback strategy based on learning neural networks has been proposed and analysed in the
framework of optimal stabilisation problems. The underlying concept is quite general and several
extensions suggest themselves. It can be of interest to carefully analyse the finite horizon case, with a state-dependent
control operator, so that bilinear control problems occur as a special case. Moreover the treatment of infinite
dimensional controlled systems is of interest.
}
\appendix
\section{Results on Nemitsky operators} \label{app:nemitsky}
In this section we provide several auxiliary results on superposition operators.
\begin{prop} \label{prop:nemitskylip}
Let~$M>0$ and~$F\in \operatorname{Lip}\left( \bar{B}_M(0);\R^{m } \right)$,~$\bar{B}_M(0)\subset \R^n$, with~$F(0)=0$ be given and define the induced Nemitsky operator as
\begin{align*}
\lbrack \F(y)\rbrack(t)= F(y(t)) \quad \forall y \in W_\infty,~\| y\|_{\mathcal{C}_b(I;\R^n)}\leq M,~t \in I.
\end{align*}
Then ~$\F(y)\in L^2(I; \R^m) \cap \mathcal{C}_b(I; \R^m)$ for all~$y \in W_\infty$,~$\| y\|_{\mathcal{C}_b(I;\R^n)}\leq M$, and we have
\begin{align*}
\|\F(y_1)-\F(y_2)\|_{L^2(I;\R^m)} \leq L_{F,M} \|y_1-y_2 \|_{L^2(I; \R^n)},
\end{align*}
as well as
\begin{align*}
\|\F(y_1)-\F(y_2)\|_{\mathcal{C}_b(I;\R^m)} \leq L_{F,M} \|y_1-y_2 \|_{\mathcal{C}_b(I;\R^n)}
\end{align*}
for all~$y_i \in W_\infty$,~$\| y_i\|_{\mathcal{C}_b(I;\R^n)}\leq M$,~$i=1,2$, where $L_{F,M}$ is the Lipschitz constant of $F$ on $\bar{B}_M(0)$.

If~$F$ is continuously Fr\'echet differentiable on~$\bar{B}_M(0)$ and~$DF\in \operatorname{Lip}(\bar{B}_M(0);\R^{m \times n})$, with Lipschitz constant $L_{DF,M}$ on $\bar{B}_M(0)$,  then~$\F:W_\infty$  to~$L^2(I;\R^n)$ is differentiable at~$y\in W_\infty$,~$\|y\|_{\mathcal{C}_b(I;\R^n)}<M$. Its Fr\'echet derivative~$D\F(y)\in \mathcal{L}(W_\infty;L^2(I;\R^m))$ satisfies
\begin{align*}
\lbrack D\F(y)\delta y \rbrack (t)= DF(y(t))\delta y(t) \quad \forall \delta y \in W_\infty
\end{align*}
and almost every~$t\in I$. Moreover we have
\begin{align*}
\|(D\F(y_1)-D\F(y_2))\delta y\|_{L^2(I;\R^m)} \leq \|\imath\|_{\mathcal{L}(W_\infty, \mathcal{C}_b(I;\R^n))} L_{DF,M} \|y_1-y_2\|_{L^2(I;\R^n)} \|\delta y\|_{W_\infty}
\end{align*}
for all~$y_i \in W_\infty$,~$\| y_i\|_{\mathcal{C}_b(I;\R^n)}< M$,~$i=1,2$, and~$\delta y\in W_\infty$.
\end{prop}
For the sake of brevity we leave the proof of this proposition and Lemma~\ref{lem:Nemitsky} to the reader.

\begin{coroll} \label{cor21a}
Let $y_k \rightharpoonup y$ in $W_\infty$, with $\|y_k\|_{\mathcal{C}_b(I;\R^n)}\le M$, and assume that
$F\in \operatorname{Lip}\left( \bar{B}_M(0);\R^{m } \right)$, with~$F(0)=0$.  Then we have
$F(y_k) \rightharpoonup F(y)$ in $L^2(I;\R^n)$.
\end{coroll}
\begin{proof}
Let $\phi \in L^2(I;\R^n)$ with $\phi=0$ on $[T_\phi,\infty)$ for some $T_\phi>0$. By compactness of $H^1(0,T_\phi;\R^n)$ in $L^2(0,T_\phi;\R^n)$ and the first assertion of Proposition \ref{prop:nemitskylip} we deduce
\begin{align*}
 (\phi, F({y}_{k}))_{L^2(I;\R^n)} \rightarrow \tcr{( \phi, F({y}))_{L^2(I;\R^n)}}.
\end{align*}
Since the set  $\{\phi\in L^2(I;\R^n): \phi=0 \text{ on  } [T_\phi,\infty) \text{ for some } T_\phi>0\} $ is dense in $L^2(I;\R^n)$, the claim follows.
\end{proof}

A result analogous to Proposition \ref{prop:nemitskylip} also  holds for superposition operators on~$\Linfboch$.
\begin{lemma} \label{lem:Nemitsky}
Let~$F\in \operatorname{Lip}\left( \bar{B}_M(0);\R^{m } \right)$, with~$F(0)=0$,
and~$\mathbf{y}\in \Linfboch$, with  $\|{\mathbf{y}}\|_{L^\infty_\mu(Y_0; L^2(I;\R^n)}\leq M$ be given. Define the induced Nemitsky operator by
\begin{align} \label{def:supf}
\lbrack \F(\mathbf{y}) \rbrack(y_0)(t)= F(\mathbf{y}(y_0)(t)) \quad \text{for}~\mu-a.e.~y_0\in Y_0,~t \in I.
\end{align}
Then~$\F(\mathbf{y})\in L^\infty_\mu(Y_0;L^2(I;\R^m))$ and
\begin{align*}
\|\F(\mathbf{y}_1)-\F(\mathbf{y}_2)\|_{L^\infty_\mu(Y_0; L^2(I;\R^m))} &\leq L^1_{F,M} \|\mathbf{y}_1-\mathbf{y}_2\|_{L^\infty_\mu(Y_0; L^2(I;\R^n))} \\
\|\F(\mathbf{y}_1)-\F(\mathbf{y}_2)\|_{L^2_\mu(Y_0; L^2(I;\R^m))} &\leq L^2_{F,M} \|\mathbf{y}_1-\mathbf{y}_2\|_{L^2_\mu(Y_0; L^2(I;\R^n))}
\end{align*}
for all~$\mathbf{y}_i \in \Linfboch$,~$\|{\mathbf{y}_i}\|_{L^\infty_\mu(Y_0;\,\mathcal{C}_b(I; \R^m))}\leq M$,~$i=1,2$.

If~$F$ is continuously Fr\'echet differentiable on~$\bar{B}_M(0)$ and~$DF\in \operatorname{Lip}(\bar{B}_M(0);\R^{m \times n})$, then~$\F: L^\infty_\mu(Y_0;W_\infty))\to  L^\infty_\mu(Y_0;L^2(I;\R^m))$  is differentiable at~$\mathbf{y}\in \Linfboch$ with~$\|{\mathbf{y}_i}\|_{L^\infty_\mu(Y_0;\,\mathcal{C}_b(I; \R^m))} <M$. Its Fr\'echet derivative satisfies
\begin{align*}
D\F(\mathbf{y})\in \mathcal{L}(\Linfboch,L^\infty_\mu(Y_0;L^2(I;\R^m))),\quad\lbrack D\F(\mathbf{y})\delta \mathbf{y} \rbrack (t)= DF(\mathbf{y}(t))\delta \mathbf{y}(t)
\end{align*}
for all~$\delta \mathbf{y} \in \Linfboch$, almost every~$t\in I$, and~$\mu$-a.e.~$y_0\in Y_0$. Moreover we have
\begin{align*}
\|(D\F(\mathbf{y}_1)&-D\F(\mathbf{y}_2))\delta \mathbf{y}\|_{L^\infty_\mu(Y_0;L^2(I;\R^m))}\\& \leq L_{DF,M} \|\imath\|_{\mathcal{L}(W_\infty, \mathcal{C}_b(I;\R^n))} \|\mathbf{y}_1-\mathbf{y}_2\|_{L^\infty_\mu(Y_0;L^2(I;\R^n))}\|\delta \mathbf{y}\|_{L^\infty_\mu(Y_0;W_\infty)}
\end{align*}
for all~$\mathbf{y}_i \in \Linfboch$,~$\|{\mathbf{y}_i}\|_{L^\infty_\mu(Y_0;\,\mathcal{C}_b(I; \R^m))}< M$,~$i=1,2$, and~$\delta \mathbf{y}\in \Linfboch$.
\end{lemma}
Next, we discuss the relation between weak and pointwise almost everywhere limits in~$\Lzwoboch$.
\begin{lemma} \label{lem:measoflimit}
Let~$\{\mathbf{y}_k\}_{k \in \N} \subset \Linfboch$ be given. Assume that there exists a constant~$M_0 >0$, an element~$\mathbf{y} \in \Lzwoboch$ such that:
\begin{itemize}
\item[(i)]
$
\Linfbochnorm{\mathbf{y}_k} \leq M_0 \quad \forall k \in \N,
$
\item[(ii)]$
\mathbf{y}_k \rightharpoonup \mathbf{y} \quad \text{in}~\Lzwoboch,$
\item[(iii)] there exists a family~$\{\widehat{\mathbf{y}}(y_0)\}_{y_0 \in Y_0}\subset W_\infty$ such that
$\mathbf{y}_k (y_0) \rightharpoonup \widehat{\mathbf{y}}(y_0) \quad \text{in}~W_\infty$ for~$\mu$-a.e.~$y_0 \in Y_0$.
\end{itemize}
Then
$ y_0 \mapsto \widehat{\mathbf{y}}(y_0)$
belongs to~$\Linfboch, \,\Linfbochnorm{\widehat{\mathbf{y}}}\leq M_0$,  and~$\widehat{\mathbf{y}}=\mathbf{y}$~$\mu$-almost everywhere.

If, moreover, $F\in \operatorname{Lip}\left( \bar{B}_M(0);\R^{m } \right)$, with~$F(0)=0$ and $M=M_0  \|\imath\|_{\mathcal{L}(W_\infty, \mathcal{C}_b(I;\R^n))} $, then
\begin{align}\label{eq:kk2}
\mathcal{F}(\mathbf{y}_k) \rightharpoonup \mathcal{F}(\mathbf{y}) \quad \text{in}~L^2_\mu(Y_0;L^2(I;\R^n)).
\end{align}
\end{lemma}
\begin{proof}
For arbitrary ~$\varphi \in W^*_\infty$ define
\begin{align*}
\widehat{\mathbf{y}}^\varphi \colon Y_0 \to \R, \quad y_0 \mapsto \langle \widehat{\mathbf{y}}(y_0), \varphi \rangle_{W_\infty, W^*_\infty}
\end{align*}
as well as
\begin{align*}
\mathbf{y}^\varphi_k \colon Y_0 \to \R, \quad y_0 \mapsto \langle \mathbf{y}_k(y_0), \varphi\rangle_{W_\infty, W^*_\infty}
\end{align*}
for~$k \in \N$. Clearly, the mapping~$\mathbf{y}^\varphi_k$ is~$\mu$-measurable for every~$k \in \N$. By assumption (iii) there exist $\mu$-zero sets~$O,O_k \in \mathcal{A}$,~$k \in \N$ such that
\begin{align*}
\mathbf{y}_k(y_0) \rightharpoonup \widehat{\mathbf{y}}(y_0)  \quad \forall y_0 \in Y_0 \setminus O
\end{align*}
as well as
\begin{align*}
\wnorm{\mathbf{y}_k(y_0)} \leq M_0 \quad \forall y_0 \in Y_0 \setminus\bigcup_{k \in \N} O_k.
\end{align*}
Setting~$\mathcal{O}=O \cup \bigcup_{k \in \N} O_k $, we have~$\mu(\mathcal{O})=0$ as well as
\begin{align*}
\langle \mathbf{y}_k(y_0), \varphi\rangle_{W_\infty, W^*_\infty} \rightarrow \langle \mathbf{y}_k(y_0), \varphi\rangle_{W_\infty, W^*_\infty} \quad \forall y_0 \in Y_0 \setminus \mathcal{O}
\end{align*}
and by (i)
\begin{align} \label{def:boundaux}
\wnorm{\widehat{\mathbf{y}}(y_0)} \leq \liminf_{k \rightarrow \infty} \wnorm{\mathbf{y}_k(y_0)} \leq M_0 \quad \forall y_0 \in Y_0 \setminus \mathcal{O}.
\end{align}
Since~$\mathcal{O}$ is a~$\mu$-zero set, the mapping~$\widehat{\mathbf{y}}^\varphi$ is the pointwise almost everywhere limit of a sequence of $\mu$-measurable functions. Thus it is~$\mu$-measurable. Since~$\varphi \in W^*_\infty$ was chosen arbitrarily, this implies the weak measurability of~$\widehat{\mathbf{y}}$. Due to the separability of~$W_\infty$ and Pettis' theorem we conclude its $\mu$-measurability. From~\eqref{def:boundaux} we further get~$\widehat{\mathbf{y}} \in \Linfboch$ and
$
\Linfbochnorm{\widehat{\mathbf{y}}}\leq M_0.
$

Next we verify that  the weak limit $\mathbf{y}$  of~$\{\mathbf{y}_k\}_{k \in \N}$ coincides $\mu$-a.e. with $\widehat{\mathbf{y}}$ . Let~$\Phi \in L^2_\mu(Y_0, W^*_\infty)$ be given. Set
\begin{align*}
\mathcal{O}_\Phi=O \cup \left \{\, y_0 \in Y_0 \;|\;\Phi(y_0)\not\in W_\infty\, \right\}.
\end{align*}
Noting that~$\mu(\mathcal{O}_\Phi)=0$, the mappings
\begin{align*}
\mathbf{y}^\Phi_k \colon Y_0 \to \R, \quad y_0 \mapsto \langle \mathbf{y}_k(y_0), \Phi(y_0) \rangle_{W_\infty,W^*_\infty},
\end{align*}
$k\in\N$, as well as
\begin{align*}
\widehat{\mathbf{y}}^\Phi \colon Y_0 \to \R, \quad y_0 \mapsto \langle \widehat{\mathbf{y}}(y_0), \Phi(y_0) \rangle_{W_\infty,W^*_\infty}
\end{align*}
are~$\mu$-measurable and integrable.
We immediately conclude~$|\mathbf{y}^\Phi_k(y_0)| \leq M_0 \|\Phi(y_0)\|_{W^*_\infty}  $ for all~$k \in \N$ and~$y_0 \in Y_0\setminus \mathcal{O}_\Phi$. Furthermore we have
\begin{align*}
\mathbf{y}^\Phi_k(y_0) \rightarrow \widehat{\mathbf{y}}^\Phi(y_0) \quad \forall~y_0 \in Y_0 \setminus \mathcal{O}_\Phi.
\end{align*}
Denote by~$\langle \cdot, \cdot \rangle$ the duality pairing between~$L^2_\mu(Y_0; W_\infty)$ and~$L^2_\mu(Y_0; W^*_\infty)$. Since~$\mu$ is finite we may now apply Lebesgue's dominated convergence theorem to conclude
\begin{align*}
\langle \mathbf{y}_k, \Phi \rangle= \int_{Y_0} \langle \mathbf{y}_k(y_0), \Phi(y_0) \rangle_{W_\infty, W^*_\infty}~\de \mu(y_0) \rightarrow \int_{Y_0} \langle \widehat{\mathbf{y}}(y_0), \Phi(y_0) \rangle_{W_\infty, W^*_\infty}~\de \mu(y_0)=\langle \widehat{\mathbf{y}}, \Phi \rangle.
\end{align*}
Due to the arbitrary choice of the test function~$\Phi \in L^2(Y_0; W^*_\infty)$ we thus get~$\mathbf{y}_k \rightharpoonup \mathbf{\widehat{y}}$ in~$L^2_\mu(Y_0;W_\infty)$.
Since weak limits are unique there holds~$\widehat{\mathbf{y}}= \mathbf{y}$~$\mu$-a.e, with $\mathbf{y}$ given in (ii).

To verify \eqref{eq:kk2},  let an arbitrary~$\Psi \in L^2_\mu(Y_0;L^2(I;\R^n))$ be given and set
\begin{align*}
\mathcal{O}_\Psi=O \cup \left \{\, y_0 \in Y_0 \;|\;\Psi(y_0)\not\in L^2(I;\R^n)\, \right\}.
\end{align*}
By construction and Corollary \ref{cor21a} there holds~$\mu(\mathcal{O}_\Psi)=0$,
\begin{align*}
( \Psi(y_0), F(\mathbf{y}_k(y_0)))_{L^2(I;\R^n)} \rightarrow \langle \Psi(y_0), F(\mathbf{y}(y_0)) \rangle_{L^2(I;\R^n)} \quad \forall y_0 \in Y_0 \setminus \mathcal{O}_\Phi
\end{align*}
as well as
\begin{align*}
(\Phi(y_0), F(\mathbf{y}_k(y_0)) )_{L^2(I;\R^n)}  \leq M L_{F,M}\|\Phi(y_0)\|_{L^2(I;\R^n)}\quad \forall y_0 \in Y_0 \setminus \mathcal{O}_\Phi,~k\in \N.
\end{align*}
Since~$\mu$ is finite, the right hand side in this estimate is~$\mu$-integrable independent of~$k\in \N$.
Lebesgue's dominated convergence theorem thus yields
\begin{align*}
 \int_{Y_0} (\Phi(y_0), F(\mathbf{y}_k(y_0) )_{L^2(I;\R^n)} ~\de \mu(y_0)
\rightarrow
\int_{Y_0} ( \Phi(y_0), F(\mathbf{y}(y_0)) )_{L^2(I;\R^n)} ~\de \mu(y_0).
\end{align*}
Due to the arbitrary choice of~$\Phi \in L^2_\mu(Y_0;L^2(I;\R^n)) $ we conclude
\begin{align*}
\mathcal{F}(\mathbf{y}_k) \rightharpoonup \mathcal{F}(\mathbf{y}) \quad \text{in}~L^2_\mu(Y_0;L^2(I;\R^n)),
\end{align*}
as desired.
\end{proof}
\section{Perturbation results} \label{app:perturbation}
\subsection{A perturbation result for the nonlinear closed-loop equation}

In this section we study the behavior of solutions to~\eqref{eq:cloloop} under additive perturbations of the dynamical system. In more detail we consider
\begin{align}
\label{eq:pertstate}
\dot{y_v}= \mathbf{f}(y_v)+ \mathcal{B} \mathcal{F}^*(y_v)+v, \quad
y_v(0)=y_0
\end{align}
where~$v \in L^2(I ; \R^n)$ is a given function.

\begin{theorem} \label{thm:existspert}
Let Assumptions~~\ref{ass:feedbacklaw}
and~\ref{ass:controlloflinearized} hold.
Then there exist an open neighbourhood~$V_1 \subset L^2(I; \R^n)$ of~$0$
as well as an open neighbourhood~$\mathbf{Y}_0$ of $Y_0$ such
that~\eqref{eq:pertstate} admits a unique solution~$y_v=\mathbf{y}^v
(y_0)\in \mathcal{Y}_{ad} $ for every pair~$(v,y_0)\in V_1 \times
\mathbf{Y}_0$. Moreover the mapping
\begin{align}
\mathbf{y}^\bullet(\bullet) \colon V_1  \times \mathbf{Y}_0 \to
\mathcal{Y}_{ad}   , \quad (v, y_0) \mapsto \mathbf{y}^v(y_0)
\end{align}
is at least continuously Fre\'chet differentiable.
\end{theorem}
\begin{proof}
The proof is based on the application of the implicit function theorem to
\begin{align*}
G \colon \mathcal{Y}_{ad} \times \mathcal{N}(Y_0) \times L^2(I; \R^n)
\to L^2(I; \R^n) \times \R^n
\end{align*}
with
\begin{align*}
G(y,y_0,v)= \left(
\begin{array}{c}
\dot{y}-\mathbf{f}(y)-\mathcal{B}\mathcal{F}^*(y)-v\\
y(0)-y_0\\
\end{array}
\right)
\end{align*}
Given an arbitrary~$\bar{y}_0 \in Y_0$ and the associated unique
solution~$\bar{y}=\mathbf{y}^*(\bar y_0)\in \operatorname{int}
\mathcal{Y}_{ad}$ (according to Assumption \ref{ass:feedbacklaw}~$\mathbf{A.3}$) to the unperturbed closed loop system,~$G(\bar{y},\bar{y}_0,0)=0$ holds. Moreover~$G$ is at least of
class~$\mathcal{C}^1$ in a neighborhood of~$(\bar{y},y_0,0)$ and there holds
\begin{align*}
D_y G(y,y_0,v)\delta y = \left(
\begin{array}{c}
\dot{\delta y}-D\mathbf{f}(y)\delta y-\mathcal{B} D
\mathcal{F}^*(y)\delta y\\
\delta y(0)\\
\end{array}
\right).
\end{align*}
Assumption~\ref{ass:controlloflinearized} now ensures
that~$D_yG(\bar{y},\bar{y}_0,0)$ is boundedly invertible. Hence applying
the implicit function theorem yields the existence of positive
constants~$\kappa_1=\kappa_1(\bar{y},\bar{y}_0)$
and~$\kappa_2=\kappa_2(\bar{y},\bar{y}_0)$, which may depend
on~$\bar{y},~\bar{y}_0$, such that for every~$y_0\in \R^n$
with~$|y_0-\bar{y}_0|<\kappa_1$ and~$|v|< \kappa_2$ there
exists~$\mathbf{y}^v(y_0) \in \mathcal{Y}_{ad}$
with~$G(\mathbf{y}^v(y_0),y_0,v)=0$ i.e.~$\mathbf{y}^v(y_0)$ is the
unique solution to~\eqref{eq:pertstate} in~$\mathcal{Y}_{ad}$. Moreover, the mapping
\begin{align*}
\mathbf{y}^\cdot(\cdot) \colon  B_{\kappa_2}(0) \times
B_{\kappa_1}(\bar{y}_0) \to \Ya, \quad (v,y_0) \mapsto
\mathbf{y}^v(y_0)
\end{align*}
is of class~$\mathcal{C}^1$.
This yields an open covering of~$Y_0$ i.e.
\begin{align*}
Y_0 \subset \bigcup_{\bar{y}_0 \in Y_0}
B_{\kappa_1(\bar{y},\bar{y}_0)}(\bar{y}_0).
\end{align*}
Since~${Y}_0$ is compact there exists a finite set of initial
conditions~$\{\bar{y}^i_0\}^N_{i=1}\subset Y_0$, including~$0$, such that
\begin{align*}
Y_0 \subset \mathbf{Y}_0 :=\bigcup^N_{i=1}
B_{\kappa_1(\bar{y}^i,\bar{y}^i_0)}(\bar{y}^i_0).
\end{align*}
Set~$V= \bigcap^N_{i=1} B_{\kappa_2(\bar{y}^i,\bar{y}^i_0)}(0) \subset
L^2(I;\R^n) $. Collecting all previous observations now yields the
existence of a~$\mathcal{C}^1$-mapping
\begin{align*}
\mathbf{y}^\cdot(\cdot) \colon V \times \mathbf{Y}_0 \to
\mathcal{Y}_{ad}, \quad \mathbf{y}^v(y_0)~\text{ uniquely
solves}~\eqref{eq:pertstate}~\text{in}~\mathcal{Y}_{ad}.
\end{align*}
\end{proof}
Denote by~$\delta \mathbf{y}^v(y_0) \in \mathcal{L}(V_1 \times
\mathbf{Y}_0, W_\infty)$ the Fr\'echet derivative
of~$\mathbf{y}^\cdot(\cdot)$ at~$(v,y_0)\in V_1 \times \mathbf{Y}_0 $.
It is evident that~$\delta y= \delta \mathbf{y}^v(y_0)(\delta v, \delta
y_0) \in W_\infty$ fulfills
\begin{align} \label{linearizedperturb}
\dot{\delta y}= D \mathbf{f}(\mathbf{y}^v(y_0))\delta y+ \mathcal{B} D
\mathcal{F}^* (\mathbf{y}^v(y_0))\delta y+ \delta v, \quad \delta
y(0)=\delta y_0
\end{align}
for every~$\delta v \in L^2(I;\R^n)$ and~$\delta y_0 \in \R^n$.

  To establish an a priori estimate for the solution to the perturbed
closed loop system~\eqref{eq:pertstate}  we require the following
auxiliary result.
\begin{coroll} \label{coroll:locallipschitzofstate}
There exists an open neighborhood~$V_2 \subset V_1 \subset L^2(I;\R^n)$
of~$0$ as well as~$c>0$ such that
\begin{align*}
\wnorm{\mathbf{y}^{v_1}(y_0)-\mathbf{y}^{v_2}(y_0)} \leq c
\|v_1-v_2\|_{L^2(0,\infty;\R^n)}
\end{align*}
holds for all~$v_1,v_2 \in V_2$ and~$y_0 \in Y_0$.
\end{coroll}
\begin{proof}
Let~$v_1,~v_2 \in V_1$ be given. By the mean value theorem we obtain
\begin{align*}
\wnorm{\mathbf{y}^{v_1}(y_0)-\mathbf{y}^{v_2}(y_0)} &\leq \sup_{s \in
[0,1]}\|\delta
\mathbf{y}^{v(s)}(y_0)(\cdot,0)\|_{\mathcal{L}(L^2(I;\R^n),W_\infty)}
\|v_1-v_2\|_{L^2(I;\R^n)}\\
& \leq
\sup_{s \in [0,1]} \max_{y_0 \in Y_0}\|\delta
\mathbf{y}^{v(s)}(y_0)(\cdot,0)\|_{\mathcal{L}(L^2(I;\R^n),W_\infty)}
\|v_1-v_2\|_{L^2(I;\R^n)}
,
\end{align*}
where~$v(s)=v_1+s(v_2-v_1)\in V_1$,~$s\in [0,1]$.
Let us now consider the mapping
\begin{align*}
h \colon V \to \R, \quad v \mapsto \max_{y_0 \in Y_0}\|\delta
\mathbf{y}^{v}(y_0)(\cdot,0)\|_{\mathcal{L}(L^2(I;\R^n),W_\infty)}.
\end{align*}
Note that $h(v)< \infty$ for all~$v\in V_1$. We now prove that~$h$ is
continuous at zero. To this end let an arbitrary sequence~$\{v_k\}_{k
\in \N}$ with~$v_k \rightarrow 0$ be given. Since the
mapping~$\mathbf{y}^\cdot(\cdot)$ is~$\mathcal{C}^1$ there exists a
sequence~$\{y^k_0\}_{k\in\N} \subset Y_0$ as well as an
element~$\bar{y}_0\in Y_0$ with
\begin{align*}
h(v_k)=\max_{y_0 \in Y_0}\|\delta
\mathbf{y}^{v_k}(y_0)(\cdot,0)\|_{\mathcal{L}(L^2(I;\R^n),W_\infty)}=
\|\delta
\mathbf{y}^{v_k}(y^k_0)(\cdot,0))\|_{\mathcal{L}(L^2(I;\R^n),W_\infty)}
\end{align*}
and~$h(0)=\|\delta
\mathbf{y}^{0}(\bar{y}_0)(\cdot,0)\|_{\mathcal{L}(L^2(I;\R^n),W_\infty)}$.
By cross-testing we obtain
\begin{align*}
g(v^k,\bar{y}_0):=\|\delta
\mathbf{y}^{v_k}(\bar{y}_0)(\cdot,0)\|_{\mathcal{L}(L^2(I;\R^n),W_\infty)}&-\|\delta
\mathbf{y}^{0}(\bar{y}_0)(\cdot,0)\|_{\mathcal{L}(L^2(I;\R^n),W_\infty)}
\leq h(v^k)-h(0)
\end{align*}
as well as
\begin{align*}
h(v^k)\!-\!h(0) \!\leq\! \|\delta
\mathbf{y}^{v_k}(y^k_0)(\cdot,0))\|_{\mathcal{L}(L^2(I;\R^n),W_\infty)}\!-\!\|\delta
\mathbf{y}^{0}(y^k_0)(\cdot,0)\|_{\mathcal{L}(L^2(I;\R^n),W_\infty)}:=
g(v^k,y^k_0).
\end{align*}
Therefore we may estimate
\begin{align*}
|h(v_k)-h(0)| \leq \max \{g(v_k, \bar{y}_0),g(v_k, y^k_0)\}.
\end{align*}
Due to the continuity of~$\delta \mathbf{y}$ the righthand side of this
inequality converges to~$0$ for~$k \rightarrow \infty$. Thus we
get~$h(v_k)\rightarrow h(0)$. Since the sequence~$\{v_k\}_{k\in\N}$ was
chosen arbitrarily we conclude the sequential continuity of~$h$ at~$0$.
Finally, we note that~$L^2(I ;\R^n)$ is a metric space. Hence,
sequential continuity and continuity in the~$\eps-\delta$ sense are
equivalent. In particular, this implies the existence of~$\kappa >0$ as
well as~$c>0$ such that
\begin{align*}
\sup_{s \in [0,1]} \max_{y_0 \in Y_0}\|\delta
\mathbf{y}^{v(s)}(y_0)(\cdot,0)\|_{\mathcal{L}(L^2(I;\R^n),W_\infty)} \leq c
\end{align*}
for all~$v_1,~v_2 \in V_1$ with~$\|v_i\|_{L^2(I;\R^n)}<
\kappa$,~$i=1,2$. Setting~$V_2:= V_1 \cap B_\kappa(0)$ finishes the proof.
\end{proof}
\begin{theorem} \label{thm:aprioriperturbed}
Let Assumptions~~\ref{ass:feedbacklaw}
and~\ref{ass:controlloflinearized} hold.
There exists a constant~$c>0$ such that
\begin{align*}
\wnorm{\mathbf{y}^v(y_0)} \leq M |y_0|+c \|v\|_{L^2(I;\R^n)} \quad
\forall y_0 \in Y_0,~v \in V_2.
\end{align*}
Here~$M$ denotes the constant from~$\mathbf{A.2}$.
\end{theorem}
\begin{proof}
We first point to~$\mathbf{y}^0(y_0)=\mathbf{y}^*(y_0)$ for all~$y_0 \in
Y_0$. Let~$v \in V_2$ be given. We estimate
\begin{align*}
\wnorm{\mathbf{y}^v(y_0)} &\leq \wnorm{\mathbf{y}^*(y_0)} +
\wnorm{\mathbf{y}^v(y_0)-\mathbf{y}^0(y_0)} \\ &\leq
M |y_0|+ c \|v\|_{L^2(I;\R^n)}
\end{align*}
Here we used Assumption~$\mathbf{A.2}$ as well as
Corollary~\ref{coroll:locallipschitzofstate} in the second inequality.
\end{proof}

\section{Smoothness of the value function} \label{app:smoothnessofvalue}

Here we provide sufficient conditions which imply \textbf{A.4} from
Assumption \ref{ass:feedbacklaw} as well as Assumption
\ref{ass:controlloflinearized} in a neighborhood of the origin.
Throughout this subsection we assume that  \textbf{A.1} from Assumption
\ref{ass:feedbacklaw} holds.

We shall assume that $A=D\mathbf{f}(0)$ is exponentially stabilizable, i.e.
\begin{equation}\label{ass:KKtemp} 
   \text{ there exists }  \hat F \in \R^{n\times m} \text{ such that
$e^{(A+B\hat F)t}$
is exponentially stable on } \R^n.
\end{equation}
Then it follows, see e.g. \cite[Theorem
6.2.7]{CurZ95}, that the algebraic Riccati equation
\begin{align}\label{eq:algebraicRiccati}
A^\top \Pi + \Pi A + I = \frac{1}{\beta} \Pi BB^\top \Pi
\end{align}
has a unique nonnegative  solution $\Pi \in \R^{n\times n}$.

Our first goal will be to show that the value function associated to
\eqref{def:openloopproblem} is smooth if ${f}$ is smooth.
It will be convenient to express \eqref{def:openloopproblem} in the form
\begin{equation} \label{def:openloopproblem2} \tag{$P_{\beta}^{y_0}$}
\begin{cases}
\begin{array}{rl}
  & \min ~\frac{1}{2} \int_{0}^{\infty} |y(t)|^2 \,dt + \frac{\beta}{2}
\int_{0}^{\infty} |u(t)|^2 \,dt \\[1.6ex]
  & \text{subject to } e(y,u) = 0,
\end{array}
\end{cases}
\end{equation}
where $e: W_\infty \times L^2(I,\mathbb{R}^m) \rightarrow
L^2(I,\mathbb{R}^n)\times \mathbb{R}^n$ is given by
$$
e(y,u) = (\dot{y} - {{f}}(y) - Bu, y(0) - y_0).
$$
Note that $e$ is $\mathcal{C}^1$, with $De(y,u): W_\infty \times L^2(I,
\mathbb{R}^m) \to L^2(I,\mathbb{R}^n) \times \mathbb{R}^n$ given by
$$
De(y,u)(z,v) = (\dot{z} - D{\bf{f}}(y)z - Bv, z(0)).
$$

We further  introduce $g:W_\infty \to L^2(I)$ as
$$
{{g}}(y) = {{f}}(y) - Ay.
$$
Here and below, contents permitting,  we shall frequently write $L^2(I)$
in place of $L^2(I;\mathbb{R}^n)$ or $L^2(I;\mathbb{R}^m)$. Moreover
balls in $\mathbb{R}^n$ of radius $\delta$ and centered at the origin
are denoted by $B_\delta$.

\begin{lemma} \label{lem:aux1}
There exists a constant $C > 0$ such that for all $\delta \in (0,1]$,
and for all $y$ and $z \in W_\infty$ with $\|y\|_{W_\infty} \le \delta$
and $ \|z\|_{W_\infty} \le \delta$, we have
$$
\|{{g}}(y) - {{g}}(z)\|_{L^2(I)} \le \delta C \|y - z\|_{W_\infty}.
$$
\end{lemma}

\begin{proof}
Due to the continuous embedding $W_\infty \to \mathcal{C}(I;\mathbb{R}^n)$ there
exists $\tilde{\delta}$ such that $|y(t)| \le
\tilde{\delta}$ for all $t \in I$ and $y \in W_\infty$ with
$\|y\|_{W^\infty} \le \delta$. Let $L_{\tilde{\delta}}$ denote the
Lipschitz constant of $g$ on the ball $B_{\tilde{\delta}}$ in
$\mathbb{R}^n$ and let
us define $g:\mathbb{R}^n \to \mathbb{R}^n$ as $g(z) = f(z) - D f(0)z$.
For all $t \geq 0$ we have,
with $y,z$ as in the statement of the lemma,
\begin{equation*}
\begin{aligned}
& |g(y(t)) - g(z(t))| = |f(y(t)) - f(z(t)) - Df(0)(y(t) -
z(t))| \\
& = |\int_0^1 (Df(sy(t) + (1-s)z(t)) - Df(0)) (y(t) -
z(t))| \\
& \le \frac{1}{2} L_{\tilde{\delta}} |y(t) - z(t)|
(|y(t)| + |z(t)| \le \|\imath\|_{\mathcal{L}(W_\infty, \mathcal{C}_b(I;\R^n))}M_{Y_0}
L_{\tilde{\delta}} \delta |y(t) - z(t)|,
\end{aligned}
\end{equation*}
where $\imath$ is the continuous embedding of $W_\infty$ into $\mathcal{C}(I;\R^n)$.
The claim now follows the integration of the above inequality.
\end{proof}

  The following result can be verified by a fixed point argument. For a
proof in a Hilbert space setting, we refer to \cite{BreKP19b}.

\begin{lemma} \label{lem:aux2}
Assume that the spectrum of ${E} \in \mathbb{R}^{n \times n}$ lies in
the left half plane and let $C$ be as in Lemma \ref{lem:aux1}. Then
there exists a constant $C_E$ such that for all $y_0 \in \mathbb{R}^n$
and all $h \in L^2(I)$ with $|y_0| + \|h\|_{L^2(I)} \le
\frac {1}{4C C_E^2}$, the system
$$
\dot{y} = Ey + g(y) + h, ~y(0) = y_0,
$$
has a unique solution in $W^\infty$ satisfying $\|y\|_{W^\infty} \le
2C_E(|y_0| + \|h\|_{L^2(I)})$.
\end{lemma}

This lemma will be applied with two different choices for E.

\begin{coroll} \label{cor:aux1}
There exists a constant $M_{\hat F} > 0$ such that for all $y_0 \in
\mathbb{R}^n$ and all $h \in  L^2(I)$ with $|y_0|_{\R^n} +\|h\|_{L^2(I)}
\le \frac{1}{4C C_{\hat F}^2}$ there exists a control $u \in L^2(I)$ such
that the system
\begin{equation}\label{eq:KK20}
\dot{y} = {{f}}(y) + Bu + h, ~y(0) = y_0,
\end{equation}
has a unique solution $y \in W_\infty$ satisfying
$$
\|y\|_{W^\infty} \le 2 C_{\hat F}(|y_0| + \|h\|_{L^2(I)})
\text{ and } \|u\|_{L^2(I)} \le 2 C_{\hat F} \|\hat F\|
(|y_0| + \|h\|_{L^2(I)}).
$$
\end{coroll}

\begin{proof}
We recall that ${{f}}(y) = {{g}}(y) + Ay$, and that by
\eqref{ass:KKtemp} there exists $\hat F$ such that $A + B \hat F$ is
exponentially stable. Now we can apply Lemma \ref{lem:aux2} to
$$
\dot{y} = (A + B \hat F)y + {{g}}(y) + h, ~y(0) = y_0,
$$
and, setting $u = -\hat Fy$, we arrive at the conclusion.
\end{proof}

\begin{coroll} \label{cor:aux2}
Let $\lambda > \|A\|$ and $(y_0,h) \in \mathbb{R}^n \times L^2(I)$.
There exists $C_\lambda > 0$ with the property:  if for $u \in L^2(I)$
the system
$$
\dot{y} = {{f}}(y) +Bu + h, ~y(0) = y_0,
$$
has a solution $y \in L^2(I)$ satisfying $|y_0| + \|h + \lambda y +
Bu\|_{L^2(I)} \le  \frac{1}{4C C_{\hat F}^2}$, then $y \in W_\infty$, and
$\|y\|_{W_\infty} \le 2M_\lambda(|y_0| + \|h + \lambda y + Bu\|_{L^2(I)})$
holds.
\end{coroll}

\begin{proof}
For $\lambda > \|A\|$ the matrix $A - \lambda I$ is exponentially
stable. This suggests to consider
$$
\dot{y} = (A - \lambda I)y + {{g}}(y) + \tilde{h},
$$
with $\tilde{h} = h + \lambda y + Bu \in L^2(I;\R^n)$. We can now apply
Lemma \ref{lem:aux2} with $E = A - \lambda I$ to assert the claim.
\end{proof}

\begin{lemma} \label{lem:aux3}
There exists $\delta_1 > 0$ such that for every $y_0 \in B_{\delta_1}$
problem \eqref{def:openloopproblem2} possesses a solution $(\bar{y},
\bar{u})$. Moreover there exists $C_1$, such that $\max(\|\bar
u\|_{L^2(I)}, \|\bar y \|_{W_\infty}) \le C_1 |y_0|$.
\end{lemma}

The proof follows standard calculus of variations arguments using
Corollaries \ref{cor:aux1}, \ref{cor:aux2}, and a-priori estimates which
are implied to hold for  minimizing sequences due to the structure of
the cost functional. The smaller the choice of $\lambda$ for the use of
Corollary \ref{cor:aux2}, the larger $M_\lambda$, and the smaller
$\delta_1$ will be, compare \cite[Lemma 8]{BreKP19b}. To pass to the
limit in the state equation, which is satisfied by the elements of
weakly convergent subsequences of state-control pairs, one uses that
$W^{1,2}(0,T)$ embeds compactly into $\mathcal{C}([0,T])$, for every $T > 0$.

\begin{prop}\label{prop:kk1}
There exists $\delta_2 \in (0, \delta_1]$ such that for all $y_0 \in
B_{\delta_2}$, and for all solutions $(\bar{y},\bar{u})$ of
\eqref{def:openloopproblem2}, there exists a unique $W_\infty$ satisfying
\begin{equation}\label{eq:kk30}
\left\{
\begin{array}l
  -\dot{p} - D{{f}}(\bar{y})^t p = \bar{y}, ~\lim_{t \rightarrow \infty}
\bar{p}(t) = 0 \\[1.6ex]
  \beta \bar{u} + B^*p = 0.
\end{array} \right.
\end{equation}
Moreover there exists a constant $C_2$ such that
\begin{equation}\label{eq:kk31}
\|p\|_{W_\infty} \le C_2|y_0|, \text{ for all } y_0 \in B_{\delta_2}.
\end{equation}
\end{prop}

\begin{proof}
We recall the formulation of problem \eqref{def:openloopproblem2} at the
beginning of the subsection and choose $C_1$  as in Lemma
\ref{lem:aux3}. Then
$$
\sup_{t > 0} \sup_{y_0 \in \overline{{B_{\delta_1}}}}
|\bar{y}(t;y_0)|\le C_I C_1 \delta_1,
$$
where $C_I$ denotes the embedding constant of $W_\infty$ into $C(I)$,
and $\bar{y}(\cdot; y_0)$ denotes a solution to
\eqref{def:openloopproblem2} with initial datum $y_0 \in B_{\delta_1}$.
Let $L$ denote the Lipschitz constant of $Df$ on the ball $B_{C_I C_1
\delta_1}$.
We next argue that $De(\bar{y}, \bar{u})$ is surjective, provided that
$\delta_2$ is sufficiently small. We choose an arbitrary pair $(r,s) \in
L^2(I, \mathbb{R}^n) \times \mathbb{R}^n$, and verify that there exists
$(z,v) \in W_\infty \times L^2(I,\mathbb{R}^m)$ such that
\begin{equation}\label{eq:KK32}
\dot{z} - {{A}}z - (D{{f}}(\bar{y}) - D{{f}}(0))z  - Bv = r, ~z(0) = s.
\end{equation}
By the Lipschitz continuity of $Df$ on $B_{C_I M \delta_1}$ we have
\begin{align*}
& \|Df(\bar{y}) - Df(0)\|_{\mathcal{L}(W_\infty, L^2(I,\mathbb{R}^n))} \le
LC_I \|\bar{y}\|_{W_\infty}.
\int_{0}^{\infty}|Df(y) - Df(0)||\delta y| dt 
\|y\|_{W_\infty} \|\delta y \|_{w_\infty}.
\end{align*}
Thus by Lemma \ref{lem:aux3} there exists $\delta_2$ such that
\begin{equation}\label{eq:KK32a}
\|Df(\bar{y}) - Df(0)\|_{\mathcal{L}(W_\infty, L^2(I,\mathbb{R}^n))} \le
LC_I \delta_2 < \frac{1}{M_{\hat F}},
\end{equation}
with $M_{\hat F}$ from Corollary \ref{cor:aux1}. We search for a
solution to \eqref{eq:KK32} with $v = \hat Fz \in L^2(I,\mathbb{R}^n)$, i.e.
\begin{equation}\label{eq:KK33}
\dot{z} - (A + B\hat F)z - (D{{f}}(\bar{y}) - D{{f}}(0))z = r, ~z(0) = s.
\end{equation}
With \eqref{eq:KK32a} holding we can apply \cite[Lemma 2.5]{BreKP18} to
conclude that \eqref{eq:KK33} admits a unique solution $z \in W_\infty$,
satisfying
\begin{equation}\label{eq:KK34}
\|z \|_{W_\infty} \le \tilde{M}(\|r \|_{L^2(I,\mathbb{R}^n)} + |s|)
\end{equation}
for a constant $\tilde{M}$ independent of $(r,s) \in L^2(I,\mathbb{R}^n)
\times \mathbb{R}^n$, and $y_0 \in B_{\delta_2}$.
The surjectivity of $De(\bar{y},\bar{u})$ implies the existence of a
unique Lagrange multiplier $(p, \mu) \in L^2(I;\mathbb{R}^n) \times
\mathbb{R}^n$ such that for all $(z,v) \in W_\infty \times
L^2(I;\mathbb{R}^m)$
\begin{equation}\label{eq:KK35}
DJ(\bar{y},\bar{u})(z,v) - ((p,\mu), De(\bar{y},\bar{u})(z,v))_{L^2(I)
\times \mathbb{R}^n} = 0.
\end{equation}
Choosing $z=0$ and $v \in L^2(I; \mathbb{R}^n)$ arbitrarily, we obtain
the second equation in \eqref{eq:kk30}:
$$
\beta\bar{u} + B^*p = 0 \text{ in } L^2(I; \mathbb{R}^m).
$$
Setting $v = 0$ we find
\begin{equation}\label{eq:KK36}
(p,\dot z)_{L^2(I)} = (p, D{{f}}(\bar{y}(\cdot)) z)_{L^2(I)} + (\bar{y},
z)_{L^2(I)} \text{ for all } z\in W_\infty.
\end{equation}
We observe that $t \to Df(\bar{y}(t)) \in L^\infty(I;\mathbb{R}^{n\times
n})$ and thus $t \to Df(\bar{y}(t))z(t)$ is in $L^2(I;\mathbb{R}^n)$ for
every $z \in L^2(I;\mathbb{R}^n)$. Moreover $S = \{z \in W_\infty :
{\text{supp}}\, z \subset I\}$ is dense in $L^2(I;\mathbb{R}^n)$. Thus
\eqref{eq:KK36} implies that $p \in W_\infty$, and the first equation in
\eqref{eq:kk30} follows.

Finally we need to derive a bound on $p$ in $L^2(I;\R^n)$. Let $r \in
L^2(I)$ and choose $(z,v)$ such that $De(\bar{y},\bar{u})(z,v) = (r,0)$.
Then, using \eqref{eq:KK35} we obtain
\begin{align*}
& (p,r)_{L^2(I)} = ((p,\mu), (r,0))_{L^2(I) \times \mathbb{R}^n} =
(De(\bar{y},\bar{u})^t (p,\mu), (z,v))_{L^2(I) \times \mathbb{R}^n}\\
& = DJ(\bar{y},\bar{u})(z,v) = (\bar{y},z) + (\bar{u},v) \le
\|\bar{y}\|_{L^2(I)} \|z\|_{L^2(I)} + \|\bar{u}\|_{L^2(I)} \|v\|_{L^2(I)} \\
& \le \hat{C}|y_0| \|r\|_{L^2(I)},
\end{align*}
for a constant $\hat{C}$ independent of $y_0 \in B_{\delta_2}$, and $r
\in L^2(I; \mathbb{R}^n)$. Here we used Lemma \ref{lem:aux3} and
\eqref{eq:KK34}. This implies that
$\|p\|_{L^2(I)} \le \hat{C}|y_0|$, for all  $y_0 \in
B_{\delta_2}$.
Now we use the first equation in \eqref{eq:kk30} and the last assertion
in Lemma \ref{eq:KK34}
to deduce \eqref{eq:kk31}.
\end{proof}

Next we carry out a sensitivity analysis for the optimality system. For
this purpose we introduce
\begin{align*}
& \Phi: X = W_\infty \times L^2(I;\mathbb{R}^m) \times W_\infty
\to Y = \mathbb{R}^n \times L^2(I;\mathbb{R}^n) \times
L^2(I;\mathbb{R}^n) \times L^2(I;\mathbb{R}^m)
\end{align*}
by
\begin{align*}
& \Phi(y,u,p) =
\begin{pmatrix}
y(0) \\
\dot{y} - {{f}}(y) - Bu \\
-\dot{p} - D{{f}}(y)p - y \\
\beta u + B^t p
\end{pmatrix},
\end{align*}
and endow $X$ and $Y$ with the $L^\infty$-product norm.

\begin{lemma} \label{lem:aux4}
Assume that ${{f}} \in \mathcal{C}^{k+1}(W_\infty, L^2(I; \mathbb{R}^n))$ with $k
\ge 1$. Then there exist $\delta_3 < 0, \,\delta_3', \, M >0$, and a
$\mathcal{C}^k$-mapping
\begin{equation*}
y_0 \in B_{\delta_3} \to (\mathcal{Y}(y_0), \mathcal{U}(y_0),
\mathcal{P}(y_0)) \in W_\infty \times L^2(I; \mathbb{R}^m) \times W_\infty
\end{equation*}
such that for each $y_0 \in B_{\delta_3}$ the triple $(\mathcal{Y}(y_0),
\mathcal{U}(y_0), \mathcal{P}(y_0))$ is the unique solution to
\begin{equation}\label{eq:KK37}
\Phi(y,u,p) = col(y_0,0,0,0) \text{ with } \|(y,u,p)\|_X \le \delta_3',
\end{equation}
and
\begin{equation}\label{eq:KK38}
\|(\mathcal{Y}(y_0), \mathcal{U}(y_0), \mathcal{P}(y_0))\|_X \le M
|y_0|.
\end{equation}
\end{lemma}

\begin{proof}
We note that $\Phi(0,0,0) = col(0,0,0,0)$ and that $\Phi \in \mathcal{C}^k$. We
argue that $D \Phi(0,0,0)$ is an isomorphism. For this purpose it
suffices to verify that for each $col(w_1, \ldots, w_4) \in Y$ there
exists a unique $(y,u,p) \in X$ such that
\begin{equation}\label{eq:KK39}
D\Phi(0,0,0)(y,u,p) = col(w_1, \ldots, w_4) ~\Leftrightarrow~~
\begin{matrix}
\begin{aligned}
y(0) &= w_1 \\
\dot{y} - Ay - Bu &= w_2 \\
-\dot{p} - A^tp - y &= w_3 \\
\beta u + B^*p &= w_4.
\end{aligned}
\end{matrix}
\end{equation}
This is the necessary and sufficient optimality system to the
linear-quadratic problem
\begin{equation}\label{eq:KK40}
\begin{cases}
\begin{array}{rl}
& \min ~\frac{1}{2} \int_{0}^{\infty} |y + w_3|^2 ~dt + \frac{\beta}{2}
\int_{0}^{\infty} |u|^2 ~dt - \int_{0}^{\infty} u^t w_4 ~dt\\[1.5ex]
& \text{ subject to } \dot{y} = Ay + Bu + w_2, ~y(0) = w_1.
\end{array}
\end{cases}
\end{equation}
By the stability assumption \eqref{ass:KKtemp} it is standard to show,
compare \cite[Proposition 3.1]{BreKP18} that \eqref{eq:KK40} has a
unique solution, (which is  the unique solution to \eqref{eq:KK39}) and that
$$
\|(y,u,p)\|_X \le \tilde{M} \|(w_1, w_2, w_3, w_4)\|_Y,
$$
for some $\tilde{M}$ independ of $w \in Y$. The inverse function theorem
implies the local existence of the $\mathcal{C}^k$-mapping $(\mathcal{Y},
\mathcal{U}, \mathcal{P})$ and \eqref{eq:KK37}. Possibly after further
reduction of $\delta_3$, estimate \eqref{eq:KK38} follows form the fact
that $\Phi(0,0,0) = 0$, and Lipschitz-continuity of $\Phi$.
\end{proof}

\begin{theorem}\label{theo:valuefuntiondiff}
Assume that \eqref{ass:KKtemp} holds and that ${{f}} \in
\mathcal{C}^{k+1}(W_\infty, L^2(I; \mathbb{R}^n))$ for $k \ge 1$. Then the value
function $\mathcal{V}$ associated to \eqref{def:openloopproblem2} is
$\mathcal{C}^k$ on $B_{\delta_4}$ for $\delta_4 = min (\delta_2, \delta_3)$.
\end{theorem}

\begin{proof}
By Lemma \ref{lem:aux3} and Proposition \ref{prop:kk1} there exists
$\delta_2 > 0$ such that for all $y_0 \in B_{\delta_2}$ problem
\eqref{def:openloopproblem2} admits a solution $(\bar{y}, \bar{u})$ with
associated adjoint $p$ such that $\|(\bar{y},\bar{u},\bar{p})\|_X \le \max
(C_1, C_2) |y_0|, \text{ and } \Phi
(\bar{y},\bar{u},\bar{p}) = col(y_0,0,0,0)$. For $y_0 \in B_{\delta_4}$
we deduce from Lemma \ref{lem:aux4} that $(\bar{y},\bar{y},\bar{p}) =
(\mathcal{Y}(y_0), \mathcal{U}(y_0), \mathcal{P}(y_0))$ is the unique
solution to $\Phi(y,u,p) = col(y_0,0,0,0)$. Hence the mapping $y_0 \to
(\bar{y}(y_0),\bar{u}(y_0))$ is $\mathcal{C}^k$ on $B_{\delta_4}$   and  the value
function $\mathcal{V}$ is $C^k$ on $B_{\delta_4}$.
\end{proof}

Finally we turn to justify   Assumption \ref{ass:controlloflinearized}
in a neighborhood of the origin.

\begin{prop}
Let Assumptions \ref{ass:feedbacklaw} and  as well as the assumptions of
Theorem \ref{theo:valuefuntiondiff} hold with $k=2$. Then there exists
$\rho > 0$ and $C>0$ such that for every $y_0 \in B_\rho$
  the linearized closed loop system
\begin{align}
\label{eq:closedloopstate}
\dot{v}= D \mathbf{f}(\mathbf{y}^*(y_0))
v+\mathcal{B}D\mathcal{F}^*({y}^*(y_0))v +\delta v, \quad v(0)=v_0
\end{align}
with  $v_0 \in \R^n, \delta y\in L^2(I;\R^n)$   admits a solution $v \in
W_\infty$
and $\wnorm{v} \leq C(\|\delta v\|_{L^2(0,\infty;\R^n)}+|\delta y_0|).$
Here $\mathbf{y}^*(y_0)) \in W_\infty$ denotes the unique solution to
\eqref{eq:cloloop}.
\end{prop}
\begin{proof}
For the linearized state equation with linearization at the origin,
the optimal feedback law is provided by $-\frac{1}{\beta} B^\top \Pi$,
with $\Pi$ given in
\eqref{eq:algebraicRiccati}. In particular the linear closed loop
system  is exponentially stable and we have for some $c>0$ that
\begin{equation}\label{aux1}
((A- \frac{1}{\beta} B B^\top \Pi) y,y)_{\R^n} \le -c |y|^2_{\R^n},
\text{ for all } y \in \R^n.
\end{equation}
For this $c>0$  we determine $\rho >0$ such that
\begin{equation}\label{aux2}
\begin{array}l
|D( \,{f}({y}^*(y_0)(t))  +{B}  {F}^* ({y}^*(y_0)(t))\,) - D({f}(0)+{B}
{{F}^*}(0))| \le \frac{c}{2}
\end{array}
\end{equation}
for all $t\ge 0$ and $y_0 \in B_{\rho}$. This choice is possible due to
the a-priori estimate \eqref{eq:apriorioptfeedbackass}. We next use the relationship
between the Riccati operator and the second derivative of the value
function as
$$
D( \mathbf{f}(0)+B {\mathcal{F}^*}(0))=  D( \mathbf{f}(0)-
\frac{1}{\beta}B B^\top \nabla V(0)) = A - \frac{1}{\beta} B B^\top \Pi.
$$
Finally we turn to the estimate the asymptotic behavior of the solutions to
\eqref{eq:closedloopstate}. Taking the inner product of
\eqref{eq:closedloopstate} with $y(t)$ we obtain by \eqref{aux1} and
\eqref{aux2}, for $y_0 \in B_\rho$
\begin{equation*}
\begin{array}l
\frac{1}{2} \frac{d}{dt} |v(t)|^2 \le   (\, D {f}( {y}^*(y_0)(t)) + B D{
F^*}( {y}^*(y_0)(t))\, \,v(t), v(t))+(\delta v(t),v(t)) \\[1.5ex]
\qquad \qquad \; \, \le -\frac{c}{2} |v(t)|^2+(\delta v(t),v(t)) .
\end{array}
\end{equation*}
Hence we have
$\frac{1}{2} \frac{d}{dt} |v(t)|^2 \le  -\frac{c}{4} |v(t)|^2+|\delta
v(t)|^2.$
This implies that $|v(t)|^2 \le e^{-\frac{c t}{2}}|v_0|^2 + 2 \int^t_0
e^{\frac{c(s-t)}{2}}  |\delta v(s)|^2 \,ds $, for all $t\ge 0$, all $y_0
\in B_\rho$, and all $v_0\in R^n, \, \delta v \in L^2(I;\R^n)$, and hence
  $\|v\|^2_{L^2(I;\R^n)}\le \frac{2}{c}(|v_0|^2 + 2 \|\delta
v\|^2_{L^2(I;\R^n)})$.
  From here the desired estimate follows.
\end{proof}
\section{Universal approximation property} \label{app:universal}
In this last section we give the technical proof of Proposition~\ref{C1app}.
\begin{proof}[Proof of Proposition~\ref{C1app}]
In many aspects we can profit from    \cite[Theorem 1]{LLPS93}, and its proof,  and from \cite{P99} adapted  to our purposes.  Here, however, we only require a mild regularity assumption for $\psi$,  and $w$ and $b$ are restricted to subsets rather than allowed to vary
in all of $\R^n$ and $\R$. These situations are also commented on in the cited references but not treated in detail.

{\em Step 1.} We recall that the set of polynomials in $n$-variables is a dense linear subspace of $\mathcal{C}^1(\mathbb{R}^n,\mathbb{R})$, i.e. for every $\varphi \in \mathcal{C}^1(\mathbb{R}^n,\mathbb{R})$ and every compact set $K \subset \mathbb{R}^n$ there exists a sequence $p_n$ of polynomials in $(x_1,\ldots, x_n)$ such that $\lim_{n \rightarrow \infty} \|\varphi - p_n\|_{C^1(K,\mathbb{R})} = 0$, see \cite[pg169]{Tr67}.\\

{\em Step 2.} We introduce
$$
\mathcal{M}(\mathcal{W}) = span \{\varphi(w\cdot \bullet) : \varphi \in \mathcal{C}^1(\mathbb{R},\mathbb{R}), w \in \mathcal{W} \} \subset  \mathcal{C}^1(\mathbb{R}^n,\R).
$$
Following \cite{LP93} we introduce
the set of homogenous polynomials of degree $k$:
$$
H_k^n = \left\{ \sum_{|{\bf{m}}|_1=k} c_{\bf m} {\bf s^m}: c_{\bf m}\in \mathbb{R} \right\},
$$
as well as the set of all homogenous polynomials in $n$ variables,
$$
H^n = \bigcup_{k=0}^{\infty} H_k^n,
$$
where the usual multi-index notation is used with ${\bf{m}} = (m_1,\ldots, m_n) \in \mathbb{Z}^n, |{\bf m}|_1 = \sum_{i=1}^{n}m_i, \text{ and } {\bf s}^m = s_1^{m_1} \ldots s_n^{m_n}.$ In \cite[Proof, Theorem 2.1]{LP93} it is verified that $H_k^n = span \{({\bf d}\cdot {\bf x})^k : {\bf d} \in \mathcal{W}\} \subset \mathcal{M}(\mathcal{W})$ for all $k$. Here it is used that $\mathcal{W}$ contains an open set and thus there does not exist a nontrivial homogenous polynomial vanishing on it.
Thus $\mathcal{M}(\mathcal{W})$ contains $H^n$ and thus all polynomials. By Step 1 we have that $\mathcal{M}(\mathcal{W})$ is dense in $\mathcal{C}^1(\mathbb{R}^n,\mathbb{R})$ in the $\mathcal{C}^1$-norm.\\

{\em Step 3.} Let $\chi$  be in  $\mathcal{C}^\infty_0(\R,\R)$, the space of $\mathcal{C}^\infty(\R,\R)$-functions with compact support. Since $\psi$ is not a polynomial, $\chi$ can be chosen such that $\psi*\chi = \int \psi(\cdot-y)\chi(y)\, dy$ is not a polynomial as well. This follows from Steps 6 and 7 of \cite[proof of Theorem 1]{LLPS93}. Let us set $\tilde \psi=  \psi*\chi$. Then $\tilde \psi \in \mathcal{C}^\infty(\R,\R)$.

We verify in this step that  $\Sigma_{\tilde 1} = span\{\tilde \psi(\lambda x + b) : \lambda \in \Lambda\setminus 0, b \in B\}$ is dense in $\mathcal{C}^1(\mathbb{R},\mathbb{R})$. Here $\Lambda$ is an open neighborhood of the origin with the property that $\lambda \in \Lambda$ implies that $\lambda w \in \mathcal{W}$ for all $ w \in \mathcal{W}$, a property, which will be used in Step 6 below.

Note  that
$$
d^h(x):
  = \frac{1}{2h}(\tilde \psi((\lambda + h) x + b) - \tilde \psi ((\lambda-h) x + b)) \in \Sigma_{\tilde 1}
  $$
for every $\lambda \in \Lambda, b \in B,$ and $|h|$ sufficiently small with $h \neq 0$.
Moreover
  $$\frac{d}{dx}d^h(x) = \frac{1}{2h}((\lambda + h) \tilde \psi'((\lambda + h)x +b) - (\lambda-h) \tilde \psi'((\lambda-h) x +b)),$$
and thus $\lim_{h \rightarrow \infty} d^h(x) = \frac{d}{d\lambda } \tilde \psi (\lambda x +b)$ in the $C^1$-norm on compact subsets of $\mathbb{R}$.
We have $\frac{d}{d\lambda}\tilde \psi(\lambda x + b) \in cl( { \Sigma}_{\tilde 1})$, where $cl({\Sigma}_{\tilde 1})$ denotes the closure of $\Sigma_{\tilde 1}$ with respect to  the $C^1$-topology. Note moreover that $\frac{d}{dx} \frac{d}{d\lambda} \tilde \psi(\lambda x + b) = \tilde \psi'(\lambda x + b) + \lambda x \tilde \psi''(\lambda x + b)$.
In a similar way we argue that
\begin{equation*}
(\frac{d}{d\lambda})^{(k)} \tilde \psi(\lambda x + b) \in cl({ \Sigma}_{\tilde 1})
\end{equation*}
for all $k=1,2,\dots$, $\lambda \in \Lambda$ and $b\in B$.
We calculate further $(\frac{d}{d\lambda})^{(k)} \tilde \psi(\lambda x + b) = x^k \tilde \psi^{(k)}(\lambda x + b)$, and $\frac{d}{dx} (\frac{d}{d\lambda})^{(k)} \tilde \psi(\lambda x + b) = k x^{(k-1)} \tilde \psi^{(k)}(\lambda x + b) + x^k \lambda \tilde \psi^{(k + 1)}(\lambda x + b)$ for $k = 2,\ldots$. In particular we have
\begin{equation}\label{eq:21}
x^k \tilde \psi^{(k)}( b) \in cl(\Sigma_{\tilde 1}), \text{ for all } k=1,2, \dots, b\in B     ,
\end{equation}
a property which will be used in Step 5 below.

Now we utilize that $\tilde \psi$ is not a polynomial and hence there exist $b_0 \in B$ such that $\tilde \psi^{(k)}(b_0) \neq 0$, see  \cite{CS54,D69},\cite[pg.156]{P99}, for all $k$. Consequently we find that
$$
x^k \tilde \psi^{(k)}(b_0)= (\frac{d}{d\lambda})^{(k)}\tilde \psi(\lambda x+b)|_{\lambda=0,b=b_0} \in cl ({\hat \Sigma}_{\tilde 1})
$$
for all $k=0,1,\dots$, where we have set $\hat\Sigma_{\tilde 1} = span\{\tilde \psi(\lambda x + b) : \lambda \in \Lambda, b \in B\}$, which differs from $\Sigma_{\tilde 1}$ only with respect to the element $\lambda=0$.
Thus $cl({\hat \Sigma}_{\tilde 1})$ contains all polynomials. Since they are dense in $C^1(\mathbb{R},\mathbb{R})$, we have that ${\hat \Sigma}_{\tilde 1}
$ is dense in $\mathcal{C}^1(\mathbb{R},\mathbb{R})$ as well. Thus for each $g\in \mathcal{C}^1(\R,\R)$ and each compact set $K\in\R$ we have the following property: For all $\eps>0$ there exist $m\in \mathbb{N}$, and $\lambda_i \in  \Lambda$, $b_i\in B$ such that $\| g- \sum_{i=1}^m \tilde \psi (\lambda_i \cdot + b_i)\|_{\mathcal{C}^1(K,\R)}\le \eps$.

Since $\tilde \psi \in \mathcal{C}^\infty(\R,\R)$ and in particular $\tilde \psi \in \mathcal{C}^1(\R,\R)$, it follow that $\Sigma_{\tilde 1}$ is dense in $\mathcal{C}^1(\mathbb{R},\mathbb{R})$. In fact, if in the above expansion there are terms with $\lambda_i=0$, they can be replaced by nontrivial, sufficiently small $\lambda_i$, such that  $\| g- \sum_{i=1}^m \tilde \psi (\lambda_i \cdot \bullet + b_i)\|_{\mathcal{C}^1(K,\R)}\le 2 \eps$.

As a final note to this step, we point out that  if  $\psi \in \mathcal{C}^\infty$, then  the regularisation by convolution is not necessary, the last estimate holds with $\tilde \psi$ replaced by $\psi$, and we can directly continue the proof at Step 6.

{\em Step 4.} Let us choose $\alpha \in(0, \bar b-\underline b$, and define $B_{-\alpha}= (\underline b +\alpha, \bar b -\alpha)$. Then
 $B_{-\alpha}$ is a nontrivial interval contained in $B_0$.
Further we choose a sequence of mollifiers $\chi_n \in \mathcal{C}^\infty_0(\R,\R)$ with the properties that
\begin{itemize}
\item $\psi*\chi_n \to \psi$ in  $L^p(K,\R)$ for some $p\in[2,\infty)$ and every compact set $K \subset \R$, and
\item with the support of $\chi_n$ contained in $ (-\alpha, \alpha)$.
\end{itemize}
We show in this step that $\psi*\chi_n \in cl(span\{\psi(\cdot +b):b \in B_0\})$, for each $n$.
  For this purpose we verify that for each $n\in \mathbb{N}$, for each compact set $K \subset \R$, and each $\varepsilon > 0$,
 there exist $m \in \mathbb{N}, \{b_i\}_{i = 1}^m \subset B_0$, and $\{\mu_i\}_{i=1}^m \subset \R$ such that
 \begin{equation}\label{eq:10}
\begin{aligned}
&\|{\psi*\chi_n} - \sum_{i = 1}^m \mu_i \psi (\cdot - b_i)\|_{W^{1,\infty}(K;\mathbb{R})} \\
 &\qquad = \|\int_{-\alpha}^\alpha \psi (\cdot - \xi) \chi_n(\xi)d\xi - \sum_{i = 1}^m \mu_i \psi (\cdot - b_i)\|_{W^{1,\infty}(K;\mathbb{R})} \le 3 \varepsilon.
\end{aligned}
\end{equation}
We set $b_i = -\alpha + \frac{2i\alpha}{m}$, for $i = 0,\ldots,m$,
\begin{equation}\label{eq:11}
\Delta_i = [b_{i-1}, b_i], \quad \mu_i = \Delta_i \chi_n (b_i), \text{ for } i = 1,\ldots,m,
\end{equation}
and choose $\delta = \delta(\varepsilon)$ such that
\begin{equation}\label{eq:121}
10 \delta \|\psi'\|_{L^{\infty}(K_\alpha; \mathbb{R})} \|\chi_n\|_{C(\mathbb{R},\mathbb{R})} \le \varepsilon,
\end{equation}
where $K_\alpha=\{ s=s_1+s_2: s_1\in K, s_2\in (-\alpha, \alpha)\}$.
By assumption there exist $r(\delta) \in \mathbb{N}$ intervals $\{I_j\}_{j = 1}^{r(\delta)}$ with $\mu (U) \le \delta$, such that $\psi$ is uniformly continuously differentiable on $K_\alpha \backslash U$. Here $\mu (U)$ denotes the Lebesgue measure of $U = \bigcup_{i = 1}^{r(\delta)} I_j$.
Next we choose $m$ sufficiently large such that
\begin{equation} \label{eq:13}
\frac{\alpha r(\delta)}{m} < \delta,
\end{equation}
and for $\xi_1$ and $\xi_2 \in K_\alpha$ with $|\xi_1 - \xi_2| \le \frac{2\alpha}{m}$
\begin{equation}\label{eq:14}
|\chi_n (\xi_1) - \chi_n(\xi_2)| \le \frac{\epsilon}{2\alpha \| \psi \|_{W^{1,\infty}(K_\alpha;\R)}},
\end{equation}
\begin{equation}\label{eq:15}
| \psi (\xi_1) - \psi(\xi_2)| \le \frac{\epsilon}{2\alpha \| \chi_n \|_{L^\infty(\R,\R)}},
\end{equation}
\begin{equation}\label{eq:16}
| \psi' (\xi_1) - \psi'(\xi_2)| \le \frac{\epsilon}{ \| \chi_n \|_{L^1(\R,\R)}}, \text{ if more over } \xi_1, \xi_2 \in K_\alpha\backslash U.
\end{equation}
To estimate \eqref{eq:10} we first use the triangle inequality
\begin{equation*}
\begin{aligned}
& \| {\psi}*\chi_n - \sum_{i=1}^m \mu_i \psi (\cdot - b_i) \|_{W^{1,\infty}(K; \mathbb{R})} \\
& \le \| \int_{-\alpha}^{\alpha} \psi(\cdot - \xi) \chi_n(\xi) d\xi - \sum_{i=1}^m \int_{\Delta_i} \psi(\cdot - b_i) \chi_n(\xi) d\xi \|_{W^{1,\infty}(K; \mathbb{R})} \\
& + \| \sum_{i=1}^m \int_{\Delta_i} \psi(\cdot - b_i) (\chi_n(\xi) - \chi_n(b_i)) \|_{W^{1,\infty}(K; \mathbb{R})} \\
& \le \| \sum_{i=1}^m \int_{\Delta_i} |\psi(\cdot - \xi) - \psi(\cdot - b_i)|\,|\chi_n(\xi)|d\xi \, \|_{L^\infty(K; \mathbb{R})}
\\
& + \| \sum_{i=1}^m \int_{\Delta_i} |\psi'(\cdot - \xi) - \psi'(\cdot - b_i)|\,|\chi_n(\xi)|d\xi \, \|_{L^\infty(K; \mathbb{R})}
\\
& + \| \sum_{i=1}^m \int_{\Delta_i} (|\psi(\cdot - b_i)| + \psi'(\cdot - b_i)|)|\chi_n(\xi) - \chi_n(b_i)|d\xi \ \|_{L^\infty(K; \mathbb{R})} \\
& \le \frac{1}{m}m\varepsilon + I + \frac{1}{m}m\varepsilon = 2 \varepsilon + I,
\end{aligned}
\end{equation*}
where we used \eqref{eq:14}, \eqref{eq:15}, and $|\Delta_i|=\frac{2\alpha}{m}$. Thus we have
\begin{equation}\label{eq:12}
\| {\psi}*\chi_n - \sum_{i=1}^m \mu_i \psi (\cdot - b_i) \|_{W^{1,\infty}(K; \mathbb{R})} \le 2 \varepsilon +I,
\end{equation}
where $I$ denotes the next to the last equality in the above estimate.  To estimate  $I$ we proceed as follows. For a.e. $x \in K$ we consider the set of intervals characterized by indices $i \in \mathcal{I}$ if and only if $(x - \Delta_i) \cap U = \emptyset$. Then for $i \in \mathcal{I}$ by \eqref{eq:16}
\begin{equation*}
\int_{\Delta_i} |\psi'(x - \xi) - \psi'(x - b_i)||\chi_n(\xi)|d\xi \le \frac{\varepsilon}{\|\chi_n\|_{L^1(\mathbb{R})}} \int_{\Delta_i} |\chi_n(\xi)|d\xi,
\end{equation*}
and thus
\begin{equation*}
\sum_{i \in \mathcal{I}} \int_{\Delta_i} |\psi'(x - \xi) - \psi'(x - b_i)||\chi_n(\xi)|d\xi \le \varepsilon.
\end{equation*}
For the complementary index set $\mathcal{I}^c = \{1,\ldots,m\} \backslash \mathcal{I}$, with the property that $(x - \Delta_i) \cap U \neq \emptyset$ for $i \in \mathcal{I}^C$ we find that the measure of such intervals satisfies $\mu(\bigcup_{i \in \mathcal{I}^C} \Delta_i) \le \delta + \frac{4\alpha}{m}r(\delta) \le 5\delta$, where we used \eqref{eq:13}. Using \eqref{eq:121} we obtain
\begin{equation*}
\sum_{i \in \mathcal{I}^c} \int_{\Delta_i} |\psi'(x - \xi) - \psi'(x - b_i)||\chi_n(\xi)|d\xi \le 10\delta \|\psi'\|_{L^\infty(K_\alpha;\R)} \|\chi_n\|_{L^\infty(\mathbb{R},\R)} \le \varepsilon.
\end{equation*}
Thus $I \le \varepsilon$ and together with  \eqref{eq:12} we obtain \eqref{eq:10}.\\

{\em Step 5.} We establish that  $\Sigma_1 = span\{\psi(\lambda \cdot + b): \lambda \in \Lambda \setminus {0}, b \in B\}$ is dense in $\mathcal{C}^1(\mathbb{R},\mathbb{R})$.  We start by observing the following inclusions, which hold for each $n$:
\begin{equation}\label{eq:20}
\begin{aligned}
x^k (\psi*\chi_n)^{(k)}(b) &\subset cl(span \{\psi*\chi_n(\lambda \cdot + b): \lambda \in \Lambda\setminus 0, b\in B_{-\alpha}\}) \\[1.5ex]
& \subset cl(span \{\psi(\lambda \cdot + b): \lambda \in \Lambda\setminus 0, b\in B_0\}) \\[1.5ex]
& \subset cl(span \{\psi(\lambda \cdot + b): \lambda \in \Lambda\setminus 0, b\in B\}),
\end{aligned}
\end{equation}
where the first inclusion one holds for each $b\in B_{-\alpha}$ and each $k=0,1,\dots.$ The last inclusion in \eqref{eq:20} is obvious.
The second inclusion is a consequence of Step 4. For $k=1,2,\dots$ the first one follows from \eqref{eq:21} in Step 3  with $B$ replaced by $B_{-\alpha}$ and  $\tilde \psi= \psi*\chi_n$. Note that  in Step 3 the requirement  that $\tilde \psi$ is not a polynomial, is only used after \eqref{eq:21}, and hence is applicable for $\tilde \psi=\psi*\chi_n$. For $k=0$ the first inclusion can be achieved by appropriate choice of small  $\lambda\neq 0$.

If $\Sigma_{1} = span\{ \psi(\lambda \cdot + b) : \lambda \in \Lambda, b \in B\}$ is not dense in $\mathcal{C}^1(\R,\R)$, then $x^{ k'}$ is not in $\bar \Sigma_{1}$ for some $k'$. From \eqref{eq:20} we conclude that $(\psi*\chi_n)^{(k')}(b)=0$ for all $b$ in  $B_{-\alpha}$ and all $n$.
This implies that $\psi*\chi_n$ is a polynomial of at most degree $k'-1$ for all $n$. By the choice of the sequence $\chi_n$ in Step 4, this implies that $\psi$ itself is a polynomial of at most degree $k'-1$. This gives a contradiction and hence $\Sigma_{1}$ is  dense in $\mathcal{C}^1(\R,\R)$.

{\em Step 6.} Now we show that $\Sigma_n$ is dense in $C^1(\mathbb{R}^n,\mathbb{R})$. Let $g \in \mathcal{C}^1(\mathbb{R}^n,\mathbb{R})$ be arbitrary and let $K \subset \mathbb{R}^n$ be an arbitrary compact set. By Step 2 for every $\epsilon > 0$ there exists $k = k(\epsilon), \,\varphi_i \in C^1(\mathbb{R})$, and $w_i \in \mathcal{W}, i = 1,\ldots,k$ such that
\begin{equation}\label{1}
\|g(\bullet) - \sum_{i = 1}^{k} \varphi_i(w_i \cdot \bullet)\|_{\mathcal{C}^1(K,\mathbb{R})} < \frac{\epsilon}{2}.
\end{equation}
 Since $\varphi \in \mathcal{C}^1(\R,\R)$  we can choose  $w_i\neq 0$ for each $i$.
Moreover, since $K$ is compact there exist intervals $[\alpha_i, \beta_i], i = 1,\ldots,k$ such that $\{w_i \cdot x : x \in K\} \subset [\alpha_i, \beta_i]$.

By Step 5 there exist indices $m_i$, and constants $c_{ij}, \lambda_{ij}\neq 0$,
and $b_{ij}$, with $i = 1,\ldots,k, j = 1,\ldots,m_i$, and $(\lambda_{ij},b_i)\in \Lambda \times B$  such that
\begin{equation}\label{eq3.9}
\|\varphi_i - \sum_{j=1}^{m_i} c_{ij} \psi
(\lambda_{ij}(\cdot) + b_{ij})\|_{W^{1,\infty}(a_i,b_i;\mathbb{R})} \le \frac{1}{\max(1, |w_i|_2)} \frac{\epsilon}{2k}.
\end{equation}
We next consider \eqref{eq3.9} with $\psi \in W^{1,\infty}(\R,\R)$ replaced by a representative $\hat \psi\in \psi$, where now $\hat \psi ' \in \mathcal{L}^\infty(\R,\R)$. By \eqref{eq3.9}
there exists a set $S_i\subset \R$ with $\mu(S_i)=0$ such that
\begin{equation}\label{eq3.10}
\begin{aligned}
\sup_{x\in [a_i,b_i] \setminus S_i} &|\varphi_i(x) - \sum_{j=1}^{m_i} c_{ij} \psi(\lambda_{ij}(x) + b_{ij})|\\[1.7ex] &  \quad + |\frac{d}{dx}(\varphi_i(x) - \sum_{j=1}^{m_i} c_{ij} \psi(\lambda_{ij}(x) + b_{ij}))|  \le \frac{1}{\max(1, |w_i|_2)} \frac{\epsilon}{2k}.
\end{aligned}
\end{equation}
We set $\mathcal{S}_i= \{x\in\R^n| w_i\cdot x \in S_i \}$. Since $w_i\neq 0$ it follows that $\mu(\mathcal{S}_i)=0$, see e.g. \cite{P87}. We therefore have
\begin{equation*}
\begin{aligned}
&\sup_{x\in K\setminus \mathcal{S}_i} |\varphi_i(w_i \cdot x) \!- \! \sum_{j=1}^{m_i} c_{ij} \psi (\lambda_{ij} w_i \cdot x + b_{ij})| \!+\! |\nabla_x(\varphi_i(w_i \cdot x) - \sum_{j=1}^{m_i}c_{ij} \psi (\lambda_{ij} w_i \cdot x + b_{ij})|_2 \\
& \le \sup_{y \in [\alpha_i,\beta_i]} \Big(|\varphi_i(y)\! -\! \sum_{j=1}^{m_i} c_{ij} \psi (\lambda_{ij} y\! +\! b_{ij})| + |w_i|_2 ~|\varphi_i'(y)\! -\! \sum_{j=1}^{m_i} c_{ij} \lambda_{ij} \psi' (\lambda_{ij} y + b_{ij})| \Big)\le \frac{\epsilon}{2 k},
\end{aligned}
\end{equation*}
where we used \eqref{eq3.10}.  This estimate together with \eqref{1} imply
\begin{equation}\label{eq3.11}
\sup_{x\in K\setminus \bigcup_{i=1}^n\mathcal{S}_i}| g(x) - \sum_{i = 1}^{k} \sum_{j = 1}^{m_i} c_{ij} \psi (\lambda_{ij} w_i \cdot x + b_{ij})| < \epsilon.
\end{equation}
From our choice of $\Lambda$  we conclude that $\lambda_{ij} w_i \in \mathcal{W}$ for all $(i,j)$.
Since $ \mu( \bigcup_{i=1}^n) =0$ and     $\hat \psi \in \psi$ was arbitrary, inequality \eqref{eq3.11} implies the desired result.
\end{proof}

\section{Semiglobal regularity of value function} \label{app:semiglobal}

This appendix corresponds to Remark \ref{rem:semiglobal}. 
 Let us consider the controlled stabilization problem
\begin{equation}\label{eq.D15}\tag{E.1}
\left\{
\begin{aligned}
& \inf_{(y,u) \in W_{\infty} \times L^2(I; \mathbb{R}^m)} \frac{1}{2} \int_0^\infty (|y(t)|^2 + \beta |u(t)|^2) dt\\
& \text{ s.t. } ~\dot{y} = Ay + g(y) + Bu, ~y(0) = y_0,	
\end{aligned}
\right.
\end{equation}
where $A \in \mathbb{R}^{n \times n}, B\in\mathbb{R}^{n \times m}, g:\mathbb{R}^n \rightarrow \mathbb{R}^n$ with $g(0) = 0$. We assume that
\begin{equation}\label{eq.D16}\tag{E.2}
\left\{
\begin{aligned}
& (A,B) \text{ is stabilizable in the sense that there exists } K \in \mathbb{R}^{m \times n}, \mu > 0 \\
& \text{ such that } \quad ((A + BK)y, y) \le - \mu|y|^2, \text{ for all } y \in \mathbb{R}^n,
\end{aligned}
\right.
\end{equation}
\begin{equation}\label{eq.D17}\tag{E.3}
\left\{
\begin{aligned}
&  g \in C^2(\mathbb{R}^n,\mathbb{R}^n), \text{ with } || Dg(x)|| \le L \text{ for some } L \in (0, \mu) \text{ independent of} \\
& x \in \mathcal{N}(Y_0). \text{ Moreover } ||Dg(x)|| \le \tilde{L} \text{ for some } \tilde{L} \text{ for all } x \in \mathbb{R}^n.
\end{aligned}
\right.
\end{equation}
For the closed loop system we have the following property.
\begin{lemma} \label{lem:D37}
For every $y_0 \in \mathcal{N}(Y_0)$ there exists a unique solution to
\begin{equation}\label{eq.D18}\tag{E.4}
\dot{y} = (A + BK)y + g(y), ~y(0) = y_0,
\end{equation}
satisfying for all $t \ge 0$
\begin{equation}\label{eq.D19}\tag{E.5}
\begin{aligned}
& |y(t)| \le |y_0| ~exp(2(L-\mu)t), \;\|y\|^2_{L^2(I;\mathbb{R}^n)} \le \frac{1}{2(\mu - L)} |y_0|^2, \\
& \|\dot{y}\|^2_{L^2(I;\mathbb{R}^n)} \le \frac{1}{2(\mu - L)} (\|A + BK\|^2 + L^2) |y_0|^2.
\end{aligned}
\end{equation}
\end{lemma}

\begin{proof}
Local existence is obvious. Global existence will follow from the following a-priori estimate. We take the inner product of \eqref{eq.D18} with $y(t)$ to obtain
\begin{equation*}
\frac{1}{2} \frac{d}{dt} |y(t)|^2 = ((A + BK)y(t), y(t)) + (g(y(t), y(t))) \le (-\mu + L) |y(t)|^2,
\end{equation*}
where we used that $g(y(t)) = \int_{0}^1 Dg(s y(t)) g(t)~ds$ and \eqref{eq.D17}.
Gronwall's lemma implies that
\begin{equation*}
|y(t)|^2  \le |y_0|^2 ~exp(2(L-\mu)t) \le |y(0)|^2,
\end{equation*}
for each $t\ge 0$. Moreover we have
\begin{equation*}
\int_0^\infty |y(t)|^2 dt \le |y_0|^2 \int_0^\infty e^{2(L-\mu)t} dt =  \frac{1}{2(\mu - L)t} |y(0)|^2,
\end{equation*}
and
\begin{equation*}
\begin{aligned}
& \int_0^\infty |\dot{y}(t)|^2 dt \le \int_0^\infty |(A + BK)y|^2 + \int_0^\infty |g(y(t)) - g(0)|^2 dt  \\
& \le (\|A + BK\|^2 + L^2)\|y\|^2_{L^2(I; \mathbb{R}^n)} \le \frac{1}{2(\mu - L)}  (\|A + BK\|^2 + L^2)|y_0|^2,
\end{aligned}
\end{equation*}
and the estimates in \eqref{eq.D19} follow.
\end{proof}

\begin{remark}
Assumption \eqref{eq.D16} can be generalized by replacing the Euclidean inner product by any inner product. As a consequence in the subsequent estimates equivalent norm subordinate to the new inner product have to be observed.
\end{remark}

The previous lemma implies that \eqref{eq.D15} admits at least one feasible control given by $u(t) = Ky(t)$.

\begin{lemma}\label{lem:D38}
For all $y_0 \in \mathcal{N}(Y_0)$ problem \eqref{eq.D15} admits at least one optimal solution $(y^*, u^*) \in W_\infty \times L^2(I;\mathbb{R}^n)$, satisfying
\begin{equation}\label{eq.D20}\tag{E.6}
\|y^*\|^2_{W_\infty} + \|u^*\|^2_{L^2(I;\mathbb{R}^m)} \le \tilde{M} |y_0|^2
\end{equation}
for some $\tilde{M} = \tilde{M}(\mu,L)$ independent of $y_0 \in \mathcal{N}(Y_0)$.
\end{lemma}

\begin{proof}
Let $y_0 \in \mathcal{N}(Y_0)$ be arbitrary and let $u_K$ be the feedback control constructed in the previous lemma with associated trajectory $y_K$. For any minimizing sequence $(y_n,u_n)$ to \eqref{eq.D15} we have for all $u$ sufficiently large
\begin{equation}\label{eq.D21}\tag{E.7}
\begin{aligned}
& \frac{1}{2} \int_0^\infty (|y_n(t)|^2 + \beta |u_n(t)|^2) ~dt ~\le ~\frac{1}{2} \int_0^\infty (|y_K(t)|^2 + \beta |u_K(t)|^2)~dt\\
& \le ~\frac{1}{4(\mu - L)} (1 + \beta\|K\|^2) |y_0|^2.
\end{aligned}
\end{equation}
By the last part of assumption \eqref{eq.D17} together with \eqref{eq.D21} we conclude that $\{y_n\}$ is bounded in $W_\infty$. Hence we can pass to the weak subsequential limit in $W_\infty$ in the equations satisfied by $(y_n,u_n)$ and conlude that there exists $(y^*,u^*)$ satisfying $\dot{y}^* = Ay^* + g(y^*) + Bu^*, y^*(0)=y_0$. By the lower semi-continuity of norms we conclude that $(y^*,u^*)$ is a solution to \eqref{eq.D15}. Estimate \eqref{eq.D20} follows from \eqref{eq.D21} and the state equation.
\end{proof}

Let us define
\begin{align*}
& \Phi: W_\infty \times L^2(I;\mathbb{R}^m) \times W_\infty
\to \mathbb{R}^n \times L^2(I;\mathbb{R}^n) \times
L^2(I;\mathbb{R}^n) \times L^2(I;\mathbb{R}^m)
\end{align*}
by
\begin{align*}
& \Phi(y,u,p) =
\begin{pmatrix}
y(0) \\
\dot{y} - (Ay + g(y)) + Bu \\
-\dot{p} - A^tp - Dg(y)^tp - y \\
\beta u + B^t p
\end{pmatrix}.
\end{align*}

\begin{lemma}\label{lem:D39}
For all $y_0 \in \mathcal{N}(Y_0)$ and for each corresponding solution $(y^*,u^*) = (y^*(y_0),u^*(y_0))$ there exists a unique costate $p = p(y_0) \in W_\infty$ such that
\begin{equation*}
\Phi(y^*(y_0),u^*(y_0),p) = \it{col}(y_0,0,0,0),
\end{equation*}
and for a constant $K_1$ independent of $y_0 \in \mathcal{N}(Y_0)$
\begin{equation*}
\|p(y_0)\|_{W_\infty} \le K_1 (\|y^*\|_{L^2(I;\mathbb{R}^n)} + \|u\|_{L^2(I;\mathbb{R}^n)}).
\end{equation*}
\end{lemma}

\begin{proof}
Let $y_0$ and $(y^*,u^*)$ be as in the statement of the previous lemma.
We consider the shifted state $z = y - y^*$ with associated shifted problem, equivalent to \eqref{eq.D15}, given by
\begin{equation}\label{eq.D22}\tag{E.8}
\left\{
\begin{aligned}
& \inf_{(z,u) \in W^0_\infty \times L^2(I;\mathbb{R}^m)} \tilde{J}(z,u) = \frac{1}{2} \int_0^\infty (|z + y^*|^2 + \beta |u|^2) dt \\
& \text{ s.t. } \dot{z} = Az + g(z+y^*) -g(y^*) + B(u-u^*), z(0) = 0,
\end{aligned}
\right.
\end{equation}
where $W^0_\infty = \{z \in W_\infty : z(0) = 0\}$. We introduce $$\tilde{e}:W_\infty \times {L^2(I;\mathbb{R}^m)} \to {L^2(I;\mathbb{R}^n)}$$ by $\tilde{e}(z,u) = \dot{z} - Az - g(z+y^*) + g(y^*) - B(u-u^*)$.
Clearly $(z^*,u^*) = (0,u^*)$ is an optimal solution for \eqref{eq.D22}. We have
\begin{equation*}
\begin{aligned}
& D\tilde{J}(0,u^*)(\xi,v) = \langle y^*, \xi \rangle_{L^2(I;\mathbb{R}^n)} + \beta \langle u^*,v \rangle_{L^2(I;\mathbb{R}^m)} \\
& D\tilde{e}(0,u^*)(\xi,v) = \dot{\xi} - A\xi -Dg(y^*)\xi - Bv.
\end{aligned}
\end{equation*}
We verify next that $D\tilde{e}(0,u^*)$ is surjective. Let $\varphi \in L^2(I;\mathbb{R}^n)$ be arbitrary and consider
\begin{equation}\label{eq.D23}\tag{E.9}
\dot{\xi} - (A + BK)\xi -Dg(y^*)\xi = \varphi, \; \xi(0) = 0.
\end{equation}
Clearly there exists a solution and we next provide a bound. Proceeding as in the proof of Lemma \ref{lem:D37} we find from \eqref{eq.D23}
\begin{equation}\label{eq.D24}\tag{E.10}
\begin{aligned}
& \frac{1}{2} \frac{d}{dt} |\xi(t)|^2 = ((A +BK) \xi(t), \xi(t)) + (Dg(y^*) \xi(t), \xi(t)) + (\varphi(t), \xi(t)) \\
& \le (L - \mu) |\xi(t)|^2 + (\varphi(t), \xi(t)) = \frac{1}{2} (L - \mu) |\xi|^2 + \frac{1}{2(\mu - L)} |\varphi|^2 .
\end{aligned}
\end{equation}
and thus
\begin{equation*}
\begin{aligned}
& \frac{d}{dt} (|\xi(t)|^2 e^{(\mu - L)t})  \le (\mu - L)^{-1} e^{(\mu - L)t} |\varphi(t)|^2, \\
& |\xi(t)|^2 \le (\mu - L)^{-1} \int_0^t e^{(L -\mu)(t - s)} |\varphi(s)|^2 ~ds,
\end{aligned}
\end{equation*}
and by a convolution argument
\begin{equation*}
\|\xi\|^2_{L^2(I;\mathbb{R}^n)}  \le (\mu - L)^{-2} \int_0^\infty e^{(L -\mu)s} ~ds \ \int_0^\infty |\varphi(s)|^2 ~ds \le (\mu - L)^{2} \|\varphi\|^2_{L^2(I;\mathbb{R}^n)} .
\end{equation*}
By \eqref{eq.D23} this implies the existence of a constant $K_2$ independent of $y_0 \in \mathcal{N}(Y_0)$ such that
\begin{equation*}
\|\xi\|_{W_\infty} \le K_2 \|\varphi\|_{L^2(I;\mathbb{R}^n)} .
\end{equation*}
In particular, $D\tilde{e}(0,u^*)$ is surjecitve with $v = K\xi$ satisfying
\begin{equation*}
\|v\|_{L^2(I;\mathbb{R}^m)} \le \|K\|\, K_2 \|\varphi \|_{L^2(I;\mathbb{R}^m)} .
\end{equation*}
Hence there exists $p \in L^2(I;\mathbb{R}^n)$ such that for all $(z,v) \in W_\infty^0 \times L^2(I;\mathbb{R}^m)$
\begin{equation}\label{eq.D25}\tag{E.11}
 D\tilde{J}(0,u^*)(z,v) - (p, D\tilde{e}(0,u^*)(z,v)) = 0.
\end{equation}
This implies that
\begin{equation*}
\begin{aligned}
&  -\dot{p} = A^tp + Dg(y^*)^t p + y^*,\\
& \beta u^* + B^t p = 0.
\end{aligned}
\end{equation*}
In particular this also implies that $p \in W_\infty$ and $\lim_{t \rightarrow \infty} p(t) = 0$.\\
Finally we verify the bound on $p$. Let $\varphi$ and $v$ be as above. Then we have
\begin{equation*}
\begin{aligned}
& (p,\varphi)_{L^2(I;\mathbb{R}^n)} = \langle p,\tilde{e}(0,u^*)(\xi,v)\rangle_{L^2(I;\mathbb{R}^n)} = D\tilde{J}(0,u^*)(\xi,v)\\
& = \langle y^*, \xi \rangle_{L^2(I;\mathbb{R}^n)} + \beta \langle u^*,v \rangle_{L^2(I;\mathbb{R}^m)} \\
& \le K_1 \|\varphi\|_{L^2(I;\mathbb{R}^n)} (\|y^*\|_{L^2(I;\mathbb{R}^n)} + \|u^*\|_{L^2(I;\mathbb{R}^m)}),
\end{aligned}
\end{equation*}
for a constant $K_1$ independent of $\varphi \in L^2(I;\mathbb{R}^n)$ and $y_0 \in \mathcal{N}(Y_0)$. This concludes the proof of the lemma.
\end{proof}

\begin{lemma}\label{lem:D40}
Assume that $g \in C^2(\mathbb{R}^n, \mathbb{R}^n)$. Then there exist three $C^1$-mappings
\begin{equation*}
y_0 \in \mathcal{N}(Y_0) \rightarrow (\mathcal{Y}(y_0),\mathcal{U}(y_0), \mathcal{P}(y_0)) \in W_\infty \times L^2(I;\mathbb{R}^m) \times W_\infty
\end{equation*}
such that for each $y_0 \in \mathcal{N}(Y_0)$ the triplet $(\mathcal{Y}(y_0),\mathcal{U}(y_0), \mathcal{P}(y_0))$ is the unique solution to
\begin{equation}\label{eq.D26}\tag{E.12}
\Phi(y,u,p) = (y_0,0,0,0)
\end{equation}
in $W_\infty \times L^2(I;\mathbb{R}^m) \times W_\infty$.
\end{lemma}

\begin{proof}
Let us introduce $\tilde{\Phi}:\mathbb{R}^n \times W_\infty \times L^2(I;\mathbb{R}^m) \times W_\infty \rightarrow \mathbb{R} \times L^2(I;\mathbb{R}^n) \times L^2(I;\mathbb{R}^n) \times L^2(I;\mathbb{R}^m)$ by
\begin{align*}
& \tilde{\Phi}(y_0,y,u,p) =
\begin{pmatrix}
y(0) - y_0 \\
\dot{y} - (Ay + g(y) + Bu) \\
-\dot{p} - A^tp - Dg(y)^tp - y \\
\beta u + B^t p
\end{pmatrix}.
\end{align*}
Since $y \in  W_\infty \subset C_b(0,\infty; \mathbb{R}^n)$ we have that $\tilde{\Phi}$ is a $C^1$-mapping. By Lemma \ref{lem:D39} there exists  for every $y_0 \in \mathcal{N}(Y_0)$  a triple  $(y^*(y_0),u^*(y_0),p^*(y_0))$ such that
\begin{equation*}
\tilde{\Phi}(y_0, y^*(y_0), u^*(y_0), p^*(y_0)) = 0.
\end{equation*}
We shall utilite the implicit function theorem to show that for every $\bar{y}_0 \in \mathcal(Y_0)$ there exists a neighborhood $\mathcal{N}(\bar{y}_0)$ and a $C^1$-triplet $(\mathcal{Y}_{\bar{y}_0},\mathcal{U}_{\bar{y}_0}, \mathcal{P}_{\bar{y}_0}): \mathcal{N}(\bar{y}_0) \rightarrow W_\infty \times L^2(I;\mathbb{R}^n) \times W_\infty$ given the unique solution to
\begin{equation}\label{eq.D27}\tag{E.13}
\tilde{\Phi}(y_0,\mathcal{Y}_{\bar{y}_0}(y_0),\mathcal{U}_{\bar{y}_0}(y_0), \mathcal{P}_{\bar{y}_0}(y_0)) = 0, \; \forall y_0 \in \mathcal{N}(\bar{y}_0).
\end{equation}
For this purpose it suffices to argue that $D_{(y,u,p)} \tilde{\Phi}(y_0, y^*(y_0), u^*(y_0), p^*(y_0)):W_\infty \times L^2(I;\mathbb{R}^m) \times W_\infty \rightarrow L^2(I;\mathbb{R}^n) \times L^2(I;\mathbb{R}^m) \times L^2(I;\mathbb{R}^n)$ is an isomorphism. There $D_{(y,u,p)} \tilde{\Phi}$ is given by
\begin{align*}
& D_{(y,u,p)} \tilde{\Phi}(y_0, y^*(y_0), u^*(y_0), p^*(y_0))(y,u,p) =
\begin{pmatrix}
y(0) \\
\dot{y} - (Ay + Dg(y^*)y + Bu) \\
-\dot{p} - A^tp - Dg(y^*)^tp - D^2g(y^*)^t(p^*,y) - y \\
\beta u + B^t p
\end{pmatrix},
\end{align*}
which is welldefinied. To investigate its inverse, let $(w_1,w_2,w_3,w_4) \in \mathbb{R}^n \times L^2(I;\mathbb{R}^n) \times L^2(I;\mathbb{R}^n) \times L^2(I;\mathbb{R}^m)$ be arbitrary and consider
\begin{equation}\label{eq.D28}\tag{E.14}
\left\{
\begin{aligned}
\begin{matrix}
y(0) = w_1 \\
\dot{y} - (Ay + Dg(y^*)y + Bu) = w_2 \\
-\dot{p} - A^tp - Dg(y)^tp - D^2g(y^*)^t(p^*,y) - y = w_3  \\
\beta u + B^t p = w_4.
\end{matrix}
\end{aligned}
\right.
\end{equation}
Let $\tilde{y} = \tilde{y}(w_1)$ be the unique solution to $\dot{\tilde{y}} = (A + BK) \tilde{y}, \tilde{y}(0) = w_1$.
By Lemma \ref{lem:D37} we have the existence of a constant $K_3$ independent of $w_1$ such that
\begin{equation}\label{eq.D29}\tag{E.15}
\|\tilde{y}(w_1)\|_{W_\infty} \le K_3 |w_1|,\quad  \|\tilde{u}(w_1)\|_{L^2(I;\mathbb{R}^m)} \le K_3 |w_1|,
\end{equation}
where $\tilde{u}(w_1) = K \tilde{y}(w_1)$. Setting $z = y - \tilde{y}(w_1)$, system \eqref{eq.D28} becomes equivalent to
\begin{equation}\label{eq.D30}\tag{E.16}
\left\{
\begin{aligned}
\begin{matrix}
z(0) = 0 \\
\dot{z} - (Az + Dg(y^*)z + Bu) = w_2 + B \tilde{u}(w_1) - Dg(y^*)\tilde{y}(w_1)\\
-\dot{p} - A^tp - Dg(y^*)^tp - D^2g(y^*)^t(p^*,z) - z = w_3 + \tilde{y} + D^2g(y^*)^t(p^*,\tilde{y}) \\
\beta u + B^t p = w_4.
\end{matrix}
\end{aligned}
\right.
\end{equation}
These are the necessary optimality conditions for the following uniformly convex linear-quadratic optimal control problem:
\begin{equation}\label{eq.D31}\tag{E.17}
\left\{
\begin{aligned}
& \inf_{(z,u) \in W_\infty^0 \times L^2(\mathcal{I};\mathbb{R}^m)} \frac{1}{2} \int_0^\infty (|\mathcal{Z}z|^2 + (g_2,z)_{\mathbb{R}^n} + \beta |u|^2 + (g_3,u)_{\mathbb{R}^m}) ~dt \\
& \text{ s.t. } \dot{z} - (Az + Dg(y^*(y_0))\;)z - B(u) - g_1 = 0,\, z(0) = 0,
\end{aligned}
\right.
\end{equation}
with
\begin{equation*}
\begin{aligned}
& g_1 = w_2 + B \tilde{u}(w_1) - Dg(y^*(y_0))\tilde{y}(w_1)\\
& g_2 = w_3 + \tilde{y} + D^2g(y^*(y_0))^t(p^*(y_0),\tilde{y}) \\
& g_3 = -w_4 \\
& \mathcal{I} = I + D^2g(y^*)^t(p^*, \cdot).
\end{aligned}
\end{equation*}
It can be established with standard techniques that \eqref{eq.D31} has a unique solution $(\bar{y},\bar{u}) \in W_\infty \times L^2(I;\mathbb{R}^m)$ for each $y_0 \in \mathcal{N}(Y_0)$, with associated adjoint $\bar{p}$, and that for some $K_4 > 0$
\begin{equation}\label{eq.D32}\tag{E.18}
\|(\bar{y},\bar{u},\bar{p})\|_{W_\infty \times L^2(I;\mathbb{R}^m) \times W_\infty} \le K_4(|w_1| + \|w_2\|_{W_\infty} + \|w_3\|_{L^2(I;\mathbb{R}^m)} + \|w_4 \|_{W_\infty}.
\end{equation}
Hence $D_{(y,u,p)} \tilde{\Phi}(y_0, y^*(y_0), u^0(y_0), p^*(y_0))$ is an isomorphism.\\
Thus, for each $y_0 \in \mathcal{D}(Y_0)$ there exists a neighborhood $\mathcal{N}(\bar{y_0})$ and mappings $(\mathcal{Y}_{\bar{y}_0},\mathcal{U}_{\bar{y}_0}, \mathcal{P}_{\bar{y}_0}): \mathcal{N}(\bar{y_0}) \rightarrow W_\infty \times L^2(I;\mathbb{R}^m) \times W_\infty$ such that \eqref{eq.D27} holds. Since $Y_0$ is compact there exists a finite subcovering of $\mathcal{N}(Y_0)$ by elements of $ \{\mathcal{N}(\bar{y_0})\}_{\bar{y_0} \in Y_0}$. Thus we have a global $C^1$-mapping $(\mathcal{Y},\mathcal{U}, \mathcal{P}): \mathcal{N}(Y_0) \rightarrow W_\infty \times L^2(I;\mathbb{R}^m) \times W_\infty$ such that $\tilde{\Phi}(y_0, \mathcal{Y}(y_0),\mathcal{U}(y_0),\mathcal{P}(y_0)) = 0, \forall y_0 \in \mathcal{N}(y_0)$ and $(\mathcal{Y}(y_0),\mathcal{U}(y_0),\mathcal{P}(y_0))$ satisfy the optimality system for \eqref{eq.D15}, established in Lemma \ref{lem:D39}. Thus \eqref{eq.D26} is satisfied.
\end{proof}

\begin{coroll}
The value function associated to \eqref{eq.D15} is continuously differentiable on $\mathcal{N}(Y_0)$ and $\nabla V (y_0) = p^*(y_0)(0)$.
\end{coroll}

\begin{proof}
By Lemma \ref{lem:D38}, problem \eqref{eq.D15} admits a solution for every $y_0 \in \mathcal{N}(Y_0)$. It satisfies the first order necessary condition $\Phi(y_0,y^*(y_0), u^*(y_0), p^*(y_0)) = (y_0,0,0,0)$ established in Lemma \ref{lem:D39}. By Lemma \ref{lem:D40} the solution to the necessary conditions is unique. Thus the value function is welldefined and the minimizers  are unique. Differentiability of the value function follows by applying the chain rule to the value function and using Lemma \ref{lem:D40}.
\end{proof}

\bibliographystyle{siam}
\bibliography{referencesnew}
\end{document}